\documentclass[12pt,a4paper]{amsart}
\usepackage{}

\usepackage{amsfonts}
\usepackage{amsmath}
\usepackage{amssymb}
\usepackage{enumerate}
\usepackage{hyperref}
\usepackage{latexsym}
\usepackage{xcolor}
\usepackage{amsthm}
\usepackage{graphicx}
\usepackage{multirow}

\def\a{\alpha}
\def\b{\beta}
\def\l{\langle} \def\r{\rangle}
\def\div{\,\,\big|\,\,}
\def\Ga{\Gamma}
\def\O{\mathbf{O}}
\def\sP{\mathsf{P}}
\def\sN{\mathsf{N}}

\def\Core{\mathbf{core}}\def\mod{{\mathrm{mod}~}}
\def\Rad{\mathbf{rad}}
\newcommand\ZZ{\mathbb{Z}} \newcommand\D{\mathrm{D}}
\newcommand\A{\mathrm{A}} \newcommand\Sy{\mathrm{S}}
\newcommand\Aut{\mathbf{Aut}}  \newcommand\Out{\mathbf{Out}}
 \newcommand\Cos{\mathbf{Cos}} \newcommand\K{\mathsf{K}}
\newcommand\soc{\mathbf{soc}}
\newcommand\C{\mathbf{C}}
\newcommand\N{\mathbf{N}} \newcommand\Z{\mathbf{Z}}

\newcommand\GL{\mathrm{GL}} \newcommand\SL{\mathrm{SL}} \newcommand\PG{\mathrm{PG}} \newcommand\AGL{\mathrm{AGL}} \newcommand\PGL{\mathrm{PGL}}  \newcommand\PSL{\mathrm{L}} \newcommand\ASL{\mathrm{ASL}} 
\newcommand\SU{\mathrm{SU}} \newcommand\GU{\mathrm{GU}} \newcommand\PGU{\mathrm{PGU}} \newcommand\PSU{\mathrm{U}}
\newcommand\GO{\mathrm{GO}} \newcommand\SO{\mathrm{SO}} 
\newcommand\Sp{\mathrm{Sp}} \newcommand\PSp{\mathrm{PSp}}

\newcommand\GammaL{\mathrm{\Gamma L}}

\newcommand\AGammaL{\mathrm{A\Gamma L}}
\newcommand\PGammaL{\mathrm{P\Gamma L}}

\newcommand\POmega{\mathrm{P\Omega}}

\newcommand\Sz{\mathrm{Sz}}  
  
 \newcommand\F{\mathrm{F}} \newcommand\G{\mathrm{G}}
\newcommand\Suz{\mathrm{Suz}}  \newcommand\J{\mathrm{J}}
\newcommand\M{\mathrm{M}}  \newcommand\Co{\mathrm{Co}}

\newtheorem{theorem}{Theorem}[section]
\newtheorem{myth}{Theorem}
\newtheorem{corollary}{Corollary}[section]
\newtheorem{remark}{Remark}[section]
\newtheorem{lemma}{Lemma}[section]

\theoremstyle{definition}

\newtheorem{hypothesis}{Hypothesis}[section]

\begin{document}

\title[Prime-valent graphs]{Symmetric graphs of prime valency \\ with a
transitive simple group}
\thanks{Supported by the National Natural Science Foundation of China
(11861012, 11971248, 11731002) and the Fundamental Research Funds for the Central
Universities.}

\author[J.J. Li]{Jing Jian Li}
\address{J.J. Li\\ College of Mathematics and Information Science\\ Guangxi University\\ Nanning  530004\\ P. R. China}
\email{lijjhx@163.com}

\author[Z.P. Lu]{Zai Ping Lu}
\address{Z.P. Lu\\ Center for Combinatorics\\
LPMC, Nankai University\\
Tianjin 300071\\
P. R. China} \email{lu@nankai.edu.cn}

\begin{abstract} A graph $\Ga=(V,E)$ is called a Cayley graph of some group $T$ if  the automorphism group $\Aut(\Ga)$   contains a  subgroup $T$ which acts on regularly on $V$.  If the subgroup $T$ is normal in $\Aut(\Ga)$ then $\Ga$ is called a normal Cayley graph of  $T$. Let $r$ be an odd prime.
Fang et al. \cite{FMW} proved that, with a finite number of exceptions for finite simple group $T$,
every connected symmetric Cayley graph of $T$ of valency $r$ is normal.
In this paper, employing maximal factorizations of finite almost simple groups, we work out a possible list of those exceptions for $T$.

\vskip 5pt

\end{abstract}
\maketitle
\parskip 2pt
\section{Introduction}\label{sct=Int}
All graphs in this paper are assumed to be finite, simple
and undirected.

Let $\Ga=(V,E)$ be a  graph with vertex set $V$ and edge set $E$. Denote by $\Aut(\Ga)$   the automorphism group of $\Ga$, and let $G\le \Aut(\Ga)$.
The graph $\Ga$ is called
$G$-vertex transitive if
$G$ acts transitively on   $V$, and called a Cayley graph of $G$ if $G$ acts regularly on $V$.
For   $\a\in V$, let $G_\a=\{g\in G\mid \a^g=\a\}$ and $\Ga(\a)=\{\b\in V\mid \{\a,\b\}\in E\}$,
 called the stabilizer and  neighborhood
 of $\a$ in $G$ and in $\Ga$, respectively.
 If $\Ga$ is $G$-vertex transitive and $G_\a$ acts transitively on
$\Ga(\a)$ then $\Ga$ is called  $G$-symmetric.

In the literature, the solution of quite a number of   problems on symmetric graphs have been reduced or partially reduced into
the class of graphs  arising from almost simple groups. For example, the reduction for symmetric graphs of prime valency \cite{Lorimer}, the reduction for $2$-arc-transitive graphs established
in \cite{Praeger-ONan}, the Weiss Conjecture \cite[Conjecture 3.12]{Weiss81a}   for non-bipartite locally primitive graphs \cite{CLP2000}, the normality of Cayley graphs of simple groups
\cite{FMW, FPW}, the existence and classification of edge-primitive graphs \cite{Giu-Li-edge-prim, 2-arc-edge-p}, and so on. Certainly, the class of graphs arising from almost simple groups plays an important role
in the theory of symmetric graphs.
In this paper, we will study  those symmetric graphs  admitting a transitive simple group.

Our main motivation comes from the study of the automorphism groups or the normality of symmetric Cayley graphs of simple  groups. (Recall that a Cayley graph $\Ga$ of a group $T$ is normal \cite{Xu-98} if $T$ is a normal subgroup of  $\Aut(\Ga)$.)
Let $\Ga=(V,E)$ be a Cayley graph of some simple group $T$ of prime valency $r\ge 3$. Assume that $\Ga$ is connected and $G$-symmetric, where $T\le G\le \Aut(\Ga)$. Then the  Weiss Conjecture is true for $\Ga$, that is, the orders of vertex-stabilizers have an upper bound  depending only on the valency $r$. This ensures that
$T$ is normal in $G$ with a finite number of exceptions, see
\cite[Theorem 1.1]{FMW}. What interests us are those exceptions for $T$. As proposed in \cite{FMW}, those  exceptions for $T$ are worth figuring out. This problem has been solved for $r\le 5$ in several papers, refer to \cite{DuFengZhou, FLX, FMW, XFWX2005}.
In a recent paper \cite{YFZhChen}, the exceptions for $T$  are classified under the assumption  that a vertex-stabilizer $G_\a$ is solvable. In this paper, we  are trying to make some progress on this problem.

Let $\Ga=(V,E)$ be a connected $G$-symmetric Cayley graph of a simple group $T$, where $T\le G\le \Aut(\Ga)$. Assume that $\Ga$ has   odd prime valency $r$.
Take an intransitive normal subgroup $K$ of $G$, and
define a graph $\Ga_{G/K}$  with
vertex set $V_K:=\{\a^K\mid \a\in V\}$ and edge set
$E_K:=\{\{\a^K,\b^K\}\mid \{\a,\b\}\in E\}$. Let
$\overline{G}$ be the group induced by $G$ on $V_K$.
Then $\Ga_{G/K}$ is a $\overline{G}$-symmetric graph of valency $r$,
$K$ is semiregular on $V$ and $\overline{G}\cong G/K$, see \cite[Theorem 9]{Lorimer} for example. Note that $T\cap K=1$. Then $\overline{G}$ has a transitive simple subgroup which is isomorphic to $T$.
These observations suggest us to investigate the symmetric  graphs
which have prime valency and admit a transitive simple group.


Let   $G^{\infty}$ be the solvable residue of $G$, i.e., the last term of the derived series of $G$. Let $\Rad(G)$ be the solvable
radical of $G$, that is, the maximal solvable normal subgroup of
$G$.
In Section \ref{sct=general},  the following result is proved.

\begin{myth}\label{mainth-1}
Let $\Ga=(V,E)$ be a connected $G$-symmetric graph of prime valency $r\ge 3$. Assume that $G$ contains a  transitive simple  subgroup $T$.
\begin{itemize}
\item[(1)] If $G$ has a transitive minimal normal subgroup $N$ then
  either $N$ is simple and $T\le N$, or $N\cong \ZZ_2^3$, $r=7$ and $T\cong\PSL_3(2)$;  for the former case,   $\Ga$ is $N$-symmetric if $T\ne N$.
\item[(2)] If  $T$ is  not normal in $G$ then
$G$ has a  characteristic subgroup $N$ such that $\Ga$ is $N$-symmetric, $N/\Rad(N)$ is simple, $|\Rad(N)|$ is indivisible by $r$,    $\Rad(N)T\lneq N=N^\infty$,  and either $\Rad(N)$ is intransitive on $V$ or $G=N\cong \AGL_3(2)$.
\end{itemize}
\end{myth}

Part (2) of Theorem \ref{mainth-1}  provides us a method to determine the normality of symmetric Cayley graphs of prime valency  of simple groups, which also suggests us to consider the symmetric graphs satisfying the following
hypothesis.
\begin{hypothesis}\label{main-hypo}
 $\Ga=(V,E)$ is  a connected graph of prime valency $r\ge 3$ such that
\begin{itemize}
\item[(I)]  $\Ga$ is $G$-symmetric, where $G$ is a nonabelian simple group;
\item[(II)]  $G=TG_\a$ for $\a\in V$, where $T$ is a   simple group with $T\ne G$.
\end{itemize}
\end{hypothesis}
\noindent Note that the simple group $T$ in (II) is nonabelian; otherwise,
$\Ga$ has odd prime order, yielding $r=2$, a contradiction.

Most part of this paper, Sections \ref{sect-sl}-\ref{sect=classical}, is devoted to analyzing the possibility of the triples
$(G,T,G_\a)$ satisfying Hypothesis \ref{main-hypo}. In Section \ref{sect-sl}, we give some limitations on $G_\a$ or $G$.
Then the case where $G_\a$ is solvable is dwelt in Section \ref{sect=sol-stab},
the case where $G$ (or $T$) is an alternating group or a sporadic simple group is dwelt in Section \ref{sect-alt}, and the case where $G$ is a classical group is left in Sections \ref{sect-tech} and \ref{sect=classical}. At last, the following result is proved in Section \ref{sect=cay}.
\begin{myth}\label{mainth-2}
Let $\Ga=(V,E)$ be a connected $X$-symmetric Cayley
graph of a finite simple group $T$
of prime valency $r\ge 7$, where $T\le X\le \Aut(\Ga)$. Assume that $T$ is not normal in $X$ and
$X_\a$ is insolvable, where $\a\in V$. Then $T$ is isomorphic to a proper subgroup of some  composition factor $L$ of $X$, and all possible triples $(L,T,r)$ are listed in Table {\rm\ref{exceptions}}.
\end{myth}
\begin{table}[ht]
\[\footnotesize
\begin{array}{|l|l|c|l|l|}\hline
\mbox{Line}&L& T&r&\mbox{Conditions}\\ \hline
1&\A_n& \A_{n-1} & \ge 7 & r\div n, r^2\nmid n, n\mbox{ not  prime}\\ \hline
2&\A_{r+1}& |T|=r+1  & \ge 7 & \Ga=\K_{r+1}\\ \hline
3&\PSU_{2m}(2^f)& \PSU_{2m-1}(2^f)  & 2^{mf}+1 &m\ge 2, mf=2^e>2 \\ \hline
4&\Omega_{2m+1}(p^f)& \POmega^-_{2m}(p^f)  & {p^{mf}-1\over p^{{m\over d}f}-1} &m\ge 3, mf=d^e, \mbox{odd primes } p,d\\ \hline
5&\POmega^-_{10}(2^f)& \PSU_{5}(2^f)  & 2^{4f}+1 & f=2^e>1 \\ \hline
6&\Omega^+_{2m}(2^f)& \Sp_{2m-2}(2^f)  & {2^{mf}-1\over 2^{mf\over d}-1} & m\ge 4,  mf=d^e, \mbox{prime } d\\ \hline
7&\Omega^+_{2m}(2^f)& \Omega^-_{2m-2}(2^f)  & {2^{mf}-1\over 2^{{m\over d}f}-1} & m\ge 5,  mf=d^e, \mbox{odd prime }d\\ \hline
8&\Omega^+_{2m}(p^f)& \POmega^-_{2m-2}(p^f)  & {p^{mf}-1\over p^{f}-1} &   \mbox{odd prime }m\ge 5,   f=m^e,\\ \hline
9&\POmega^+_{2m}(p^f)& \Omega_{2m-1}(p^f)  & {p^{mf}-1\over p^{{m\over d}f}-1} & m\ge 5,  mf=d^e, \mbox{odd primes } p,d \\ \hline
10&\Sp_{4}(2^f)& \SL_2(2^{2f})  & 2^f+1 &   f=2^e>2\\ \hline

11&\Sp_{6}(2^f)& \G_2(2^{2f}), \SL_2(2^{3f})  & 2^{2f}+1 &  f=2^e\ge2\\ \hline

12&\Sp_{2m}(2^f)& \Omega^-_{2m}(2^f)  & {2^{mf}-1\over 2^{{m\over d}f}-1} &  m\ge 4,  mf=d^e, \mbox{prime } d \\ \hline

13&\PSL_{n}(p^f)& \PSL_{n-1}(p^f)  & {p^{nf}-1\over p^{{n\over d}f}-1} & n> d,  nf=d^e,   \mbox{odd prime } d \\ \hline
14&\PSL_{n}(p^f)& \PSL_{m}(p^{{n\over m}f})  & {p^{(n-1)f}-1\over p^{f}-1} & 1<m<n, \mbox{ and } (*):\, m\div n,  \\
&&&& f=(n-1)^e,  \mbox{odd prime }n-1 \\
 & & \PSp_{m}(p^{{n\over m}f})  &
 & (*), \mbox{even }m\ge 4,  p=2 \mbox{ if }m=n\\
 & & \POmega^-_{m}(p^{{n\over m}f})  &
  &  (*), \mbox{even }{m\over 2}\ge 4,   p=2 \mbox{ if }m=n\\
 & & \G_{2}(p^{{n\over m}f})  &
 &  (*), m=6, p=2\\ \hline
15&\PSU_6(2)& \PSU_5(2)  & 7 &    \\ \hline
 16&\Omega^+_{12}(2)& \Sp_{10}(2)  & 7 &    \\ \hline
 17&\Sp_{12}(2)& \Omega^-_{12}(2)  & 7 &    \\ \hline
\end{array}
\]
{\caption{\small Exceptions for $T$.}\label{exceptions}}
\end{table}

\vskip 20pt

\section{Preliminaries}\label{sct=stab}

In this section, we  collects a few of preliminary results which will be used frequently in the following sections.

\subsection{Primitive divisors}
Let $\Phi_{m}(q)$ be the $m$-th cyclotomic polynomial evaluated at a prime power $q=p^f$. Then
\[q^m-1=\prod_{d\mid m}\Phi_d(q),\,\,q^m+1=\prod_{d\mid 2m,d\nmid m}\Phi_d(q).\]
Noting that
\[\prod_{d\mid mf, d\nmid f}\Phi_d(p)={q^m-1\over q-1}=\prod_{1<d\mid m}\Phi_d(q),\]
the next lemma is easily shown.
\begin{lemma}\label{r}
If ${q^d-1\over q-1}$ is a prime then $d$ is a prime, $(d,q-1)=1$,  $f$ is a power of $d$, and either $d$ is odd or $d=p=2$.
\end{lemma}

Recall that a primitive prime divisor of $q^m-1$  is a prime which divides $q^m-1$ but does not  divide $q^l-1$ for any $0<l<m$. If $r$ is a  primitive prime divisor of $q^m-1$, then
$q$ has order $m$ modulo $r$, and thus $m$ is a divisor of $r-1$; if further $r\div (q^n-1)$ with $n\ge 1$ then $m\div n$.
Let $\Phi_{m}^*(q)$ be the product of all primitive prime divisors of $q^m-1$ (counting multiplicities). Then the following lemma holds.
\begin{lemma}\label{phi}
If $\Phi_{m}^*(q)\ne 1$ then
every  prime divisor of $\Phi_{m}^*(q)$ is  no less than $m+1$, and $(\Phi_{m}^*(q),q^n-1)\ne 1$ yields that $m\div n$.
\end{lemma}
It is easily shown that $\Phi_{m}^*(q)$ is a divisor of
$\Phi_{m\over d}^*(q^d)$, where $d$ is a divisor of $m$. In particular, $\Phi_{mf}^*(p)$ is a divisor of $\Phi_{m}^*(q)$.
Moreover, by Zsigmondy's Theorem,  we have
the following fact.
\begin{lemma}\label{prime-div}
$\Phi_{m}^*(q)=1$ if and only if $(m,q)=(6,2)$ or $m=2$ and $q=p=2^s-1$ for some prime $s\ge 2$.
\end{lemma}

\subsection{Factorizations}
All groups in this paper are assumed to be finite, and our group-theoretic notation mostly consistent with \cite{max-factor}.

Let $G$ be a group.
An expression $G=XY$  is called a factorization of $G$ if $X$ and $Y$ are proper subgroups of $G$, where $X$ and $Y$ are called factors. A factorization $G=XY$ is called
maximal  if  both $X$ and $Y$ are maximal subgroups of  $G$.
The next lemma collects a few of simple facts, most of them follow from Lagrange's Theorem.
\begin{lemma}\label{Order}
Let $G=XY$ be  factorization of $G$. Then the following hold.
\begin{itemize}
\item[(1)] $YX=G=A^gY^h$ for $X\le A\le G$ and $g,h\in G$.
\item[(2)] $|G:X|=|Y:(X\cap Y)|$ and $|G:Y|=|X:(X\cap Y)|$; in particular, $|X|$ and $|Y|$ are divisible by $|G:Y|$ and $|G:X|$, respectively.

\item[(3)] $|G|_p\le |X|_p|Y|_p$, $|X|_p\ge |G:Y|_p$ and $|Y|_p\ge |G:X|_p$, where $p$ is a prime and $|\cdot|_p$ means the largest $p$-power dividing $|\cdot|$.
\item[(4)] There is a maximal factorization $G=AB$ with $X\le A$ and $Y\le B$.
\end{itemize}
\end{lemma}

Let $X$ be a subgroup of a group $G$.
The intersection $\cap_{g\in G}X^g$ is called the core of $X$ in $G$, denoted by $\Core_G(X)$. If $\Core_G(X)=1$, i.e., $X$ contains no normal subgroup of $G$ other than $1$, then $X$ is called a core-free subgroup of $G$.

Given a factorization $G=XY$,
we have a
transitive action of $G$ on   $[G:Y]:=\{Yg\mid g\}$ by right multiplication, where the kernel is $\Core_G(Y)$. Clearly, $X$ acts   transitively on   $[G:Y]$. If $Y$ is core-free in $G$ then $G$ (and hence $X$) acts faithfully on $[G:Y]$, so all orbits of a normal subgroup of $X$ have equal length greater than $1$. This implies the next fact.
\begin{lemma}\label{core-free-f}
Let $G$ be a group and $G=XY$ be a factorization. If $Y$ is core-free in $G$, then
$X\cap Y$ is core-free in both $X$ and $G$.
\end{lemma}

A factorization $G=XY$ is called core-free if   $X$ and $Y$ are core-free in $G$.
Let $G=XY$ be a core-free factorization. Assume that   $X<A\le G$. Then we have a factorization $A=X(A\cap Y)$. By Lemma \ref{core-free-f}, since $G=AY$, we know that $A\cap Y$ is core-free in $A$. In general, $X$ is not necessarily core-free in $A$, and the factorization $G=AY$ is not necessarily core-free too.

For a simple group $S$ and a group $G$, by $S\preccurlyeq G$ or $G\succcurlyeq S$ we means that
$G$ has a composition factor (isomorphic to) $S$.

\begin{lemma}\label{sub-factor}
Let  $G=XB$ be a core-free factorization and $X\le A\le G$.
Assume that $N\unlhd A$, $A/N$ is  almost simple, and $X^\infty\not\le N$. Let $\overline{A}=A/N$, $\overline{X}=XN/N$ and $\overline{Y}=(A\cap B)N/N$.
Then either $\overline{A}=\overline{X}\,\overline{Y}$ is core-free, or  one of the following holds:
\begin{itemize}
\item[(1)]
$\overline{A}^\infty\preccurlyeq X$;

\item[(2)] $|G:B|$ is a divisor of $|\overline{A}:\overline{A}^\infty||N|$.
\end{itemize}
If  $\Core_G(A)=1$ and $N=\Rad(A)$ normalize $A\cap B$
then
$\overline{A}^\infty\not \preccurlyeq A\cap B$.
\end{lemma}
\proof
Assume that either $\overline{X}$ or  $\overline{Y}$ is not core-free in  $\overline{A}$.
If   $\overline{X}$ is not core-free in $\overline{A}$ then part (1) of this lemma follows.
Suppose that $\overline{Y}$ is not core-free in $\overline{A}$.
Then $NA^\infty\le N(A\cap B)$. Thus $|G:B|=|A:(A\cap B)|=|A:N(A\cap B)||N(A\cap B):(A\cap B)|$, and so $|G:B|$ is a divisor of $|A:NA^\infty||N|$. Then
part (2) of this lemma occurs.

Now let  $\Core_G(A)=1$ and $\Rad(A)$ normalize $A\cap B$.
Then $((A\cap B)\Rad(A))^\infty=(A\cap B)^\infty$.
Suppose that $S\preccurlyeq A\cap B$.
Then
$(A\cap B)\Rad(A)\ge A^\infty\Rad(A)$, and so
$A^\infty\le ((A\cap B)\Rad(A))^\infty=(A\cap B)^\infty$.
Thus $A^\infty\le B\cap A$, which contradicts Lemma \ref{core-free-f}. Then the lemma follows.
\qed

\begin{lemma}\label{m-factor}
Let  $G=XY$ be a core-free factorization. Suppose that 
$G$ is almost simple.
Then there exists   a core-free maximal factorization $G_0=A_0B_0$ such that
\begin{itemize}
\item[(1)]  $G\ge G_0=(A_0\cap X)(B_0\cap Y)$, $G^\infty\le G_0$, $X^{\infty}\le A_0$ and $Y^{\infty}\le B_0$;
\item[(2)] $|X:(A_0\cap X)|$ and  $|Y:(B_0\cap Y)|$ are divisors of $|G:G_0|$.
\end{itemize}
\end{lemma}
\proof
If $G^\infty X=G=G^\infty Y$, taking maximal subgroups $A_0$ and $B_0$ of $G$ with $X\le A_0$ and $Y\le B_0$,  then $G=A_0B_0$ is a desired factorization. Thus, with out loss of generality, we assume next that $G\ne G^\infty X$.

Let $G_1=G^\infty X$ .
Then $G=G_1Y$, and $G_1=X(G_1\cap Y)$. Let $Y_1=G_1\cap Y$. It is easy to see that $X$ and $Y_1$ are core-free in $G_1$.  Since $G=G_1Y$,
by Lemma \ref{Order}, $|Y:Y_1|=|G:G_1|$. Moreover, $Y^\infty\le G^\infty=G_1^\infty$, and so  $Y^\infty\le G_1\cap Y=Y_1$, yielding $Y^\infty=Y_1^\infty$.
Noting that $|G_1|<|G|$, by induction,
there are $G_0\le G_1$ and a core-free maximal factorization $G_0=A_0B_0$ such that $G_1^\infty\le G_0$, $X^{\infty}\le A_0$ $Y_1^{\infty}\le B_0$,  $|X:(A_0\cap X)|$ and  $|Y_1:(B_0\cap Y_1)|$ are divisors of $|G_1:G_0|$. Since $Y_1\le Y$,
we have $|Y:(B_0\cap Y)||(B_0\cap Y):(B_0\cap Y_1)|=|Y:(B_0\cap Y_1)|=|Y:Y_1||Y_1:(B_0\cap Y_1)|$, it follows that $|Y:(B_0\cap Y)|$ is a divisor of $|G:G_0|$. Clearly, $G^\infty=G_1^\infty\le G_0$ and $Y^\infty=Y_1^\infty\le B_0$. Then the factorization $G_0=A_0B_0$ is desired. \qed

\subsection{Stabilizers}
Let $\Ga=(V,E)$ be a connected
$G$-symmetric graph  of prime valency $r\ge 3$. Denote by $\soc(G)$  the
socle of $G$,
and by $\pi(G)$ the set of prime divisors of $|G|$.
By
\cite[Lemma 1.1]{CLP2000}, \cite[Theorem 9]{Lorimer} and
\cite[Theorem 1]{Lorimer90},  the following lemma holds.

\begin{lemma}\label{N-p}
$|G_\a|_r=r=\max\pi(G_\a)$.
If  $N\unlhd G$  then either $N$ is semiregular on $V$, or $N_\a$ is transitive on $\Ga(\a)$
and $N$ has at most two orbits on $V$.
\end{lemma}

Take $\{\a,\b\}\in E$. Then the arc-stabilizer $G_{\a\b}$ is a Hall $r'$-subgroup of $G_\a$. Since $\Ga$ is $G$-symmetric,
$(\a,\b)^x=(\b,\a)$ for some $x\in G$.
Clearly, such an  $x$ has even order and lies in the normalizer $\N_{G}(G_{\a\b})$ of $G_{\a\b}$. Thus the element $x$ may be chosen as a $2$-element. Moreover, since $\Ga$ is connected, $\l x, G_\a\r=G$.
\begin{lemma}\label{arc-stab-normalizer}
Let $\{\a,\b\}\in E$. Then $\l x, G_\a\r=G$ for  some $2$-element $x\in \N_{G}(G_{\a\b})$ with $x^2\in G_{\a\b}$.
\end{lemma}

Note that $G_\a$ induces a transitive permutation group say $G_\a^{\Ga(\a)}$ (on
$\Ga(\a)$) of  degree $r$. Then $G_\a^{\Ga(\a)}$ and $r$ are known as in Table \ref{tgp-d-r}, refer to
 \cite[pp. 99]{Dixon-book}. Moreover, if $r={p^{fd}-1\over p^f-1}$ for some prime   $p$ then,  by Lemma \ref{r}, $d$ is a prime and $f$ is a power of $d$.

\begin{table}[ht]
\[\begin{array}{|l|l|l|l|l|l|l|} \hline
\soc(G_\a^{\Ga(\a)}) & \ZZ_r& \A_r& \PSL_d(q)&\PSL_2(11)&\M_{11}&\M_{23}\\  \hline
r&r&r\ge 5&{q^d-1\over q-1}  &11 &11&23   \\ \hline
\end{array}\]
{\caption{Transitive groups of prime degree.}\label{tgp-d-r}}
\end{table}
Let $G_\a^{[1]}$ be the kernel of
the action of $G_\a$ on $\Ga(\a)$. Then $G_\a^{\Ga(\a)}\cong G_\a/G_\a^{[1]}$.
 Take   $\beta \in \Ga(\a)$, and set $G_{\a\b}^{[1]}=G_\a^{[1]}\cap G_\b^{[1]}$. Then
 $G_{\a\b}^{[1]}$ has order a prime power, see \cite{Gardiner-73}.
 For  $r=3$, by   \cite{DDM, Tutte47,Tutte59}, $G_\a$ and $G_{\a\b}^{[1]}$ are listed in
 Table \ref{cubic-stab}.
\begin{table}[ht]
\[\small
\begin{array}{|l|l|l|l|l|l|}\hline
G_\a&\ZZ_3&\Sy_3&\D_{12}&\Sy_4&\Sy_4\times \ZZ_2\\ \hline
G_\a^{[1]}& 1 &1 & \ZZ_2& \ZZ_2^2& \ZZ_2^3\\ \hline
G_{\a\b}^{[1]}&1&1& 1& \ZZ_2& \ZZ_2^2\\ \hline
G_{\a\b}& 1& \ZZ_2&\ZZ_2^2&\D_8&\D_8\times\ZZ_2\\ \hline
\end{array}
\]
{\caption{\small Stabilizers of symmetric cubic graphs.}\label{cubic-stab}}
\end{table}

In general, we have
\begin{align}\label{stab-extension}
\begin{split}
G_{\a\b}/G_{\a\b}^{[1]}\lesssim (G_\a^{\Ga(\a)})_\b\times (G_\b^{\Ga(\b)})_\a,\,\, (G_\a^{[1]})^{\Ga(\b)}\unlhd
(G_\b^{\Ga(\b)})_\a\cong (G_\a^{\Ga(\a)})_\b\\
G_\a^{[1]}=G_{\a\b}^{[1]}.(G_\a^{[1]})^{\Ga(\b)},\,\,
G_\a=(G_{\a\b}^{[1]}.(G_\a^{[1]})^{\Ga(\b)}).G_\a^{\Ga(\a)}.
\end{split}
\end{align}
Then, combining with \cite{Gardiner74,Weiss81a}, we have the following result.
\begin{lemma}
\label{stab-lem}
\begin{itemize}
\item[(1)] $G_\a$ is soluble if and only if
$G_\a^{\Ga(\a)}\le \AGL_1(r)$.
\item[(2)] Either $\PSL_d(q)\unlhd G_\a^{\Ga(\a)}$ and $r={q^d-1\over q-1}$, or $G_\a$ is listed in Table {\rm \ref{stab-r>3}}.
\end{itemize}
\end{lemma}
\begin{table}[ht]
\[\footnotesize
\begin{array}{|l|l|l|l|l|l|l|}\hline
G_\a^{\Ga(\a)}&\ZZ_r&\Sy_r&\A_r&\PSL_2(11)&\M_{11}&\M_{23}\\ \hline
G_{\a}&(\ZZ_{l'}{\times} \ZZ_r){.}\ZZ_l&\Sy_r,\Sy_{r-1}{\times } \Sy_r & \A_r& \PSL_2(11)&\M_{11}&\M_{23}\\
&l'
\div l\div r-1& (\A_{r-1}{\times} \A_r){.}2& \A_{r-1}{\times} \A_r&
\A_5{\times}\PSL_2(11)& \M_{10}{\times} \M_{11}&\M_{22}{\times\M_{23}}\\

&& \ZZ_2^2{\times} \Sy_5& \ZZ_2^2{\times} \A_5 & &\A_6{\times } \M_{11}&  \\ \hline
\end{array}
\]
{\caption{\small Some stabilizers with $G_{\a\b}^{[1]}=1$.}\label{stab-r>3}}
\end{table}



\section{The proof of Theorem \ref{mainth-1}}\label{sct=general}
In this section,  $\Gamma=(V,E)$ is a connected
 graph of prime valency $r\ge 3$.
\begin{lemma}\label{t-mnsg-G}
Let $\Gamma=(V,E)$ be a connected
$G$-symmetric graph of  prime valency $r\ge 3$. Then
$G$ has at most one transitive minimal normal subgroup.
\end{lemma}
\proof
Suppose that $G$ has two distinct  transitive minimal normal subgroups, say $M$ and $N$.
Then $M\cap N=1$, and so $M$ and $N$ centralize each other. Thus
$M$ and $N$ are  nonabelian and regular on $V$, refer to
  \cite[pp.108-109, Lemma 4.2A and Theorem 4.2A]{Dixon-book}.
Clearly, $X:=MN$ is not regular on $V$. By Lemma \ref{N-p}, $\Ga$ is $X$-symmetric; in particular, $r\in \pi(X_\a)$ for $\a\in V$. Since $NX_\a=X=MX_\a$, we have $M\cong X/N\cong X_\a\cong X/M\cong N$. Then $|M|=|N|$ is divisible by $r$ but not by $r^2$. Noting that $M$ and $N$ are direct products of isomorphic simple groups, it follows that $M\cong T\cong N$ for some nonabelian simple group $T$. Then the actions of $N$ and $M$
 on $V$ are equivalent to the actions of $T$  on $T$ by right and left multiplications.
Identifying $V$ with $T$ and letting $\a=1\in T$,
the neighborhood $\Ga(\a)$ is in fact a conjugacy class of $T$ with length $r$,
 which is impossible by \cite[pp.36, Theorem 3.7]{Isaacs}.
This completes the proof.
\qed

\begin{corollary}\label{Cayley-1}
Let $\Gamma=(V,E)$ be a connected
$G$-symmetric graph of   prime valency $r\ge 3$.
Assume that $G$  has a  transitive simple subgroup   $T$.
If $T$ is normal in a characteristic subgroup of $G$ then $T$ is normal in $G$.
\end{corollary}
\proof
Assume $N$ is a characteristic subgroup of $G$, and $T\unlhd N$.
Then $T^g\unlhd N$ for each $g\in G$.
Suppose that $T\ne T^g$ for some $g$.
Then $N$ is not regular on $V$, and so $\Ga$ is $N$-symmetric.
Noting that $T\cap T^g=1$, it follows that $T$ and $T^g$ centralize each other, and hence they are regular on $V$. Thus $N$ has two transitive minimal normal subgroups, which contradicts Lemma \ref{t-mnsg-G}. Then the corollary follows.
\qed


\vskip 5pt

Recall that a permutation group
  $G$ is quasiprimitive if its minimal normal subgroups are all transitive.


\begin{lemma}\label{quasi-G-T}
Let $\Gamma=(V,E)$ be a connected
$G$-symmetric graph of   prime valency $r\ge 3$, and let $N=\soc(G)$.
Assume that $G$ is quasiprimitive on $V$ and $G$ has a  transitive simple subgroup $T$.
Then either  $r=7$, $T\cong \PSL_3(2)$ and $G\cong \AGL_3(2)$, or $N$ is simple with  $T\le N$. Moreover, for the latter case,  if $T\ne N$ then $\Ga$ is $N$-symmetric.
\end{lemma}
\proof
By Lemma \ref{t-mnsg-G}, $N$ is the unique minimal normal subgroup of $G$. Since $\Ga$ has odd valency, $|V|$ is even, and then $T$ is a nonabelian simple group.

Assume first that $N$ is abelian. Then $N\cong \ZZ_p^d$ and $G\lesssim \AGL_d(p)$ for
some prime $p$ and integer $d\ge 1$. In this case, $T\lesssim\GL_d(p)$, and $N$ is regular on $V$.
Since $\Ga$ has odd valency, $|N|=|V|$ is even, and so
  $p=2$. Since $T$ is transitive on $V$, we have $|T:T_\a|=2^d$.
 By \cite{Gur}, $d\ge 3$ and either $T\cong \A_{2^d}$ or $\PSL_n(q)$ with ${q^n-1\over q-1}=2^d$ and $n$ a prime. Note that $\A_{2^d}\not\lesssim\GL_d(2)$, see \cite[pp. 186, Proposition 5.3.7]{KL-book}. Then $T\cong  \PSL_n(q)$ and ${q^n-1\over q-1}=2^d$. By Lemma \ref{prime-div}, $n=2$ and $q=2^d-1$ is a prime.
By \cite[pp. 188, Theorem 5.3.9]{KL-book}, we have $d\ge {q-1\over (2,q-1)}=2^{d-1}-1$, yielding $d=3$. Then $N\cong \ZZ_2^3$, $|V|=8$ $r=7$, $T\cong \PSL_3(2)$  and $G\cong \AGL_3(2)$.

Next let $N$ be  nonabelian. Write $N=T_1\times\cdots\times T_k$ for integer $k\ge 1$ and isomorphic nonabelian simple groups $T_1,\ldots, T_k$.
Suppose that $T\not\le N$. Then $T\cap N=1$, and $TN/N\cong T$. Since $N$ is the unique minimal normal subgroup of $G$, we have
$\C_G(N)=1$, and thus $TN\cong TN/\C_{TN}(N)\lesssim\Aut(N)$. It follows that $T\cong TN/N\lesssim\Out(N)$, the outer automorphism group of $N$.
Noting that $\Out(N)\cong (\Out(T_1)\times\cdots\times \Out(T_k)){:}\Sy_k$, it follows that
$T\lesssim\Sy_k$. In particular,
$|T|$ is a divisor of $k!$.
Since $N$ is transitive on $V$, we have $G=NG_\a$ for $\a\in V$, and so
 $T\cong TN/N\le G/N\cong G_\a/(G_\a\cap N)$.
Then $G_\a$ is insolvable, and so $\Gamma$ is
$(G,2)$-arc-transitive. By \cite[Theorem 2]{Praeger-ONan},
$G$ satisfies III(b)(i) or III(c) given in
\cite[pp.229]{Praeger-ONan}. It follows  that there is a prime
divisor $p$ of $|N|$ such that $|V|$ is divisible by $p^k$. Since
$T$ is transitive on $V$, it follows that $p^k$ is a divisor of
$|T|$. Then $k!$ is  divisible by $p^k$, and so $k\le
\sum_{i=1}^\infty\lfloor{k\over p^i}\rfloor$, which is impossible.
Therefore, $T\le N$.

Assume that $T\ne N$. Then $N$ is not regular on $V$, and so  $\Gamma$ is $N$-symmetric by  Lemma \ref{N-p}. Suppose that $k>1$.
Then, by Lemma \ref{t-mnsg-G}, every $T_i$ is intransitive on $V$.
Noting that $\Ga$ is not bipartite, by Lemma \ref{N-p},
$|T_i|$ is a proper divisor of $|V|$, and so $|T_i|$ is a proper divisor of $|T|$ as $T$ is transitive
on $V$, where $1\le i\le k$.  Now consider the
projections $\phi_i: N\rightarrow T_i, x_1\cdots x_k\mapsto x_i$,
where $x_j\in T_j$ for $1\le j\le k$, and $1\le i\le k$. Then
$\phi_i(T)\ne 1$ for some $i$, and thus $T\lesssim T_i$, which is impossible.  Thus $k=1$, and $N$ is simple. This completes the proof. \qed

\vskip 5pt

The following corollary completes the proof of part (1) of Theorem \ref{mainth-1}.

\begin{corollary}\label{G-with-t-mnsg}
Let $\Gamma=(V,E)$ be a connected
$G$-symmetric graph of  prime valency $r\ge 3$.
Assume that $G$  has a transitive minimal normal subgroup $N$ and a transitive simple  subgroup $T$. Then either $r=7$, $T\cong \PSL_3(2)$ and $G\cong \AGL_3(2)$, or $T\le N$ and $N$ is simple.
\end{corollary}
\proof Choose a   maximal intransitive normal subgroup $M$ of $G$. Then $T\cap M=N\cap M=1$;
in particular, $MN=M\times N$. If $M=1$ then $G$ is quasiprimitive
on $V$, and so  the corollary holds by Lemma \ref{quasi-G-T}. Thus
let $M\ne 1$. This implies that $N$ is nonabelian as $M\le
\C_G(N)\ne N$.  Write $N=T_1\times\cdots\times T_k$ for integer
$k\ge 1$ and isomorphic nonabelian simple groups $T_1,\ldots, T_k$.
Then $G$ acts transitively on $\{T_1,\ldots,T_k\}$ by conjugation.
Let $\overline{G}=G/M$ and $\overline{T_i}=T_iM/M$ for $1\le i\le
k$. It is easily shown that $\overline{G}$ acts transitively on
$\{\overline{T_1},\ldots,\overline{T_k}\}$ by conjugation.
Let $\overline{N}=NM/M$.
Then $\overline{N}$ is a minimal normal subgroup of
$\overline{G}$. Let $\overline{T}=TM/M$.
Note that
$\overline{N}\cong N$ and $\overline{T}\cong T$.

Since $\Ga$ is not bipartite, $M$ has at least three orbits
on $V$. By Lemma \ref{N-p}, $|M|$ is a proper divisor of $|V|$, and so $|M|$ is a proper divisor of $|T|$.
Identify $\overline{G}$ with a subgroup of $\Aut(\Ga_{G/M})$.
Then Lemma \ref{quasi-G-T} works for $\Ga_{G/M}$.
We have that $\overline{T}\le \overline{N}$ and $\overline{N}$ is simple.
This yields that $T\le TM\le NM$. Suppose that $T\not\le N$. Then
$N\cap T=1$, and so $T\cong TN/N\le MN/N\cong M$. In particular,
$|T|$ is a divisor of $|M|$, a contradiction. Therefore, $N\ge T$, and the
corollary follows. \qed

\begin{corollary}\label{G-with-int-mnsg}
Let $\Gamma=(V,E)$ be a connected
$G$-symmetric graph of odd prime valency $r$. Let $K$ be a maximal intransitive normal subgroup of $G$.
Assume that $G$  has a  transitive simple subgroup $T$. Then either $K=1$, or $T\cong TK/K\le \soc(G/K)$.
\end{corollary}
\proof
Considering $\Ga_{G/K}$,  $TK/K$ and $G/K$, by Lemma \ref{quasi-G-T},
either the lemma holds, or  $\soc(G/K)\cong \ZZ_2^3$.
Assume that the latter case occurs. Then
$G/K\cong \AGL_3(2)$,
$T\cong TK/K\cong \PSL_3(2)$ and $r=7$. We next show $K=1$.

 Let $\soc(G/K)=N/K$, where $N\unlhd G$. Then $N$ is a regular subgroup of $G$.
Noting that  $|N|=|V|$, since $T$ is transitive on $V$, we know that $|V|$ is a divisor of $|T|=168$. It follows that  $|K|$ is a divisor of $21$. Thus $K$ is a subgroup of the cyclic group or the Frobenius group of order $21$. In particular, $|\Aut(K)|$ is indivisible by $8$.
Take a Sylow $2$-subgroup $P$ of $N$ and consider the conjugation of $P$ on $K$. Noting that $P\cong \ZZ_2^3$, it follows that
$\C_N(K)\ge \C_P(K)\ne 1$. Clearly, $\C_N(K)\unlhd G$.
Then $G/K$ has a (nontrivial) normal subgroup $\C_N(K)K/K$.
Recalling that $N/K$ is
the unique minimal normal subgroup of $G/K$, we have $N/K=\C_N(K)K/K$, yielding $N=\C_N(K)K$. It follows that $P\le \C_N(K)$ and $N=K\times P$.

Suppose that $K\ne 1$. Then we may choose a normal Hall subgroup $H$ of $N$   such that $P\le H$ and $|N:H|\in \{3,7\}$.
Note that $H$ is a characteristic subgroup of $N$. This implies that $H\unlhd G$. Then we get a symmetric graph $\Ga_{G/H}$ of valency $7$ and order $3$ or $7$, which is impossible. Thus   $K=1$, and the lemma follows.
\qed

\vskip 5pt

Now we are ready to prove part (2) of Theorem \ref{mainth-1}.

\begin{lemma}\label{Cayley-2}
Let $\Gamma=(V,E)$ be a connected
$G$-symmetric graph of  prime valency $r\ge 3$.
Assume that $G$  has a  transitive simple subgroup $T$, and   $T$ is not normal in $G$.  Then  either $G\cong \AGL_3(2)$, or $G$ has
  a characteristic subgroup $N$  such that
  $\Rad(N)$ is intransitive on $V$, $N/\Rad(N)$ is simple, $r\not\in\pi(\Rad(N))$ and
$\Rad(N)T\lneq N=N^\infty$.
\end{lemma}
\proof
We choose a minimal member $N$ from those characteristic subgroups of $G$ which
contain  $T$. Note that $T\le N^\infty$, and so $N=N^\infty$ by the choice of $N$.
Then $T$ is not normal in $N$ by Corollary \ref{Cayley-1}, and  $\Ga$ is $N$-symmetric by Lemma \ref{N-p}.

Assume first that $N$ is quasiprimitive on $V$. Then,  by the choice of $N$ and Lemma \ref{quasi-G-T}, either $N$ is simple or $N\cong \AGL_3(2)$. For the latter, noting that $\Rad(N)$ is normal in $G$, we have  $G\cong \AGL_3(2)$ by Corollary \ref{G-with-t-mnsg}. Then the lemma is true.

Assume that $N$ is not quasiprimitive on $V$.
Let $K$ be a maximal intransitive normal subgroup of $N$. Then $K\ne 1$ and, by Lemma \ref{quasi-G-T} and Corollary \ref{G-with-int-mnsg}, $T\cong TK/K\le \soc(N/K)$ and $\soc(N/K)$ is simple.
Recalling that $N=N^\infty$, it is easily shown that $\soc(N/K)=N/K$.
By the choice of $K$ and Lemma \ref{N-p}, $|K|$ is a proper divisor of $|V|$, and thus $|K|$ is a proper divisor of $|T|$.

For $\sigma\in \Aut(G)$, we have $K^\sigma\unlhd N^\sigma=N$, and so $K^\sigma K/K\unlhd N/K$. If $K^\sigma \ne K$ then $T\cong TK/K\le N/K=K^\sigma K/K\cong K^\sigma/(K^\sigma\cap K)$, and so $|T|$ is a divisor of $|K|$, a contradiction. Therefore, $K^\sigma =K$. It follows that $K$ is characteristic in $G$.

Set $X=KT$ and let $\a\in V$. Noting $K\cap T=1$ and $KT=X=TX_\a$, we have $|K|=|X:T|=|X_\a:T_\a|$.
Take two characteristic subgroups $I$ and $J$ of $K$
with $I\le J$. It is easy to see that  both $I$ and $J$ are   characteristic in $G$. Then $JT/I=(J/I)(TI/I)$, and $\C_{N/I}(J/I)$ is characteristic in $G/I$. Set $C/I=\C_{N/I}(J/I)$. Then $C$ is characteristic in $G$. Choose $J$ and $I$ such $J/I$ is a direct product of isomorphic simple groups, and write $J/I=S_1\times\cdots\times S_l$, where $S_i$ are simple.

Suppose that $K$ is insolvable. Then we may choose $J$ and $I$ such that $S_1\cong \cdots\cong S_l\cong S$ for some nonabelian simple group $S$.
In this case, $(C/I)\cap (J/I)=1$, and then $C/I$ is properly contained in $N/I$. By the choice of $N$, we have $l>1$ and
$TI/I\not \le C/I$; otherwise, we have $N/I=C/I$, yielding that $J/I$ is abelian, a contradiction. Then $T\cong TI/I\cong (TI/I)(C/I)/(C/I)\lesssim \Aut(J/I)\cong \Aut(S)\wr\Sy_l$. Recalling that $|K|$ is a proper divisor of $|T|$, it follows that $|S|$ is  a proper divisor of $|T|$, and thus we have $T\lesssim\Sy_l$.
Let $p$ be an odd prime of $|S|$. Then we have $p^l\le |K|_p\le |T|_p\le  |\Sy_l|_p$, yielding $l\le
\sum_{i=1}^\infty\lfloor{l\over p^i}\rfloor$, which is impossible.
Thus $K$ is solvable, and so $K=\Rad(N)$.

Suppose that $r\in \pi(K)$.
Then we may choose $J$ and $I$ such that $S_1\cong \cdots\cong S_l\cong \ZZ_r$. Recalling that $|K|=|X_\a:T_\a|$, we have $r^l\le |K|_r=|X_\a:T_\a|_r$.
This yields that $\Ga$ is $X$-symmetric as $r\in\pi(X_\a)$, and so $|X_\a|_r=r$ by Lemma \ref{N-p}. Then $|K|_r=r$, and so $l=1$ and $J/I\cong \ZZ_r$.
Thus $JT/I=(J/I)(TI/I)=(J/I)\times (TI/I)$, and then $TI/I\le C/I$. By the choice of $N$, we have $C=N$  and  $N/I=\C_{N/I}(J/I)$. This forces that $\Z(N/I)\ge J/I\cong \ZZ_r$.
Since $|K|_r=r$, we know that $J/I$ is a Sylow $r$-subgroup group
of $K/I$. Take a Sylow $r$-subgroup $Q/I$ of $TI/I$. Then $(J/I)\times (Q/I)$ is a  Sylow $r$-subgroup of $X/I$ and has order $r|T|_r$. By Lemma \ref{N-p}, $|N_\a|_r=r$, since $N=TN_\a$, we have $|N|_r= r|T|_r$. It follows that $(J/I)\times (Q/I)$ is a  Sylow $r$-subgroup of $N/I$. By Gasch\"{u}tz' Theorem (refer to \cite[pp. 31, (10.4)]{Aschbacher}), $N/I$ splits over $J/I$.
Let $M/I=(N/I)'$. Then $T\le M$, $M$ is  characteristic in $G$, and
$|M|_r=|T|_r<r|T|_r=|N_r|$, which contradicts the choice of $N$.
Therefore, $|K|_r=1$.

Finally, if $TK=N$ then $r|T|_r=|T|_r|N_\a|_r=|N|_r=|T|_r|K|_r$, yielding $|K|_r=r$, a contradiction. Thus $TK\ne N$, and the lemma follows.
\qed

\vskip 20pt

\section{Some limitations on the vertex-stabilizers}\label{sect-sl}
\begin{lemma}\label{Ga-insolvable}
Let $\Ga$, $G$, $T$ and $r$ be as in Hypothesis {\rm \ref{main-hypo}}. Then
$G$ is not an  exceptional simple group of Lie type, and
$G_\a$   has at most one insolvable composition factor.
\end{lemma}
\proof
Suppose that $G_\a$   has at least two  insolvable composition factors. By (\ref{stab-extension}) and Lemma \ref{stab-lem}, $G_\a$ has exactly two  insolvable composition factors, one is $S=\soc(G_\a^{\Ga(\a)})$ and the other one, say $S_1$, is a composition factor of $(G_\a^{\Ga(\a)})_\beta$ for $\b\in \Ga(\a)$.
We list them  as follows.
\[\small
\begin{array}{l|l|l|l|l|l}
S& \A_r &\PSL_d(q) &\PSL_2(11)& \M_{11}&\M_{23}\\ \hline
S_1& \A_{r-1}& \PSL_{d-1}(q) &\A_5&\A_6&\M_{22}
\end{array}
\]
Checking the factorizations given in \cite[Theorem 1.1]{Li-Xia19},
we conclude that there is no a simple group $T$ such that $G=TG_\a$,
a contradiction. Then  $G_\a$ has at most one insolvable composition factor. Noting  the factorization
$G=TG_\a$, our lemma follows from \cite[Theorem 1]{Factor-excep}. \qed

\vskip 5pt

In view of Lemma  \ref{Ga-insolvable},
we make some assumptions for the following sections.

\begin{hypothesis}\label{hypo-graphs}
 $\Ga=(V,E)$ is  a connected graph of prime valency $r\ge 3$ such that
\begin{itemize}
\item[(I)]  $\Ga$ is $G$-symmetric, where $G$ is a nonabelian simple group;
\item[(II)]  $G=TG_\a$ for $\a\in V$, where $T$ is a nonabelian simple group with $T\ne G$;

  \item[(III)] $G$ is not an exceptional simple group of Lie type, and $G_\a$ has at most one  insolvable composition factor.

\end{itemize}
\end{hypothesis}

By \cite{Trofimov}, (\ref{stab-extension}) and Lemmas \ref{r}, \ref{stab-lem} and \ref{Ga-insolvable}, we have the following lemma.

\begin{lemma}\label{stab-lem-insol}
Assume that Hypothesis {\rm \ref{hypo-graphs}} holds, and $G_\a$ is insolvable. Let $S=\soc(G_\a^{\Ga(\a)})$. Then  $\pi(G_\a)=\pi(S)$, and one of the following holds.
\begin{itemize}
\item[(1)]  $(S,G_\a,r)$ is one of $(\A_r,\A_r,r)$, $(\A_r,\Sy_r,r)$,
$(\PSL_2(11),\PSL_2(11),11)$, $(\M_{11},\M_{11},11)$ and $(\M_{23},\M_{23},23)$.

\item[(2)] $S\cong \PSL_d(q)$
 and
 $(r,|\O_p(G_\a)|,d,q)$ is listed as in Table {\rm  \ref{stab-O-p}}. In this case, if
 $\O_p(G_\a)=1$ then $G_\a\lesssim\PGammaL_d(q)$, if
 $\O_p(G_\a)\ne 1$ then
 $|G_\a:\O_p(G_\a)|$ is a divisor of $lf^2|\GL_d(q)|$, where $l=1$ or  $(l,d,q)$ is one of $(6,3,2)$ and $(24,3,3)$.
\end{itemize}
\end{lemma}
\begin{table}[ht]
\[\small
\begin{array}{|c|c|l|l|}\hline
r& \O_p(G_\a)&d&q=p^f\\ \hline
q+1& 1,\,q,\, q^2,\,q^3& 2& p=2, f={2^t}\\ \hline

{q^d-1\over q-1}& 1,\, q^{d-1},\,q^{d(d-1)\over 2}& \mbox{prime }d\ge 5  &     f=d^t\\ \hline
{q^d-1\over q-1}&  1,\, q^{d-1},\, q^{d}&  \mbox{prime }d\ge 3& f=d^t \\ \hline

 2^d-1&  1,\, 2^{d-1},\,2^{d+1}&  \mbox{prime }d\ge 3& q=2 \\ \hline
q^2+q+1& 1,\,q^2,\,q^6&  3& q\in \{2,3\}\\ \hline

7& 1,\,4,\, 2^{20}& 3& q=2\\ \hline

13& 1,\,9,\,3^{6}& 3& q=3 \\ \hline
31 & 1,\,2^4,\,2^{30}&  5&q=2\\ \hline
\end{array}
\]
{\caption{}\label{stab-O-p}}
\end{table}
\begin{remark} {\rm
Given Hypothesis \ref{hypo-graphs}, we
consider the (faithful) transitive action of
$G_\a$ on $[G:T]$ by right multiplication.  Let $n=|G:T|$.
Since $G$ is simple, we view $G$   as a transitive  subgroup of the alternating group $\A_n$, and thus we may let $T=\A_{n-1}\cap G$.
Then there is a general construction based on the coset graphs
for connected graphs satisfying the above (I) and (II), shown below.
Let $H$ be a
group which is a possible vertex-stabilizer of some
connected symmetric graph of valency $r$.
For any given faithful transitive permutation representation of $H$ on $\Delta=\{1,2,\ldots,n\}$ and a Hall $r'$-subgroup $K$ of $H$,
let $\overline{H}$ and $\overline{K}$ be the images of $H$ and $K$ respectively.
Determine the normalizer $\N_{\A_n}(\overline{K})$ of $\overline{K}$ in $\A_n$, and then seek the elements $x$ in $\N_{\A_n}(\overline{K})\setminus \overline{H}$ which have the property that $x^2\in \overline{K}$ and $\l x,\overline{H}\r$ is a simple group. Note that the simplicity of $\l x,\overline{H}\r$ ensures that $\overline{H}\le \A_n$.
Suppose that  there is such an $x$, and let $\overline{G}=\l x,\overline{H}\r$.
Then $\overline{G}$ has a factorization $\overline{G}=(\A_{n-1}\cap \overline{G})\overline{H}$.
If further $\A_{n-1}\cap \overline{G}$ is simple then the coset graph $\Cos(\overline{G},\overline{H},x)$ is a desired graph.
In view of this, we sometimes deal with $G$ and $G_\a$ as two transitive subgroups of $\A_n$.
}\qed
\end{remark}



\begin{lemma}\label{Ga-primsb-A-n}
Let $\Ga$, $G$, $T$ and $r$ be as in Hypothesis {\rm \ref{hypo-graphs}}. Let $n=|G:T|$ and   view  $G$ as a transitive subgroup of $\A_n$ acting     on $\Delta=\{1,2,\ldots,n\}$.
 Then either $G_\a$ is imprimitive on $\Delta$, or one of following holds.
\begin{itemize}
\item[(1)]   $G$, $T$ and $G_\a$ are listed in Table {\rm \ref{tab=over-Ga}} up to isomorphism.
\item[(2)] $G=\A_n$, $T=\A_{n-1}$, and $G_\a$ is almost simple but not $2$-transitive on $\Delta$.
\item[(3)] $n=r\ge 7$, $G=\A_r$, $T=\A_{r-1}$, $|G:T|_r=r$ and  $G_\a=\ZZ_r{:}\ZZ_{l}$, where $l$ is a divisor of $r-1\over 2$ with $(r,l)\ne (7,3)$.
\end{itemize}
\begin{table}[ht]
\[\tiny
\begin{array}{|l|l|c|c|c|l|l|}\hline
\mbox{Line}&r& G_\a  & G  & T  &|G:T|_r &\Ga \\ \hline
1&11& \A_{11}& \A_{12}& \M_{12}& 1&\K_{12} \\ \hline
2&23& \A_{23}& \A_{24}& \M_{24}&1 &\K_{24} \\ \hline
3&r=2^{2d-1}\pm 2^{d-1}-1,d\ge 3&\A_r&\A_{r+1}& \Sp_{2d}(2)& 1 &\K_{r+1}\\ \hline
4&r\not\in\{7,11,17,23\}&\A_r&\A_{r+1}& \PSL_{2}(r)&1 &\K_{r+1}\\ \hline

5&11& \M_{11}& \M_{12} &\M_{11} & 1 &\K_{12}\\ \hline
6&11& \PSL_2(11)& \A_{12} &\A_{11} &1 & \mbox{unique}\\ \hline
7&11& \ZZ_{11}& \PSL_2(11) &\A_{5} &11 & \mbox{Cayley}\\ \hline
8&11& \ZZ_{11}{:}\ZZ_5& \PSL_2(11) &\A_{5} &11 & \K_{12}\\ \hline
9&23& \ZZ_{23}& \M_{23} &\M_{22} &23 & \mbox{Cayley}\\ \hline

\end{array}
\]
{\caption{}\label{tab=over-Ga}}
\end{table}
\end{lemma}
\proof
Suppose that $G_\a$ is primitive on $\Delta$.
Since $G_\a$ has at most one insolvable composition factor, $G_\a$ is almost simple or of affine type.

{\bf Case 1}. Assume that $G_\a$ is almost simple.
Since $G$ is primitive on $\Delta$  with a simple point-stabilizer $T$, either $G=\A_n$ or   $(G,G_\a)$ can be read out from \cite[Tables III-VI]{maxsub-An} by \cite[Proposition 6.1]{Praeger90}.
Suppose  that $G\ne \A_n$. Then, combining Lemma \ref{stab-lem-insol}, either $(G,T,G_\a)$ is desired as
 at Lines 1-5 of Table \ref{tab=over-Ga}  and $\Ga\cong \K_{r+1}$, or $r=11$, $G_\a\cong \PSL_2(11)$, $G\cong  \M_{12}$ and $T\cong \M_{11}$. The latter is excluded by Lemma \ref{arc-stab-normalizer} and calculation of $\N_G(G_{\a\b})$, where $\b\in \Ga(\a)$.

Now let $G=\A_n$. Then $T=\A_{n-1}$ by the assumption in this lemma.
We assume that $G_\a$ is $2$-transitive on $\Delta$, and
check  $2$-transitive subgroups of $\A_n$, refer to\cite[pp. 197, Table 7.3]{Cameron}.
Combining Lemma \ref{stab-lem-insol}, $n$, $\soc(G_\a)$ and $r$  are listed as follows:
\[\tiny
\begin{array}{l|l|l|l}
n & \soc(G_\a)& r&\\ \hline
6,5&\PSL_2(4)& 5& \\ \hline
8,7&\PSL_3(2)&7& G_\a=\PSL_3(2) \\ \hline
12,11& \PSL_2(11)& 11& G_\a=\PSL_2(11)\\ \hline
12,11& \M_{11}& 11& G_\a=\M_{11}\\ \hline
23& \M_{23}& 23& G_\a=\M_{23}\\ \hline
15&\A_7&7&G_\a=\A_7\\ \hline
{q^d-1\over q-1}& \PSL_d(q)& {q^d-1\over q-1}>7& d\mbox{ a prime}\\ 
\end{array}\]
For all lines but the last one, since $\Ga$ has odd valency, $|V|=|G:G_\a|$ is even, we have $(n,r)\ne (5,5)$, $(7,7)$.
By calculation using GAP \cite{GAP}, we conclude that either $(n,r)=(6,5)$ or $n=12$ and $G_\a\cong \PSL_2(11)$. Noting that $\A_5\cong \PSL_2(5)$, the former case is included at Line 4 of Table \ref{tab=over-Ga}. The latter case implies Line 6 of Table \ref{tab=over-Ga}.

 We  suppose next that $n=r={q^d-1\over q-1}>7$ and $\soc(G_\a)=\PSL_d(q)$, and deduce a contradiction.
 Let $q=p^f$ for some prime $p$, and let $\sigma$ be the Frobenius automorphism of the  field of order $q$. Note that, by Lemma \ref{r}, $d$ is a prime, $(d,q-1)=1$ and $f$ is a power of $d$.
 Thus, by Lemma \ref{stab-lem}, we may let $G_\a=\SL_d(q){:}\l \sigma^l\r\le \GammaL_d(q)\cap \A_n$, where  $l$ is a divisor  of $f$.
  Without loss of generality, we identify $\Delta$ with the point-set of the projective geometry $\PG_{d-1}(q)$, and let $G_{\a\b}=q^{d-1}{:}\GL_{d-1}(q){:}\l\sigma^l\r$.
 Fix $\b\in \Ga(\a)$, and choose a $2$-element  $x\in \N_{\A_n}(G_{\a\b})$  with $(\a,\b)^x=(\b,\a)$ and $\A_n=\l x,G_\a\r$.

Suppose   that $d=2$, in this case, $p=2$ and $f=2^s$ for some $s\ge 2$. Then $G_{\a\b}$ is the  stabilizer  of some projective point $i$, and thus
 $x$ fixes $i$. Then
$G_{\a\b}\l x\r$ is a solvable $2$-transitive group of degree
$q$. It follows that $G_{\a\b}\l x\r\le \AGammaL_1(q)=\AGL_1(q){:}\l \sigma \r$. Thus $\l x, G_\a\r\le \SL_2(q){:}\l \sigma\r$, yielding $\l x, G_\a\r\ne \A_n=G$, a contradiction.

Now $d$ is an odd prime, and so $f$ is odd. Then $G_\a=\SL_d(q){:}\l \sigma^l\r\le \GammaL_d(q)\le \A_n$, and $\N_{\A_n}(G_{\a\b})$ has a subgroup $P{:}K\cong \AGammaL_{d-1}(q)$, where $P=\O_p(G_{\a\b})$ and $K\cong \GammaL_{d-1}(q)$ with $G_{\a\b}\le PK$.
Let  $C=\C_{\A_n}(P)$ and $N=\N_{\A_n}(P)$ . Then $P\le C$, $C\cap K=1$, $N\ge CK=C{:}K$, and so $K\lesssim N/C\lesssim\Aut(P)=\GL_{(d-1)f}(p)$.
Noting that $\N_{\A_n}(G_{\a\b})\le N$, we have 
$\N_{\A_n}(G_{\a\b})=
\N_N(G_{\a\b})$. Then $K\cong PKC/C\le \N_{\A_n}(G_{\a\b})C/C=\N_N(G_{\a\b})C/C\le \N_{N/C}(G_{\a\b}C/C)\lesssim\N_{\GL_{(d-1)f}(p)}(\GL_{d-1}(q){:}\l\sigma^l\r)\cong \GammaL_{d-1}(q)\cong K$. It follows $C{:}K=\N_N(G_{\a\b})C\ge \N_N(G_{\a\b})$. Write $G_{\a\b}=P{:}H$ with $H\le K$.
Recalling that $K\le \N_N(G_{\a\b})\le C{:}K$, we have $\N_N(G_{\a\b})=(C\cap \N_N(G_{\a\b}){:}K=P(C\cap H)K=PK$.
Then $x=yz\sigma^t$, where $y\in P$, $t$ is a divisor of $f$ and $z$ is contained in a subgroup of $H$ isomorphic to $\GL_{d-1}(q)$.
Since $f$ is odd,  the order of $x$
is not   a $2$-power unless $t=f$.   Then $x=yz \in G_{\a\b}$, a contradiction.

{\bf Case 2}.
Suppose  that $G_\a$ is not almost simple.
Then $4<n=p^m$ for some prime $p$ and integer $m\ge 1$, and
$G_\a\le \AGL_m(p)$. By   \cite[Proposition 6.2]{Praeger90},
one of the following holds: $G=\A_n$;  or $n=p=r=11$, $G=\PSL_2(11)$, $T\cong \A_5$ and $\ZZ_{11}\lesssim G_\a\lesssim\ZZ_{11}{:}\ZZ_5$; or $n=p=r=23$, $G=\M_{23}$, $T=\M_{22}$ and $\ZZ_{23}\lesssim G_\a\lesssim\ZZ_{23}{:}\ZZ_{11}$. Calculation shows that the last two case give  Lines 7-9 of Table \ref{tab=over-Ga}.

Assume that $G=\A_n$. Then $T=\A_{n-1}$. If
$G_\a$ is solvable, then $n=r=p$ and $G_\a=\ZZ_r{:}\ZZ_{l}$, where $l$ is a divisor of $r-1\over 2$. In this case, by calculation, we exclude the case $(r,l)=(7,3)$, and  part (3) of this lemma follows.
Thus we assume further that $G_\a$ is insolvable.
In particular, $m>1$, and $n\ge 8$.

Clearly, $\soc(G_\a)\le \O_p(G_\a)\ne 1$, and so $\Z(\O_p(G_\a))\ne 1$.
It follows that $\Z(\O_p(G_\a))$ is a regular normal subgroup of
$G_\a$ (acting on $\Delta$), and thus
$\soc(G_\a)=\Z(\O_p(G_\a))$.
Then, by Lemma \ref{stab-lem}, we may let
 $G_\a=\O_p(G_\a){:}((\ZZ_l\times \SL_d(q)){:}\l\sigma^e\r)$ and $r={q^d-1\over q-1}$, where
 $l$ is a divisor of $q-1$,  $\sigma$ is the Frobenius automorphism of the  field of order $q$, and  $q$ is a power of the prime $p$.
(Recall that, since $r$ is a prime,  $d$ is a prime and $(d,q-1)=1$.)
Replacing $T$ by its a conjugation in $G$, we let $(\ZZ_l\times \SL_d(q)){:}\l\sigma^e\r\le T$. Then $\Ga$ is $T$-symmetric
and $T_\a=(\ZZ_l\times \SL_d(q)){:}\l\sigma^e\r$, and thus
$T_\a^{[1]}\cong \ZZ_l$ and $T_{\a\b}^{[1]}=1$, where $\b\in \Ga(\a)$. By (\ref{stab-extension}), it is easily shown that $l=1$. Then we have $G_\a=\O_p(G_\a){:}(\SL_d(q){:}\l\sigma^e\r)$. It follows that
$p^m=|\O_p(G_\a)|=q^d$, and then $\ASL_d(q)\le G_\a\le \AGammaL_d(q)$.
 In particular, $G_\a$ is $3$-transitive on $\Delta$.
Then, noting that $\A_n=\A_{n-2}G_\a$, we  have a factorization
$T=\A_{n-2}T_\a$ with $T_\a=\SL_d(q){:}\l\sigma^e\r$.
Moreover, $T_\a$ is a $2$-transitive subgroup of $\A_{n-1}$.
Recalling that $n\ge 8$, we have $n-2\ge 6$.
Then the argument in Case 1
works for the triple $(n-1,n-2,r)$. Thus, noting that $n=q^d$, the only possibility is that $n=8$, $r=7$, and $T_\a\cong \PSL_3(2)$.
Then $G_\a=2^3{:}\PSL_3(2)$, which has odd index in $\A_8$, a contradiction.
\qed


\vskip 20pt

\section{Graphs with solvable vertex-stabilizers}\label{sect=sol-stab}

In this section, we  deal with the graphs
described as in {\rm Hypothesis \ref{hypo-graphs}} under the further assumption that  $G_\a$  is solvable. We always let $n=|G:T|$, and  view  $G$ as a transitive subgroup of $\A_n$ (acting     on $\Delta=\{1,2,\ldots,n\}$).

By \cite[Theorem 1.1]{Li-Xia}, combining   Table \ref{cubic-stab} and Lemma \ref{stab-lem}, $G$, $T$ and $G_\a$ are listed as in Table \ref{tab=Ga-sol} up to isomorphism.
\begin{table}[ht]
\[\footnotesize
\begin{array}{|l|l|c|c|c|l|}\hline
\mbox{Line} &G  & T  & G_\a & r &   \\ \hline
1&\PSL_2(r)& \A_5& \ZZ_r{:}\ZZ_{r-1\over 2}& r\in \{11,19,29,59\}&  \\ \hline
2&\PSL_2(11)& \A_5& \ZZ_{11} & 11&  \\ \hline
3&\PSL_2(19)& \A_5& \ZZ_{19}{:}\ZZ_3& 19&  \\ \hline

4&\PSL_2(29)& \A_5& \ZZ_{29}{:}\ZZ_7& 29&  \\ \hline

5&\A_c&\A_{c-1}& &    &  c\ge 6\\ \hline


6&\A_6&\A_5&   & r=3& \\ \hline

7& \M_{c}& \M_{c-1}&   & & c\in \{12,24\} \\ \hline

8& \M_{23}& \M_{22}& \ZZ_{23},\ZZ_{23}{:}\ZZ_{11} &23 &  \\ \hline
\end{array}
\]
{\caption{}\label{tab=Ga-sol}}
\end{table}
Note that, since  $G_\a$ is solvable,  if $r\ge 5$ then $G_\a$
has a unique Sylow $r$-subgroup, which is semiregular on $\Delta$.

\begin{theorem}\label{Ga-solvable}
Assume that {\rm Hypothesis \ref{hypo-graphs}} holds  and   $G_\a$ is solvable.
Let $n=|G:T|$.
Then $n$ is divisible by $r$, and one of the following holds.
\begin{itemize}
\item[(1)]  $\Ga\cong \K_{r+1}$ for $r\in \{11,19,29,59\}$, $G=\PSL_2(r)$, $T\cong \A_5$ and $G_\a\cong \ZZ_r{:}\ZZ_{r-1\over 2}$.
\item[(2)] $n=11$, $r=11$, $G=\PSL_2(11)$,  $T\cong \A_5$ and $G_\a\cong \ZZ_{11}$.

\item[(3)] $n=203$, $r=29$, $G=\PSL_2(29)$,  $T\cong \A_5$ and $G_\a\cong \ZZ_{29}{:}\ZZ_7$.

\item[(4)]
 $n=r=23$, $G=\M_{23}$,  $T=\M_{22}$ and $G_\a\cong \ZZ_{23}$.

\item[(5)] $G=\A_n$, $T=\A_{n-1}$, $n\ge 7$, $r\ge 5$ and $r^2$ is not a divisor of $n$.

\item[(6)] $G=\A_{48}$, $T=\A_{47}$, $r=3$ and $G_\a\cong \ZZ_2\times \Sy_4$.
\end{itemize}
\end{theorem}
\proof
First, Lines 1, 2 and 4 of Table \ref{tab=Ga-sol} give parts (1), (2) and (3) of this lemma, respectively. For Line 3 of Table \ref{tab=Ga-sol},   $G$ has a unique conjugacy class of elements of order $3$. This yields that $T\cap G_\a$ contains
an element of order $3$, see \cite[Lemma 2.2]{sps-factorization} for example. Then
$3420=|G|={|T||G_\a|\over |T\cap G_\a|}\le 1140$, a contradiction. For Line 8 of Table \ref{tab=Ga-sol}, if
$G_\a$ has order $253$, then $G_{\a\b}\cong \ZZ_{11}$ for $\b\in \Ga(\a)$,  and $|\N_G(G_{\a\b})|$ has order $55$ (confirmed by GAP), which  contradicts Lemma \ref{arc-stab-normalizer}. Thus $G_\a\cong \ZZ_{23}$, and part (4) of this lemma follows.

Suppose that Line 6 of Table \ref{tab=Ga-sol} occurs. Noting that $G_\a$ is known as in Table \ref{cubic-stab}, we have
$G_\a\cong\Sy_3$ or $\Sy_4$, and
$|V|=60$ or $15$ respectively.
Since $\Ga$ has odd valency, $|V|$ is even, and so $G_\a\cong\Sy_3$. Then $|G|=|T||G_\a|$, and so
$T\cap G_\a=1$.
Thus $G_\a$ acts regularly on $\Delta$. This yields that
$G_\a$ contains odd permutations of $\Sy_6$,   a contradiction.

Suppose that Line 7 of Table \ref{tab=Ga-sol} occurs. Then
$n=c\in \{12,24\}$.
If $r\geqslant5$ then the Sylow $r$-subgroup of $G_\a$ is semiregular on $\Delta$, and so $r$ is a divisor of $n$,
which is impossible. Thus we have  $r=3$.
If $T\cap G_\a=1$ then $\Ga$ is isomorphic to a symmetric
Cayley graph of  $T$, and so $T$ is normal in $G$ by \cite{CubicC}, which is impossible.  Therefore,  $T\cap G_\a\ne 1$.
Recalling that $n=|G:T|=|G_\a:(T\cap G_\a)|$, it follows that $n$ is a proper divisor of $|G_\a|$. By Table \ref{cubic-stab},
$(n, G_\a)$ is one of
$(12, \Sy_4)$, $(12, \Sy_4\times \ZZ_2)$ and $(24, \Sy_4\times \ZZ_2)$. 
With the help of   GAP, up to conjugation, $\A_n$ has at most three transitive subgroups isomorphic to $G_\a$, and none of them
has a Sylow $2$-subgroup $P$ such that $|\l G_\a, x\r|=|\M_{n}|$ for some $x\in \N_{\A_{n}}(P)$. Noting that $P=G_{\a\b}$ for some $\b\in \Ga(\a)$,
we get a contradiction, see Lemma \ref{arc-stab-normalizer}.

Finally, suppose that Line 5 of Table \ref{tab=Ga-sol}  occurs. Then $n=c$, and $n=|G_\a:(T\cap G_\a)|$  as
$G_\a$ is a transitive subgroup of $\A_n$.
If $r\ge 5$ then the  Sylow $r$-subgroup of $G_\a$ is semiregular on $\Delta$, and so $n$ is divisible by $r$, yielding part (5) of this lemma.
Assume that $r=3$. Then
$G_\a\cong \Sy_3,\,\D_{12}$, $\Sy_4$ or $\Sy_4\times\ZZ_2$. Suppose that  $T\cap G_\a\ne 1$. Recalling that $6\le n=|G_\a:(T\cap G_\a)|$, we have $G_\a\cong\D_{12},\, \Sy_4$ or $\Sy_4\times\ZZ_2$. It follows that $(n, G_\a)$ is one of $(6,\D_{12})$, $(6,\Sy_4)$,
$(8,\Sy_4)$, $(8, \Sy_4\times \ZZ_2)$, $(16, \Sy_4\times \ZZ_2)$, $(12, \Sy_4)$, $(12, \Sy_4\times \ZZ_2)$ and $(24, \Sy_4\times \ZZ_2)$. With the help of   GAP, a similar argument as in the above paragraph leads to a contradiction. Thus   $T\cap G_\a=1$. Then $\Ga$ is a Cayley graph of $\A_{n-1}$, and so part (6) of this lemma holds by \cite{CubicC}.
\qed

\begin{remark}
{\rm For   $r\in \{5,7\}$,  there is a graph satisfying Theorem \ref{Ga-solvable}~(5) if and only if $(n,r,|G_\a|)$ is one of $(20,5,40)$, $(40,5,40)$, $(40,5,80)$, $(80,5,80)$, $(7,7,7)$, $(7,7,21)$, $(21,7,21)$, $(21,7,63)$, $(63,7,63)$ and $(84,7,84)$. This is confirmed by GAP.
}
\end{remark}

\begin{corollary}\label{prime-dim}
Assume that {\rm Hypothesis \ref{hypo-graphs}} holds, and that  $G$ is isomorphic to a linear simple group of prime dimension. Then one of the following holds.
\begin{itemize}
\item[(1)]  $\Ga\cong \K_{r+1}$ for $r\in \{11,19,29,59\}$, $G=\PSL_2(r)$, $T\cong \A_5$ and $G_\a\cong \ZZ_r{:}\ZZ_{r-1\over 2}$.
\item[(2)] $n=11$, $r=11$, $G=\PSL_2(11)$,  $T\cong \A_5$ and $G_\a\cong \ZZ_{11}$.

\item[(3)] $n=203$, $r=29$, $G\cong \PSL_2(29)$,  $T\cong \A_5$ and $G_\a\cong \ZZ_{29}{:}\ZZ_7$.
\end{itemize}

\end{corollary}
\proof
Since $G=TG_\a$ and $T$ is nonabelian simple,  by \cite[Corollary 3.4]{Li-Xia},  $G_\a$ is solvable. By Theorem \ref{Ga-solvable},
either one of parts (1)-(3) follows, or
$G=\A_n$, $T=\A_{n-1}$, $r\ge 5$ and $r$ is a divisor of $n$. For the latter,  we have $G\cong \PSL_2(9)$ but $n=6$, a contradiction.
\qed

\begin{corollary}\label{G-U-2m+1}
Assume that {\rm Hypothesis \ref{hypo-graphs}} holds.
Then $G$ is not (isomorphic to) a simple unitary group of odd dimension.
\end{corollary}
\proof
Suppose that $G\cong \PSU_{2m+1}(q)$.
Then, by Theorem \ref{Ga-solvable}, $G_\a$ is insolvable.
Thus  $G=AB$ for maximal subgroups $A$ and $B$ of $G$ such that
  one of $A$ and $B$ contains $T$ and the other one contains $G_\a$.  Then, checking the pairs $(A,B)$   given in
\cite[pp. 13, Table 3]{max-factor}, the only possibility is that $G\cong \PSU_{9}(2)$, $A\cong \J_3$ and $B\cong [2^{15}]{:}\PSU_7(2)$. If $T\le B$ then $G=AT$;
however,
$|A|_2|T|_2\le 2^{28}<2^{36}=|G|_2$, a contradiction. Thus $T\le A$ and $G_\a\le B$.
It is easy to check that $|B|$ is indivisible by $19$ or $17$, and $|A|$ is indivisible by $43$ or $11$. Check the maximal subgroups of $A$ and $\PSU_7(2)$, refer to the Atlas \cite{Atlas}. We conclude that $T=A$ and $\PSU_7(2)\preccurlyeq G_\a$,
which contradicts Lemma \ref{stab-lem-insol}.
This completes the proof.
\qed

\vskip 20pt

\section{Graphs arising from alternating groups}\label{sect-alt}

In this section, using some known factorizations of simple groups, we investigate the graphs in Hypothesis \ref{hypo-graphs} under the further assumption that $G_\a$ is insolvable and  either $G$ or $T$   is isomorphic to an alternating group.  Let $n=|G:T|$, view  $G$ as a transitive subgroup of $\A_n$
acting     on $\Delta=[G:T]$.

\begin{lemma}\label{Ga=A-r}
Assume that {\rm Hypothesis \ref{hypo-graphs}} holds.
Let $n=|G:T|$.
Assume that $|\Ga(\a)|=r\ge 5$ and $G_\a \cong \A_r$. Then
one of the following holds.
\begin{itemize}
\item[(1)] $G=\A_6$, $T=\A_5$, $G_\a=\PSL_2(4)$ and $\Ga\cong \K_6$.

\item[(2)] $n>{r\choose 2}$,
and either $n\ge {r\choose {r-1\over 2}}$  or $r\ge 11$ and   $n$ is divisible by $r\choose l$ for some $l$ with $2\le l<{r-1\over 2}$. Moreover, one  the following holds:
\begin{itemize}
\item[(i)]
 $\Ga\cong \K_{r+1}$ for $r\ge 7$, $G\cong \A_{r+1}$,  and $T$ is isomorphic to a transitive subgroup of $\A_{r+1}$;
 \item[(ii)]  $G=\A_n$,  $T=\A_{n-1}$ and $G_\a$ is a transitive subgroup of $\A_n$.
 \end{itemize}
  \end{itemize}
\end{lemma}
\proof
Note that $n>r$, $2n$ is a divisor of $r!$ and $\A_r$ has a subgroup $H$ of index $n$. Let $k={r-1\over 2}$.
Suppose that $n<{r\choose k}$. By \cite[Theorem 5.2A]{Dixon-book}, either
$\A_{r-l}\lesssim H\lesssim (\Sy_l\times \Sy_{r-l})\cap \A_r$ for
some $l$ with $2\le l<k$, or $(r,k,n,H)$ is one of $(5,2,6,\D_{10})$ and $(7,2,15,\PSL_3(2))$.
For $(r,k,n,H)=(5,2,6,\D_{10})$, we have $G=\A_6$, $T=\A_5$, $G_\a=\PSL_2(4)$ and $\Ga\cong \K_6$ as in part (1) of this lemma. Assume that $(r,k,n,H)=(7,2,15,\PSL_3(2))$.
Up to equivalence, $\A_7$ has two permutation representations of degree $15$. Given such a representation of $G_\a$, using GAP, we search the $2$-elements $x$ in $\A_{15}\setminus G_\a$ which normalize a subgroup $\A_6$ of $G_\a$. Computation shows that
$\l x,G_\a\r\cong \A_8$ for each a desired $x$. We have $G\cong \A_8$; however, a subgroup of index $15$ in $\A_8$ is isomorphic to $2^3{:}\PSL_3(2)$, which can not play the role of $T$, a contradiction. Thus we assume next that either $n\ge {r\choose {r-1\over 2}}$, or $n<{r\choose {r-1\over 2}}$ and $n$ is divisible by $r\choose l$ for
some $l$ with $2\le l<{r-1\over 2}$.
Moreover, for the latter case, $r\ge 11$. In fact, if  $r\in \{5,7\}$ and $n<{r\choose {r-1\over 2}}$ then $(n,r)=(21,7)$, which is excluded by calculation in a similar way as above.

Suppose that $n={r\choose 2}$. Then either $r=5$ or $r\ge 11$.
Let $M$ be a maximal subgroup of $G_\a$ with $H\le M$.
Assume that $M$ is a primitive subgroup of $\A_r$. Then $|H|\le |M|<4^r$ by
\cite{order-primgp}, and so ${r\choose 2}=n> {r!\over 4^r2}$, yielding $r\le 13$. Then $(n,r)=(10,5)$, $(21,7)$, $(55,11)$ and $(78,13)$. However, for such an $r$, every primitive subgroup of $\A_r$ has index indivisible by ${r\choose 2}$, a contradiction.
Thus $M$ is an intransitive or imprimitive subgroup of $\A_r$.
It follows from \cite{maxsub-An} that
$M$ is the stabilizer of some $2$-subset of $\{1,2,\ldots,r\}$, and thus $H=M\cong \Sy_{r-2}$.
Then the action of $G_\a$ on $\Delta=[G:T]$ is equivalent that
on the   $2$-subsets of an $r$-set.
Thus, for $\b\in \Ga(\a)$, the arc-stabilizer $G_{\a\b}$ has
two orbits on $\Delta$, say $\Delta_1$ and $\Delta_2$, which have length $r-1$ and $(r-1)(r-2)\over 2$ respectively.
Moreover, it is easily shown that $G_{\a\b}$ acts primitively on
both $\Delta_1$ and $\Delta_2$. Noting that $|\Delta_1|\ne |\Delta_2|$, we know that
 $\C_G(G_{\a\b})$  fixes both  $\Delta_1$ and $\Delta_2$ set-wise. By \cite[Theorem 4.2A]{Dixon-book}, $\C_G(G_{\a\b})$ acts
 trivially on both $\Delta_1$ and $\Delta_2$, yielding
 $\C_G(G_{\a\b})=1$. Recalling $r\ne 7$, we have $\N_G(G_{\a\b})\lesssim
 \Sy_{r-1}$. Note that the transposition $(1\,2)$ induces a product of $r-1$ transpositions on $\Delta$. This implies that
 $\Sy_r\lesssim\N_G(G_{\a})$, and thus $\N_{G_\a}(G_{\a\b})\gtrsim\Sy_r$. Then $\N_G(G_{\a\b})\le G_\a$, a contradiction. Thus  $n>{r\choose 2}$.

Since $G_\a\cong \A_r$ and $G=TG_\a$, by \cite[Theorem 1.1]{XIA},
one of the following occurs:
\begin{itemize}
\item[(i)] $G\cong \A_{r+1}$, $G_\a\cong \A_{r}$, $T$ is isomorphic to a transitive subgroup of $\A_{r+1}$;
\item[(ii)] $G=\A_n$, $G_\a\cong \A_{r}$,  $T=\A_{n-1}$ and $G_\a$ is a transitive subgroup of $\A_n$;
\item[(iii)] $G\cong \A_{r+k}$, $G_\a\cong \A_{r}$, and $T$ is isomorphic to a $k$-transitive subgroup of $\A_{r+k}$, where $2\le k\le 5$;
\item[(vi)] the triple
$(G,r,T)$ is one of $(\M_{12},5,\M_{11})$, $(\Omega_8^+(2),7,\Sp_6(2))$, $(\A_{15},13,\A_{13})$, $(\A_{15},13,\A_{8})$  and $(\A_{15},7,\A_{13})$.
\end{itemize}
For the case (i) above, $\Ga$ is the complete graph $\K_{r+1}$.
If $r=5$ then $G\cong \A_6$, yielding $T\cong \A_5$, and so
$6=n>{5\choose 2}=10$, which is impossible. Then $r\ge 7$, desired as in part (i) of this lemma.
For the case (ii) above, we have part (ii) of this lemma.
To complete the proof, we next exclude the left two cases.

{\bf Case 1}. Assume that case (iii) occurs. Consider the action of $G$ on an $(r+k)$-set,
and take a $G_\a$-orbit with maximal length, say $m$. Then
$m\ge r$, and $G_\a$ has a subgroup $K$ of index $m$.
Clearly, $m\le r+k\le 2r\le {r\choose 2}$.
Suppose that $m={r\choose 2}$. Then $r(r-1)=2m$, the only possibility is that $r=5$ and $m=10$. This implies that
$k=r=5$.  Thus $T$ is (isomorphic to) a $5$-transitive subgroup of $\A_{10}$. This yields that $T=\A_{10}=G$, a contradiction.
Therefore, $|G_\a:K|=m<{r\choose 2}$ and $m\le r+5$.

By  \cite[Theorem 5.2A]{Dixon-book},
either $(r,m)=(5,6)$, or $K\lesssim \A_{r-1}$. For the latter case,
$r$ is a divisor of $m=|G_\a:K|$, since $m\le r+k$ and $2\le k\le 5$, we have  $r=k=5$, which yields that $T$ is a $5$-transitive subgroup of $\A_{10}$, a contradiction. If $r+k\ge 9$ and $r=5$ then  $k\ge 4$, and we have a similar contradiction as above.
All in all,  $(r,m)=(5,6)$  and $7\le r+k<9$.
Computation using GAP shows that, for $\b\in \Ga(\a)$, there is no $2$-element $x\in \N_G(G_{\a\b})$ with $G=\l G_\a,x\r$, which contradicts Lemma \ref{arc-stab-normalizer}.

{\bf Case 2}. Assume that case (iv) occurs.
Suppose that $(G,r,T)=(\Omega_8^+(2),7,\Sp_6(2))$. Then, as a permutation group on $\Delta$, the group $G$ has  $2$ transitive  subgroups $\A_7$ up to conjugacy, confirmed by GAP.  Taking $G_\a$ as any one of them, further calculation shows that there is no $2$-element  $x\in \N_G(G_{\a\b})$ with $x^2\in G_{\a\b}$ and $G=\l G_\a,x\r$,
 which contradicts Lemma \ref{arc-stab-normalizer}. Suppose that $(G,r,T)=(\M_{12},5,\M_{11})$. Then
 $G_\a$ is a transitive subgroup of $\A_{12}$  acting on a $12$-set.  Note that, up to equivalence, $G_\a$ has a
unique permutation representation of degree $12$, and
 $G_{\a\b}$ is regular under this representation, where $\b\in \Ga(\a)$. Computation  shows that there is no $x\in \N_{\A_{12}}(G_{\a\b})$ with $x^2\in G_{\a\b}$ and $|\M_{12}|=|\l G_\a,x\r|$,
 a contradiction.

 For the left three triples, we conclude that $G_\a\cong \A_r$ has a   $2$-transitive permutation representation of degree $15$. It follows that $r=7$.  Considering the transitive permutation representations of $\A_7$ of degree $15$, we get a similar contradiction by calculation.
This completes the proof.
\qed


\begin{theorem}\label{G-or-T-alte}
Let $\Ga$, $G$, $T$ and $r$ be as in Hypothesis {\rm\ref{hypo-graphs}}, and let $n=|G:T|$.
Assume that  $G$ or $T$ is isomorphic to an alternating group, and that $G_\a$ is insolvable, where $\a\in V$. Then one of the following holds.
\begin{itemize}

\item[(1)] $G=\A_6$, $T=\A_5$, $G_\a=\PSL_2(4)$, $r=5$ and $\Ga\cong \K_6$.
\item[(2)] $G\cong\M_{12}$, $T\cong \A_5$, $G_\a\cong \M_{11}$, $r=11$ and $\Ga\cong \K_{12}$.

\item[(3)] $G$, $T$ and $G_\a$ are given in Table {\rm \ref{tab=over-Ga}}.

\item[(4)] $G=\A_n$, $T=\A_{n-1}$,  $G_\a$ is   a transitive but not $2$-transitive  subgroup of $\A_n$.

\item[(5)] $n>{r\choose 2}$ for $r\ge 7$, $\Ga\cong \K_{r+1}$, $G\cong \A_{r+1}$, $G_\a\cong \A_r$, $T$ is isomorphic to a transitive subgroup of $\A_{r+1}$, and either $n\ge {r\choose {r-1\over 2}}$  or $r\ge 11$ and   $n$ is divisible by $r\choose l$ for some $l$ with $2\le l<{r-1\over 2}$.
\end{itemize}
\end{theorem}
\proof
If $G_\a\cong \A_r$ then, by Lemma \ref{Ga=A-r} one of parts (1), (4) and (5) occurs. Thus  we  assume next that $G_\a\not\cong \A_r$, and
 discuss in two cases: $T\cong \A_c$ for some $c\ge 5$;  $G\cong \A_{c+1}$ for some $c\ge 5$, but   $T$ is not isomorphic to an alternating group.

{\bf Case 1}. Assume that $T\cong \A_c$ for some $c\ge 5$. By \cite[Theorem 1.1]{XIA}, combining Lemma \ref{stab-lem-insol}, one line of the following table holds.
\[\tiny
\begin{array}{|l|l|l|l|}\hline
G& T& G_\a &r \\ \hline
\M_{12}&\A_5&\M_{11}&  11\\ \hline
\A_{15}&\A_7,\A_8&\Sy_{13}& 13\\ \hline
\PSL_4(3)&\A_6&3^3{:}\PSL_3(3)&13\\ \hline
\Omega_7(3)& \A_8,\A_9& [3^6]{:}\PSL_3(3)&13\\ \hline
\Omega_7(3)& \A_9& 3^3{:}\PSL_3(3)&13\\ \hline
\Omega^+_8(2)& \A_9& [2^4]{:}\SL_2(4)\le G_\a\le [2^6]{:}\GammaL_2(4) \mbox{ (see \cite{567})}&5\\ \hline
\A_{c+k}& \A_{c}& k\mbox{-transitive subgroup of }\A_{c+k}, 1\le k\le 5&\\ \hline
\end{array}
\]
If $G\cong \M_{12}$ then $\Ga$ is the complete graph $\K_{12}$, and part (2) of this lemma follows.

 For
$G\cong \A_{15}$ and $G_\a\cong \Sy_{13}$, the action of $G$ on $V$ is equivalent to that on the $2$-subsets of a $15$-set, and thus $G_\a$ has no orbit of length $13$ on $V$, a contradiction.
For
$G\cong \PSL_4(3)$ and $G_\a\cong 3^3{:}\PSL_3(3)$, the action of $G$ on $V$ is $2$-transitive, yielding $\Ga\cong \K_{40}$, which is not the case.

Suppose that $G\cong \Omega_7(3)$. If $G_\a\cong [3^6]{:}\PSL_3(3)$ then $G$ is a primitive group on $V$ with subdegrees $1$, $39$, $351$ and $729$, and so $G_\a$ has no orbit of length $13$ on $V$, a contradiction. Let $G_\a\cong 3^3{:}\PSL_3(3)$. Then $G_{\a\b}\cong [3^5]{:}2\Sy_4$ for $\b\in \Ga(\a)$. Checking the maximal subgroups of $G$ in the Atlas \cite{Atlas}, $G_\a$ is contained in a maximal subgroup $[3^6]{:}\PSL_3(3)$ or $\PSL_4(3).2$. Computation shows that $|\l x, G_\a\r|\le 2|\PSL_4(3)|<|G|$ for
$x\in \N_G(G_{\a\b})$ with $x^2\in G_{\a\b}$, a contradiction.
For $G\cong \Omega^+_8(2)$, we know that $G_\a$ is one of $[4^2]{:}\SL_2(4)$, $[4^2]{:}\GL_2(4)$, $[4^2]{:}\GammaL_2(4)$, $[4^3]{:}\GL_2(4)$ and $[4^3]{:}\GammaL_2(4)$, see \cite[Theorem 3.4]{567} for example.
Confirmed by GAP, all subgroups of $\Omega^+_8(2)$ isomorphic to $G_\a$
are conjugate under $\Aut(G)$. Then calculation gives a similar contradiction as above.

Now $G\cong \A_{c+k}$, $T\cong  \A_{c}$ and $G_\a$ is (isomorphic to) a $k$-transitive subgroup of $\A_{c+k}$, where $1\le k\le 5$. If $k=1$ then $n=c+1$, and thus either part (3) of this lemma occurs by Lemma \ref{Ga-primsb-A-n} when $G_\a$ is  a $2$-transitive subgroup of $\A_{n}$ or, otherwise we have part (4) of this lemma.
Now  let $k>1$. Then Lemma \ref{Ga-primsb-A-n} works for the triple $(\A_{c+k}, \A_{c+k-1},G_\a)$. In this case, the only possibility is that $c+k=12$ and $G_\a \cong \PSL_2(11)$,
and so $\Ga$ is desired as in  part (3) of this lemma.

{\bf Case 2}.
Assume that $G\cong \A_{c+1}$ and $T$ is not (isomorphic to) an alternating group. By \cite[Theorem D]{max-factor}, combining  Lemma \ref{stab-lem-insol},    $T$ is a $k$-homogeneous subgroup of $\A_{c+1}$, and either $r>5$ and $G_\a\cong \Sy_r$ or
 $r=5$ and $c+1<9$,
where $1\le k\le 5$ and $c+1=r+k$. If $k=1$ then $c=r$, and so $\Sy_5\lesssim\A_6$, which is impossible. Then $k\ge 2$.

Consider the action of $G$ on a $(c+1)$-set, and take a $G_\a$-orbit with maximal length, say $m$. Then
$m\ge r$, and $G_\a$ has a subgroup $K$ of index $m$.
Clearly, $m\le c+1\le 2r\le {r\choose 2}$.
Suppose that $m={r\choose 2}$. Then $r(r-1)=2m$, the only possibility is that $r=5$ and $m=10$. This implies that
$k=r=5$, and $c+1=10$.  Thus $T$ is (isomorphic to) a $5$-transitive subgroup of $\A_{10}$, yielding $T\cong \A_{10}\cong G$, a contradiction.
Therefore, $|G_\a:K|=m<{r\choose 2}$ and $m\le r+5$.

By  \cite[Theorem 5.2B]{Dixon-book},
  $(r,m)=(5,6)$  or $K\lesssim \Sy_{r-1}$. For the latter case,
$r$ is a divisor of $m=|G_\a:K|$, since $m\le c+1$ and $2\le k\le 5$, we have  $r=k=5$, which yields that $T$ is a $5$-transitive subgroup of $\A_{10}$, a contradiction. Thus we have $(r,m)=(5,6)$  and $7\le r+k<9$, yielding $c+1=8$.
Computation   shows that, for $\b\in \Ga(\a)$, there is no $x\in \N_G(G_{\a\b})$ with $x^2\in G_{\a\b}$ and $G=\l G_\a,x\r$, which contradicts Lemma \ref{arc-stab-normalizer}.
\qed

\vskip 5pt

In the following theorem, we consider the case where $G$ is  a sporadic simple group.
\begin{theorem}\label{G-sporadic}
Assume that {\rm Hypothesis \ref{hypo-graphs}} holds, and
 let $n=|G:T|$.
Assume that
$G$ is isomorphic to a sporadic simple group. Then one of the following holds.
\begin{itemize}
\item[(1)] $n=r=23$, $G=\M_{23}$, $T=\M_{22}$, $G_\a\cong \ZZ_{23}$.
\item[(2)] $n=12$, $r=11$, $G=\M_{12}$, $T=\M_{11}$, $G_\a\cong \M_{11}$, and $\Ga\cong \K_{12}$.

\item[(3)] $r=11$,  $G\cong\M_{12}$, $T\cong \A_5$, $G_\a\cong \M_{11}$, and $\Ga\cong \K_{12}$.

\item[(4)]   $r=11$, $G\cong \M_{12}$, $T\cong \PSL_2(11)$, $G_\a\cong \M_{11}$, and $\Ga\cong \K_{12}$.
 \item[(5)] $r=23$, $G\cong\M_{24}$, $T\cong\PSL_3(2)$, $G_\a\cong \M_{23}$, and $\Ga\cong \K_{24}$.
\item[(6)] $r=23$, $G\cong\M_{24}$, $T\cong\PSL_2(23)$, $G_\a\cong \M_{23}$, and $\Ga\cong \K_{24}$.
\item[(7)] $n=24$, $r=7$, $G=\M_{24}$, $T=\M_{23}$, $G_\a\cong \PSL_3(2)$, and $\Ga$ is unique up to isomorphism.
\end{itemize}
\end{theorem}
\proof
By Lemmas \ref{Ga-primsb-A-n} and \ref{Ga=A-r} and Theorems  \ref{Ga-solvable}  and \ref{G-or-T-alte}, one of parts (1)-(3) holds if one of the following cases occurs: (i) $G_\a$ is primitive on $\Delta$; (ii) $G_\a$ is solvable; (iii) one of   $T$ and $G_\a$ is an alternating group. Thus we suppose that none of (i)-(iii) occurs.
By \cite[Theorem 1.1]{sps-factorization}, combining Lemma \ref{stab-lem-insol}, the following cases are left:
\begin{itemize}
\item[(iv)] $r=c-1$, $G\cong \M_{c}$, $G_\a\cong \M_{c-1}$, $c\in \{12,24\}$, and $T$ is a transitive   subgroup of $M_c$ acting on $c$ points.
\item[(v)]  $G=\M_{n}$, $T=\M_{n-1}$, $n\in \{12,24\}$, and $G_\a$ is a transitive subgroup of $M_n$.
\end{itemize}

 Assume that (iv) occurs. Then $\Ga\cong \K_c$ and, in this case, $n>c$ as $G\le \A_n$. By \cite{Gentchev}, $(c,T)$ is one of $(12,\PSL_2(11))$, $(24,\PSL_2(7))$ and $(24,\PSL_2(23))$, which gives one of parts (4)-(6) of this lemma, respectively.

 Assume that (v) occurs. Combining Lemma \ref{stab-lem-insol}, calculation using GAP shows that either $n=12$ and $G_\a\cong \Sy_5$, or $n=24$ and $G_\a$ is isomorphic to one of the following groups: $\Sy_5$ (twice up to the conjugacy in $\M_{24}$),
 $2^2\times\A_5$ (twice), $(2^2\times\A_5){:}2$ (twice),
 $\A_4\times\A_5$,  $2^4\times\A_5$, $(\A_4\times\A_5){:}2$,
 $2^6{:}\GL_2(4)$, $2^6{:}\GammaL_2(4)$, $\PSL_3(2)$, $2^6{:}(3\times \PSL_3(2))$ and $2^6{:}(3\times \PSL_3(2)){:}2$.
 Consider the normalizer $N$ of a Hall $r'$-subgroup of $G_\a$ in $\M_n$. Since $N$ contains a $2$-element which together with $G_\a$ generate $\M_n$, it follows from a GAP calculation
 that all cases but $G_\a\cong \PSL_3(2)$ are excluded.
 The group $\PSL_3(2)$ has two subgroups  of index $7$ up to conjugation, which are isomorphic to $\Sy_4$.
 Embedding $\PSL_3(2)$ in $\M_{24}$ as a transitive subgroup, one $\Sy_4$ of $\PSL_3(2)$ is self-normalized in $\M_{24}$, and the second one has normalizer $\ZZ_2\times \Sy_4$. Then we have a unique
 $\Ga$ arising from the second $\Sy_4$, and thus part (7) of this lemma follows.
 \qed


\vskip 10pt

By the argument in Section \ref{sect=sol-stab} and this section, our task is reduced to those graphs which satisfy
the following assumptions.

\begin{hypothesis}\label{hypo-graphs-classical}
 $\Gamma=(V,E)$, $G$, $T$ and $r$ satisfy Hypothesis \ref{hypo-graphs}, and
\begin{itemize}

\item[(IV)] $G$ is   a classical simple group over some finite field,
  $G_\a$ has a unique insolvable composition factor, say $S=\soc(G_\a^{\Ga(\a)})$, where $\a\in V$;
\item[(V)]  $G\not\cong \A_5,\,\A_6,\,\A_8$, $\PSL_m(q)$ (with $m$ a prime) or $\PSU_{2m+1}(q)$, and neither $T$ nor
$G_\a$ is isomorphic to an alternating group.

\end{itemize}
\end{hypothesis}

\vskip 20pt

\section{Some technique lemmas}\label{sect-tech}

The following two lemmas are derived from
\cite[Theorem 3.1]{Bam-Pen}, which bring   a lot of convenience for us to determine the composition factor $S$ of $G_\a$ in Hypothesis {\rm \ref{hypo-graphs-classical}}.


\begin{lemma}\label{linear-stab-1}
Let $S$ be a nonabelian simple group with a subgroup of prime index $r$.
Let $q=p^f$ for a prime $p$ and integer $f\ge 1$, and
let $H\le \GL_m(q)$ for some $m\ge 2$.
 Suppose that $H$ has a unique insolvable composition factor $S$.
If
$\Phi^*_{mf}(p)\ne 1$ and $|H|$ is divisible by $\Phi^*_{mf}(p)$ then one of the following holds.
\begin{itemize}
\item[(1)] $S=\PSL_d(q^{m\over d})$,  $|H|_p$ is a divisor of ${m\over d}|S|_p$, and either $(m,d,q)=(2,2,11)$ or $d$ is a prime and $mf$ is a power of $d$.
\item[(2)] $S$, $m$ and $q$ are listed in Table {\rm \ref{tab=exc-1}}.
\end{itemize}
\end{lemma}
\begin{table}[ht]
\[\tiny
\begin{array}{|l|l|l|l|l|l|l|l|l|l|}\hline
S &\A_5  & \A_7  & \A_{11} & \A_{13} & \A_{19} & \M_{11} &\PSL_3(2)&\PSL_2(11)&\PSL_3(3)  \\ \hline
(m,q)&(4,2)& (6,3)& (10,2)&(12,2)&(18,2)&(10,2)&(6,3)& (10,2)&(12,2) \\ 
& (4,3)& (6,5)&&&&& (6,5)&(5,4)&\\
&& (4,2)&&&&&(3,9)&&\\
&&(4,7)&&&&&(3,25)&&\\
&&(3,25)&&&&&&&\\ \hline
\end{array}
\]
{\small\caption{}\label{tab=exc-1}}
\end{table}
\proof
Assume that $\Phi^*_{mf}(p)\ne 1$ and $|H|$ has a divisor $\Phi^*_{mf}(p)$.
By  \cite[Theorem 3.1]{Bam-Pen}, ignoring those $H$ which does not produce an suitable $S$, one of the following holds:
\begin{itemize}
\item[(1.1)] $H$ has a normal subgroup $\SL_m(q)$, $\Sp_2(q)$ (with $m=2$)  or $\Omega^-_4(q)$ (with $m=4$);
\item[(1.2)] $H\le \GL_{m\over b}{.}[b]$, where $b$ is a proper divisor of $m$;
\item[(1.3)] $S$, $m$ and $q$ are listed in Table \ref{tab=exc-1}.
\end{itemize}
We only need to deal with cases (1.1) and (1.2).

Assume that (1.1) holds. Then either $S=\PSL_2(11)$, or $r$ is one of ${q^m-1\over q-1}$, $q+1$ and $q^2+1$ respectively.
By Lemma \ref{r}, for the latter case, $m$ is a prime and
$mf$ is a power of $m$, or $m=4$ and $2f$ is a power of $2$.
Moreover, it is easy to see that $|H|_p=|S|_p$.
Then part (1) of this lemma occurs.

Assume that case (1.2) occurs. Let $H_1=H\cap \GL_{m\over b}(q^b)$. Then $S\preccurlyeq H_1$, and $|H|_p$ is a divisor of $|H_1|_pb$. By \cite[pp. 38,Proposition B]{max-factor}, $\Phi^*_{mf}(p)$  is a divisor of $|H_1|$. Applying \cite[Theorem 3.1]{Bam-Pen}
to $H_1$ and $\Phi^*_{{m\over b}bf}(p)$, we have
\begin{itemize}
\item[(2.1)] $H_1$ has a normal subgroup $\SL_{m\over b}(q^b)$, $\Sp_2(q^b)$ (with $b={m\over 2}$)  or $\Omega^-_4(q^b)$ (with $b={m\over 4}$); or
\item[(2.2)] $H_1\le \GL_{m\over ab}{.}[a]$, where $a$ is a proper divisor of ${m\over b}$; or
\item[(2.3)] $(S, {m\over b},q^b)$ is one of $(\A_7,3,25)$, $(\PSL_3(2),3,9)$, $(\PSL_3(2),3,25)$, $(\PSL_2(11),5,4)$.
\end{itemize}
First, all triples in (2.3) are included in Table \ref{tab=exc-1}.

Assume that (2.1) occurs. Then $q^b\ne 11$ as $b>1$, and $|H_1|_p=|S|_p$. Then $|H|_p$ is a divisor of $b|S|_p$.
By Lemma \ref{r}, $S$ is described as in part (1) of this lemma. 

Assume that case (2.2) occurs. By \cite[Theorem 3.1]{Bam-Pen} and an induction on $m$, we conclude that $S=\PSL_d(p^{m\over d})$, $\PSp_2(q^{m\over 2})$ or $\POmega^-_4(q^{m\over 4})$. Then,  by Lemma \ref{r}, $S$ is described as in part (1) of this lemma.
\qed

\begin{lemma}\label{linear-stab-2}
Let $S$ be a nonabelian simple group with a subgroup of prime index $r$.
Let $q=p^f$ for a prime $p$ and integer $f\ge 1$, and let
 $H\le \GL_m(q)$ for some $m\ge 3$. Suppose that $H$ has a unique insolvable composition factor $S$.
 If
$\Phi^*_{(m-1)f}(p)\ne 1$ and $|H|$ is divisible by $\Phi^*_{(m-1)f}(p)$ then one of the following holds.
\begin{itemize}
\item[(1)] $S=\PSL_m(q)$,  $m$ is a prime and $f$ is a power of $m$.
\item[(2)] $S=\PSL_d(q^{m-1\over d})$, $|H|_p$ is a divisor of ${m-1\over d}q^{m-1}|S|_p$, and  either $(m-1,q)=(2,11)$ or $d$ is a prime and $(m-1)f$ is a power of $d$.
\item[(3)] $S$, $m$ and $q$ are listed in Table {\rm \ref{tab=exc-2}}.
\end{itemize}
\end{lemma}
\begin{table}[ht]
\[\tiny
\begin{array}{|l|l|l|l|l|l|l|l|l|l|l|}\hline
S &\A_5  & \A_7  & \A_{11} & \A_{13} & \A_{19} & \M_{11}&\M_{23} &\PSL_3(2)&\PSL_2(11)&\PSL_3(3)  \\ \hline
(m,q)&(5,2)& (7,3)& (11,2)&(13,2)&(19,2)&(5,3)& (11,2)& & &  \\ 
& (5,3)& (7,5)&&(11,2)&&& & & &\\
&& (5,2)&&&&&& &&\\
&&(5,3)&&&&&& &&\\
&&(4,2)&&&&&&&&\\
&&(4,9)&&&&&&&&\\
&&(4,25)&&&&&&&&\\ \hline

(m-1,q)&(4,2)& (6,3)& (10,2)&(12,2)&(18,2)&(10,2)&&(6,3)& (10,2)&(12,2) \\ 
& (4,3)& (6,5)&&&&& &(6,5)&(5,4)&\\
&& (4,2)&&&&&&(3,9)&&\\
&&(4,7)&&&&&&(3,25)&&\\
&&(3,25)&&&&&&&&\\
&&(3,2)&&&&&&&&\\ \hline
\end{array}
\]
{\caption{}\label{tab=exc-2}}
\end{table}
\proof
Assume that $\Phi^*_{(m-1)f}(p)\ne 1$ and $|H|$ has a divisor $\Phi^*_{(m-1)f}(p)$.
By  \cite[Theorem 3.1]{Bam-Pen}, one of the following holds:
\begin{itemize}
\item[(1.1)] $H$ has a normal subgroup $\SL_m(q)$;
\item[(1.2)] $H\le G_1:=q^{m-1}{:}((q-1)\times \GL_{m-1}(q))$;
\item[(1.3)] $S$, $m$ and $q$ are listed in Table \ref{tab=exc-2}.
\end{itemize}
Case (1.1) gives part (1) of this lemma occurs. If  (1.3) occurs then $(S,m,q)$ is given as in  part (3) of this lemma.
For case (1.2), by Lemma \ref{linear-stab-1}, either $S$ is described as in part (2) of this lemma, or $(S,m-1,q)$ is listed in Table \ref{tab=exc-2}.
\qed

\vskip 10pt

Under Hypothesis \ref{hypo-graphs-classical}, let $A$ and $B$ be maximal subgroups of $G$ such that one of them contains $T$ and the other one contains $G_\a$.
Then we have a maximal factorization
$G=AB$  with two factors $A$ and $B$ insolvable.


 Employing \cite[Theorem A]{max-factor},
we will work on the possible candidates for the triple $(G,T,G_\a)$.
 An essential part of our approach is determine the  inclusion of $G_\a$ in $A$ or $B$.
Suppose that $G_\a\le B$ and $B$ has an almost simple quotient $B/N=L.O$, where $L$ is simple.
Assume further that $|G:A|$ has a prime divisor which does not divide $|O||N|$. By Lemma \ref{sub-factor}, either
$S\preccurlyeq B$, or we have a core-free
maximal factorization of $B/N$ with a factor admitting  $S$ as a composition factor. Then, by  Lemma \ref{m-factor} and \cite[Theorem A]{max-factor}, we may determine $S$ in a  recursive process.
We   use some notation
defined in \cite[pp.5]{max-factor}, such as   $\sP_i,\,\sN_i,\,\sN_i^-,\,\sN_i^+$.

\begin{lemma}\label{linear-stab-3}
Let $S$ be a nonabelian simple group with a subgroup of prime index $r$.
Let $q=p^f$ for a prime $p$ and integer $f\ge 1$,  $L=\PSU_{2m}(q)$ and $\soc(G_0)=L$ for some $m\ge 2$.
Suppose that $G_0=A_0B_0$ is a core-free maximal factorization,  $A_0\cap L=\sN_1$ and $G_0=A_0Y$ for $Y\le B_0$.
Suppose that $S$ is the unique insolvable composition factor of $Y$, and $\pi(S)=\pi(Y)$.   Then   one of the following holds.
\begin{itemize}
\item[(1)] $S=\PSL_2(q^m)$, $p=2$, $r=q^m+1$ and $mf$ is a power of $2$.
\item[(2)]   $S$, $2m$ and $q$ are listed in Table {\rm \ref{tab=exc-3}}.
\end{itemize}
\begin{table}[ht]
\[\tiny
\begin{array}{|l|l|l|l|l|l|l|l|l|l|}\hline
S &\A_5  & \A_7  & \A_{11} & \A_{13} & \A_{19} & \M_{11} &\PSL_3(2)&\PSL_2(11)&\PSL_3(3)  \\ \hline
(2m,q)&&&&&&&(6,2)&&(6,2)\\ \hline

(2m,q)&(4,2)& (6,3)& (10,2)&(12,2)&(18,2)&(10,2)&(6,3)& (10,2)&(12,2) \\ 
& (4,3)& (6,5)&&&&& (6,5)& &\\
&&(4,7)&&&&& &&\\ \hline

(m,q^2)& & (3,25)&&&&& (3,25)&(5,4)&\\
& &  &&&&& (3,9)& &\\\hline
\end{array}
\]
{\caption{}\label{tab=exc-3}}
\end{table}
\end{lemma}
\proof
By \cite[Theorem A, Tables 1 and 3]{max-factor}, $|G_0:A_0|=q^{2m-1}{q^{2m}-1\over q+1}$ and $L\cap B_0$ is described as follows:
\[\tiny
\begin{array}{l|l|l|l|l|l|l}
L\cap B_0& \sP_m& \PSp_{2m}(q){.}a&  {}^{\bf\hat{}}\SL_m(4).2& {}^{\bf\hat{}}\SL_m(16).[6]&\PSU_4(3){.}2,\M_{22}& \Suz\\ \hline
& & a\le 2 & q=2,m\ge 3&q=4, G_0\ge L{.}4&q=2,m=3&q=2,m=6\\
\end{array}
\]
Noting that $Y^\infty\unlhd L \cap Y$, it follows that $S\preccurlyeq L\cap Y$.
Since $LY\le G_0$, it follows that $|Y:(L\cap Y)|$ is a divisor of $o:=|G_0:L|$.
By Lemma \ref{Order}, $|Y|$ is divisible by $|G_0:A_0|$.
Then $|L\cap Y|$ has a divisor $q^{2m-1}{q^{2m}-1\over o(q+1)}$.

Assume that $L\cap B_0\cong \Suz$. Then $|(L\cap B_0)\cap Y|$ is divisible by $2^{11}\cdot5\cdot7\cdot13$. By the Atlas \cite{Atlas}, we have $(L\cap B_0)\cap Y\lesssim \G_2(4)$. Again by  the Atlas \cite{Atlas}, $\G_2(4)$ has no maximal subgroup of order divisible $2^{11}\cdot5\cdot7\cdot13$. It follows that $(L\cap B_0)\cap Y=\G_2(4)$, and then $S\cong \G_2(4)$, which is impossible.

Assume that $(m,q)=(3,2)$. Then $|L\cap Y|$ is divisible by $2^{5} \cdot7$, and $L\cap B_0\cong 2^9{:}\PSL_3(4)$, $\Sp_6(2)$, $\PSL_3(4){.}2$, $\PSU_4(3){.}2$ or $\M_{22}$. Noting that $S\preccurlyeq (L\cap Y)\Rad(L\cap B_0)/\Rad(L\cap B_0)$,
with the help of GAP, we search the  subgroups of $(L\cap B_0)/\Rad(L\cap B_0)$ with order divisible by $2^{5} \cdot7\over (2^{5}, |\Rad(L\cap B_0)|)$. It follows that $S\cong \PSL_3(2)$ or $\PSL_3(3)$, desired as in part (2) of this lemma.

For the left candidates for $L\cap B_0$, we have $\Phi^*_{2m}(q)\ne 1$, and $|L\cap Y|$ has a divisor $q^{2m-1}\Phi^*_{2m}(q)$.
Let $\overline{B}=(L\cap B_0)/\Rad(L\cap B_0)$ and $\overline{Y}=(L\cap Y)\Rad(L\cap B_0)/\Rad(L\cap B_0)$.
Then $S\preccurlyeq \overline{Y}$, and $|\overline{Y}|$ is divisible by $\Phi^*_{2m}(q)$.
By \cite[pp. 38, Proposition B]{max-factor}, we conclude that  both $|\overline{B}^\infty|$ and $|\overline{B}^\infty\cap \overline{Y}|$ are divisible by $\Phi^*_{2m}(q)$.
Thus Lemma \ref{linear-stab-1} works for the preimages of
$\overline{B}^\infty\cap \overline{Y}$ and $\overline{B}^\infty$ in $\GL_m(q^2)$ (or in $\GL_{2m}(q)$ when $\overline{B}^\infty\cong \PSp_{2m}(q)$).
Then the lemma follows.
\qed

\begin{lemma}\label{linear-stab-4}
Let $S$ be a nonabelian simple group with a subgroup of prime index $r$.
Let $q=p^f$ for a prime $p$ and integer $f\ge 1$,  $L=\POmega^+_{8}(q)$ and $\soc(G_0)=L$.
Suppose that $G_0=A_0B_0$ is a core-free maximal factorization,  $A_0\cap L\cong \Omega_7(q)$ and $G_0=A_0Y$ for $Y\le B_0$.
Suppose that $S$ is the unique insolvable composition factor of $Y$, and $\pi(S)=\pi(Y)$.
 Then either $S\cong \PSL_2(q^2)$ and $r=q^2+1$ with $p=2$ and $f$   a power of $2$, or $(S,q)$ is one of  $(\A_5,3)$, $(\A_7,2)$, $(\A_7,3)$ and $(\A_7,7)$.
\end{lemma}
\proof
By \cite[Table 4]{max-factor},  $L\cap B_0$ is described as follows:
\[\tiny
\begin{array}{l|l|l|l|l|l}
L\cap B_0& (L\cap A_0)^\tau& \sP_i&  {}^{\bf\hat{}}({q\mp1\over e}{\times}\Omega^\pm_6(q)){.}2^e&(\PSp_2(q)\otimes\PSp_4(q)){.}2 &\Omega^-_8(\sqrt{q})  \\ \hline
&\mbox{triality }\tau &i\in \{1,3,4\} &e=(2,q-1)&q\mbox{ odd} &q\mbox{ square}\\ \hline
L\cap B_0& \A_9& (\PSL_2(4)\times\PSL_2(4)){.}2^2&(\PSL_2(16)\times\PSL_2(16)){.}2^2&\Omega^+_8(2)\\ \hline
&q=2& q=2&q=4, G_0\ge L{.}2&q=3
\end{array}
\]
It is easily shown that $S\preccurlyeq L\cap Y$, and $|L\cap Y|$ has a divisor ${1\over o }|G_0:A_0|$, where $o=|G_0:L|$. Note that  $|G_0:A_0|={1\over (2,q-1)}q^3(q^4-1)$.

Suppose that $q=3$ and $L\cap B_0\cong \Omega^+_8(2)$.
Then $|L\cap Y|$ has a divisor $2^3\cdot3^2\cdot5$.
Checking the orders of maximal subgroups of $\Omega^+_8(2)$,
it follows that $L\cap Y\lesssim\Sp_6(2)$, $2^6{:}\A_8$, $\A_9$, $(3\times\PSU_4(2)){:}2$ or $(\A_5\times\A_5){.}2^2$. Inspecting the subgroups of $\Sp_6(2)$, $\A_9$, $\A_8$, $\PSU_4(2)$ and  $\A_5\times\A_5$, we conclude that $S\cong  \A_5$ or $\A_7$.
If $L\cap B_0\cong \A_9,\,  (\PSL_2(4)\times\PSL_2(4)){.}2^2$ or $(\PSL_2(16)\times\PSL_2(16)){.}2^2$ then we have $(S,q)=(\A_7,2)$, $(\PSL_2(q^2),2)$ or $(\PSL_2(q^2),4)$.

Next we deal with the left factorizations $G_0=A_0B_0$.

{\bf Case 1}. Assume that $L\cap B_0$ is one of $\sP_i,\,  {}^{\bf\hat{}}({q\mp1\over 2}{\times}\Omega^\pm_6(q)){.}2^e,\,(\PSp_2(q)\otimes\PSp_4(q)){.}2$ and $\Omega^-_8(\sqrt{q})$.
Let $\overline{B}=(L\cap B_0)/\Rad(L\cap B_0)$ and $\overline{Y}=(L\cap Y)\Rad(L\cap B_0)/\Rad(L\cap B_0)$. Then $S\preccurlyeq \overline{Y}$, and $|\overline{Y}|$ is divisible by ${1\over a}|G_0:A_0|$, where $a$ is a divisor of $o|\Rad(L\cap B_0)|$.
Then $S\preccurlyeq \overline{B}^\infty\cap \overline{Y}$, and it is easily shown that $|\overline{B}^\infty\cap \overline{Y}|$ is divisible by $\Phi^*_{4}(q)$ (and by $\Phi^*_{8}(\sqrt{q})$ if $q$ is a square). If $\overline{B}^\infty$ is neither $\PSU_4(q)$ nor
$\PSp_2(q)\times\PSp_4(q)$ then, applying Lemma \ref{linear-stab-1} to the preiamge of $\overline{B}^\infty\cap \overline{Y}$ in $\GL_4(q)$ (or $\GL_8(\sqrt{q})$), we conclude that either $S\cong \PSL_2(q^2)$, $p=2$ and $f$ is a power of $2$, or $(S,q)$ is one of $(\A_5,2)$, $(\A_5,3)$, $(\A_7,2)$, $(\A_7,3)$ and $(\A_7,7)$.

Assume that $\overline{B}^\infty\cong \PSU_4(q)$.
Then $|\overline{B}^\infty\cap \overline{Y}|$ has a divisor ${1\over (2,p)}q^3\Phi^*_{4}(q)$. Checking   maximal subgroups
of $\PSU_4(q)$ (refer to \cite[Tables 8.10 and 8.11]{Low}), we conclude that $\overline{B}^\infty\cap \overline{Y}\le [q^4]{:}[a_1]{.}\PSL_2(q^2){.}[a_2]$, where $a_1a_2={(q-1)(2,q^2-1)\over (4,q+1)}$. Applying Lemma \ref{linear-stab-1} to a suitable quotient of $\overline{B}^\infty\cap \overline{Y}$, we have $S\cong \PSL_2(q^2)$, $p=2$ and $f$ is a power of $2$.

Now let $\overline{B}^\infty\cong \PSp_2(q)\times\PSp_4(q)$ with odd $q>3$. Then $|\overline{B}^\infty\cap \overline{Y}|$ is divisible by $q^3\Phi^*_{4}(q)$. Let $N$ be the direct   factor of $\overline{B}^\infty$ with $N\cong \PSp_2(q)$.
Noting that $S\preccurlyeq \overline{B}^\infty\cap \overline{Y}$,
we have $\pi(S)=\pi(\overline{B}^\infty\cap \overline{Y})$.
It follows that $S\not\preccurlyeq N$, and so $S\preccurlyeq (\overline{B}^\infty\cap \overline{Y})N/N$.
 Clearly, $(\overline{B}^\infty\cap \overline{Y})N/N$
 has order divisible by  $\Phi^*_{4}(q)$.
Since $q^2+1$ is not a prime as $q$ is odd, it follows from
Lemma \ref{linear-stab-1} that $q=7$ and $S\cong \A_7$.

{\bf Case 2}.
Assume that $L\cap A_0\cong\Omega_7(q)\cong L\cap B_0$. Then $B_0=A_0^\tau$ for some triality $\tau$.
By \cite[Proposition 3.1.1]{Kleidman-1},
$L\cap A_0\cap B_0\cong \G_2(q)$.
Considering  the limitations on $S$, we know that  $\G_2(q)\not\preccurlyeq Y$. Since $G_0=A_0Y$, we have $B_0=(A_0\cap B_0)Y$.
By Lemma \ref{m-factor}, we have a maximal core-free factorization $G_1=A_1B_1$ with $B_0^\infty\le G_1$, $(A_0\cap B_0)^\infty\le A_1$ and $Y^\infty\le B_1$.
 By \cite[Theorem A]{max-factor},   $B_1$ can be read out from \cite[Tables 2 and 3]{max-factor}:
\[\tiny
\begin{array}{l|l|l|l|l|l}
B_1\cap B_0^\infty& [q^5]{:}({q-1\over 2}{\times} \PSp_4(q)){.}2& \Omega^\pm_6(q){.}2& ({q\pm1\over 2}{\times}  \PSp_4(q)){.}2^2& \Sp_6(2), q=3&\Sy_9, q=3\\ \hline
B_1\cap B_0^\infty&\Omega^\pm_6(q){.}2&[q^5]{:}({q-1\over 2}{\times} \PSp_4(q))&\Sp_2(q){\times} \PSp_4(q)
\end{array}
 \]
For $\Sp_6(2)$ and $\Sy_9$,
we have $S\cong  \A_5$ or $\A_7$  by checking their subgroups.
For the left candidates for $B_1$, by a similar argument as in Case 1, either $S\cong \PSL_2(q^2)$, $p=2$ and $f$ is a power of $2$, or $(S,q)$ is one of $(\A_5,2)$, $(\A_5,3)$, $(\A_7,2)$, $(\A_7,3)$ and $(\A_7,7)$.
\qed

 \begin{lemma}\label{linear-stab-5}
 Let $S$ be a nonabelian simple group with a subgroup of prime index $r$.
Let $q=p^f$ for a prime $p$ and integer $f\ge 1$,  $L=\POmega^+_{2m}(q)$ and $\soc(G_0)=L$ for some $m\ge 5$.
Suppose that $G_0=A_0B_0$ is a core-free maximal factorization,  $A_0\cap L=\sN^-_2$ and $G_0=A_0Y$ for $Y\le B_0$.
Suppose that $S$ is the unique insolvable composition factor of $Y$, and $\pi(S)=\pi(Y)$.
Then   one of the following holds.
\begin{itemize}
\item[(1)] $S=\PSL_m(q)$, $r={q^m-1\over q-1}$,  $m$ is a prime and $f$ is a power of $m$.
\item[(2)] $(m,q)=(6,2)$, $S\cong \PSL_5(2)$ and $r=31$.
\end{itemize}
\end{lemma}
\proof
By \cite[Tables 1-3]{max-factor}, $|G_0:A_0|={q^{2m-2}(q^m-1)(q^{m-1}-1)\over 2(q+1)}$, and
$B_0$ is of type $\sP_m$, $\sP_{m-1}$, ${}^{\bf\hat{}}\GL_m(2){.}2$ or ${}^{\bf\hat{}}\GL_m(4){.}2$.
Let $\overline{B_0}=B_0'/\O_p(B_0')$ and $\overline{Y}=(Y\cap B')\O_p(B_0')/\O_p(B_0')$. Then $S\preccurlyeq \overline{Y}$ and
$|\overline{Y}|$ is divisible by $\Phi^*_m(q)\Phi^*_{m-1}(q)$.
Consider the preimage of $\overline{Y}$ in $\GL_m(q)$, and note that $\pi(S)=\pi(\overline{Y})$.
If $\Phi^*_m(q)=1$ then $(m,q)=(6,2)$, and thus $S\cong \PSL_5(2)$ by Lemma \ref{linear-stab-2}.
Assume that $\Phi^*_m(q)\ne 1$. Then $S$ is described as in
 Lemma \ref{linear-stab-1}. For those $S$ in Table \ref{tab=exc-1}, we  have $\Phi^*_{m-1}(q)\ne 1$ and $(|S|,\Phi^*_{m-1}(q))=1$. Thus $L\cong \PSL_d(q^{m\over d})$,
 where $d$ is a prime and $mf$ is a power of $d$.
 The only possibility is that $m=d$, and the lemma follows.
 \qed

 \begin{lemma}\label{linear-stab-6}
 Let $S$ be a nonabelian simple group with a subgroup of prime index $r$.
Let $q=p^f$ for a prime $p$ and integer $f\ge 1$,  $L=\POmega^+_{2m}(q)$ and $\soc(G_0)=L$ for some $m\ge 5$.
Suppose that $G_0=A_0B_0$ is a core-free maximal factorization,  $A_0\cap L=\sN_1$ and $G_0=A_0Y$ for $Y\le B_0$.
Suppose that $S$ is the unique insolvable composition factor of $Y$, and $\pi(S)=\pi(Y)$.
Then   one of the following holds.
\begin{itemize}
\item[(1)] $S=\PSL_d(q^{m\over d})$, $r={q^m-1\over q^{m\over d}-1}$,  $d$ is a prime and $mf$ is a power of $d$.
\item[(2)] $(S,m,q)$ is one of the following triples: \\ $(\A_7,6,3)$,  $(\A_7,6,5)$,  $(\A_{11},10,2)$, $(\A_{13},12,2)$, $(\A_{19},18,2)$, $(\M_{11},10,2)$, \\ $(\PSL_3(2),6,3)$, $(\PSL_3(2),6,5)$, $(\PSL_2(11),10,2)$, $(\PSL_2(11),5,4)$,  $(\PSL_3(3),12,2)$,\\
    $(\A_7,6,2)$, $(\PSL_3(2),6,2)$ and  $(\PSL_3(3),6,2)$.
\end{itemize}
\end{lemma}
\proof
By \cite[Tables 1-3]{max-factor},  $L\cap B_0$ is described as follows:
\[\tiny
\begin{array}{l|l|l|l|l|l|l}
L\cap B_0&  \sP_m \mbox{ or } \sP_{m-1}&
{}^{\bf\hat{}}\GU_m(q){.}2&
(\PSp_2(q){\times}\PSp_{m}(q)){.}c&
{}^{\bf\hat{}}\GL_m(q){.}2& \Omega^+_m(4).2^2 & \Omega^+_m(16).2^2 \\ \hline
& &m\mbox{ even}& q>2,m\mbox{ even},c=(2,q-1)&G_0\ge L.2\mbox{ if}&q=2 &q=4, m\mbox{ even} \\
&&&&m\mbox{ odd}&m\mbox{ even}&G_0\ge L.2\\ \hline
L\cap B_0& \Omega_9(q).a& &\Co_1\\ \hline
& m=8, a\le 2&& m=12, q=2

\end{array}
\]
Note that $|Y|$ is divisible by $|G_0:A_0|={1\over (2,q-1)}q^{m-1}(q^m-1)$,  $S\preccurlyeq L\cap Y$,
and $|L\cap Y|$ has a divisor ${1\over o }|G_0:A_0|$, where $o=|G_0:L|$.

Assume that $L\cap B_0\cong \Co_1$.  Checking the maximal subgroups of $\Co_1$ with order divisible by $2^{10}\cdot3^2\cdot 5\cdot7\cdot 13$, we have $L\cap Y\le B_1\cong 3{.}\Suz{.}2$ or $(\A_4{\times}\G_2(4)){.}2$. Clearly, $S\preccurlyeq (L\cap Y)\Rad(B_1)/\Rad(B_1)$.  Checking
the subgroups of $B_1/\Rad(B_1)$ with order divisible by $2^{8}\cdot3\cdot 5\cdot7\cdot 13$, we conclude that $(L\cap Y)\Rad(B_1)/\Rad(B_1)\succcurlyeq\G_2(4)$ or $\Suz$. This implies that $S\cong \G_2(4)$ or $\Suz$, a contradiction.

{\bf Case 1}. Assume that $L\cap B_0\cong \Omega_9(q).a$,
$\sP_m,\,\sP_{m-1},\,
(\PSp_2(q){\times}\PSp_{m}(q)){.}c$ or
${}^{\bf\hat{}}\GL_m(q){.}2$.
Let $B_1=(L\cap B_0)'$ and $Y_1=Y\cap B_1$.
Then $S\preccurlyeq Y_1$ and $|Y_1|$ is divisible by $\Phi^*_m(q)$.
If $L\cap B_0\cong \Omega_9(q).a$  then
Lemmas \ref{linear-stab-1} or \ref{linear-stab-2} work for $B_1$ and $Y_1$, and thus $S\cong \PSL_2(q^4)$ and $r=q^4+1$ with $p=2$ and $f$ a power of $2$, desired as in part (1) of this lemma with $d=p=2$ and $m=8$. (Note $\Omega_9(q)=\Sp_8(q)$ for even $q$.)
For the left cases,  choose $N\unlhd B_1$ such that $B_1/N$ is almost simple,  $S\preccurlyeq Y_1N/N$ and $|Y_1N/N|$ is  divisible by $\Phi^*_m(q)$.
Then
Lemma \ref{linear-stab-1} works for the preimages of $B_1/N$ and $Y_1N/N$ in $\GL_m(q)$. It follows that $S$ is desired as in part (1) or (2) of this lemma.

{\bf Case 2}. Assume that $L\cap B_0\cong \GU_m(q){.}2$ for even $m$. Noting that $B_0=(A_0\cap B_0)Y$,
by Lemma \ref{sub-factor}, we have a core-free factorization
 $\overline{B_0}=\overline{(A_0\cap B_0)}\,\overline{Y}$,
where
$\overline{B_0}=B_0/\Rad(B_0)$, $\overline{(A_0\cap B_0)}=(A_0\cap B_0)\Rad(B_0)/\Rad(B_0)$ and $\overline{Y}=Y\Rad(B_0)/\Rad(B_0))$.
It is easily shown that $|\overline{(A_0\cap B_0)}|$ is divisible by $\Phi^*_{2m-2}(q)$. By Lemma \ref{m-factor}, we may
choose a core-free maximal factorization $G_1=A_1B_1$ such that
$\soc(G_1)\cong \PSU_m(q)$, $|A_1|$ is divisible by $\Phi^*_{2m-2}(q)$, $\overline{Y}^\infty\le B_1$ and $G_1=A_1(B_1\cap \overline{Y})$. By \cite[Theorem A]{max-factor},
either $A_1$ is a subgroup of $G_1$ with type $\sN_1$, or
$(m,q)=(6,2)$, $A_1\cap \soc(G_1)\cong \M_{22}$ and  $B_1\cap \soc(G_1)\cong \PSU_5(2)$. By Lemma \ref{linear-stab-3}, we conclude that $S$ is desired as in part (1) (with $d=p=2$) or (2) of this lemma.

Assume that $(m,q)=(6,2)$, $A_1\cap \soc(G_1)\cong \M_{22}$ and  $B_1\cap \soc(G_1)\cong \PSU_5(2)$. In this case, $|G_1:A_1|=2^8\cdot3^4$. It is easily shown that $S\preccurlyeq \overline{Y}\cap B_1\cap \soc(G_1)$, and $|\overline{Y}\cap B_1\cap \soc(G_1)|$ is divisible by $2^7\cdot3^3$.
Checking the maximal subgroups of $\PSU_5(2)$, we conclude that
$\overline{Y}\cap B_1\cap \soc(G_1)$ is solvable, a contradiction.

{\bf Case 3}. Assume that $L\cap B_0\cong  \Omega^+_m(q^2).2^2$ for $q\in \{2,4\}$ and even $m$.
Suppose that $m=6$. Then $(L\cap B_0)'\cong \PSL_4(q^2)$.
It is easily shown that $S\preccurlyeq Y\cap (L\cap B_0)'$, and
$|Y\cap (L\cap B_0)'|$ is divisible by ${1\over 8}q^5(q^6-1)$.
Assume that $q=4$. Then $S$ is given in Lemma \ref{linear-stab-2}, and the only possibility is that
 $S\cong \PSL_3(16)$. Then $r={16^3-1\over 16-1}$, which is not a prime. Thus $q\ne 4$, and so $q=2$. Checking the subgroups of maximal subgroups of $\PSL_4(4)$, we get
  $S\cong \PSL_3(2)$ or $\A_7$, which is included in part (2) of this lemma.

Next let $m\ge 8$. Similarly as in Case 2, we have a core-free maximal factorization $G_1=A_1B_1$ such that
$\soc(G_1)\cong \Omega^+_m(q^2)$, $|A_1|$ is divisible by $\Phi^*_{m-2}(q^2)$, $S\preccurlyeq \overline{Y}^\infty\le B_1$ and $G_1=A_1(B_1\cap \overline{Y})$, where $\overline{Y}=Y\Rad(B_0)/\Rad(B_0))$.

Assume thta $m=8$. If $A_1$ has a normal subgroup isomorphic to
$\Sp_6(q)$ then   $S\cong \PSL_2(q^4)$ by Lemma \ref{linear-stab-4}, desired as in part (1) of this lemma.
Thus, by \cite[Table 4]{max-factor}, we may let
$|G_1:A_1|={1\over2}q^{12}(q^6-1)(q^4+1)(q^2-1)$, and
$B_1\cap \soc(G_1)\cong [q^{12}]{:}\GL_4(q^2)$ or $\GL_4(4){.}2$.
Using Lemmas \ref{linear-stab-1} and \ref{linear-stab-2}, by a similar argument as in the proof of Lemma \ref{linear-stab-5}, we conclude that $S\cong \PSL_2(q^4)$ or $\PSL_3(q^2)$. Noting that
$\pi(S)=\pi(B_1\cap \overline{Y})$ and $|B_1\cap \overline{Y}|$ is divisible by $|G_1:A_1|$, it follows that $q=2$ and $S\cong \PSL_2(16)$, desired as in part (1) of this lemma.

Let $m\ge 10$. Then, by \cite[Theorem A and Proposition 2.5]{max-factor}, $A_1$ has type $\sN_1$, $\sN^-_2$ or ${}^{\bf\hat{}}\GU_{m\over 2}(q^2){.}2$ (with ${m\over 2}$ even). For    $\sN^-_2$, by Lemma \ref{linear-stab-5}, $S\cong \PSL_{m\over 2}(q^2)$ with
 ${m\over 2}$ a prime and $mf$ a power of ${m\over 2}$,
 yielding $m=2$, a contradiction.
 Assume that $A_1\cap \soc(G_1)\cong {}^{\bf\hat{}}\GU_{m\over 2}(q^2){.}2$.  By \cite[Table 1]{max-factor}, $B_1$ has type $\sP_1$ or $\sN^+_2$. It follows that $(|B_1|,\Phi^*_{m}(q))=1$.
 Clearly, $\Phi^*_{m}(q)\ne 1$.
 Recalling that $\pi(S)=\pi(Y)$ and $|Y|$ is divisible by $\Phi^*_{m}(q)$, since $S\preccurlyeq B_1\cap \overline{Y}$, it follows that $(|B_1|,\Phi^*_{m}(q))\ne 1$, a contradiction.
Thus  $A_1$ is of type $\sN_1$.
Then, by an  induction on $m$, the lemma follows.
\qed

 \begin{lemma}\label{linear-stab-7}
 Let $S$ be a nonabelian simple group with a subgroup of prime index $r$.
Let $q=2^f$,  $L=\Sp_{2m}(q)$ and $\soc(G_0)=L$, where $m\ge 3$ and $mf>3$.
Suppose that $G_0=A_0B_0$ is a core-free maximal factorization,  $A_0\cap L\cong \GO_{2m}^-(q)$ and $G_0=A_0Y$ for $Y\le B_0$.
Suppose that $S$ is the unique insolvable composition factor of $Y$, and $\pi(S)=\pi(Y)$.
Then  $B_0$ is given as in \cite[Table 1]{max-factor},
and one of the following holds.
\begin{itemize}
\item[(1)]   $S\cong \PSL_d(q^{m\over d})$, where $d$ is a prime and $mf$ is a power of $d$.

\item[(2)]  $(S,m,q)$ is one of  $(\A_7,6,2)$,   $(\PSL_3(2),6,2)$  and  $(\A_{13},12,2)$.
    \item[(3)] $B_0=\sP_m$,  $(S,m,q)$ is one of $(\A_7,4,2)$ and  $(\PSL_5(2),6,2)$.
\end{itemize}
\end{lemma}
\proof
Assume that $B_0$ is described as in \cite[Table 3]{max-factor}.
Then $(m,q)=(4,2)$, $G_0=\Sp_8(2)$, $B_0\cong \Sy_{10}$ and $A_0\cap B_0\cong \Sy_3\times \Sy_7$, see \cite[pp.97, 5.1.9]{max-factor}. By \cite[Theorem D]{max-factor}, $Y$ is a primitive group of degree $10$. Noting that $\A_{10}\not\preccurlyeq Y$, it follows that $S=\soc(Y)\cong \A_6$ or $\A_5$. Noting the limitations on $S$, we have $Y\cong \A_5$ or $\Sy_5$. Then $Y\cong \Sy_5$ as $|Y|$ is divisible by $|G_0:A_0|=2^3\cdot 3\cdot 5$, and part (1) of this lemma follows.

Assume that $B_0$ is described as in \cite[Table 2]{max-factor}. Then $L\cap B_0\cong \G_2(q)$.
Noting that $B_0=(A_0\cap B_0)Y$, by Lemma \ref{sub-factor}, we have a core-free factorization $B_1=X_1Y_1$ with $\soc(B_1)\cong \G_2(q)$, and $S\preccurlyeq Y_1$. It follows from
 \cite{Factor-excep} that $q=4$ and $S\cong \J_2$, $\PSU_3(4)$ or  $\PSU_3(3)$, which is impossible.

Assume next that $B_0$ is given as in \cite[Table 1]{max-factor}.
We discuss in several cases according the structure of $B_0$.
Note that $|Y|$ is divisible by $|G_0:A_0|={1\over 2}q^{m}(q^m-1)$.

{\bf Case 1}.
Let $B_0=\sP_m$, $\overline{B_0}=B_0/\O_2(B_0)$ and $\overline{Y}=Y\O_2(B_0)/\O_2(B_0)$.
Then $\overline{B_0}\cong \GL_m(q)$, $S\preccurlyeq\overline{Y}$ and  $|\overline{Y}|$ is divisible by $q^m-1$.
Assume that $(m,q)\ne (6,2)$. Then, by Lemma \ref{linear-stab-1},  either part (1) of this lemma occurs, or
$(S,m,q)$ is one of  $(\A_5,4,2)$,  $(\A_7,4,2)$, $(\A_{11},10,2)$, $(\A_{13},12,2)$,  $(\A_{19},18,2)$, $(\M_{11},10,2)$,  $(\PSL_2(11),10,2)$,  $(\PSL_2(11),5,4)$ and  $(\PSL_3(3),12,2)$. For the latter case, considering $\pi(S)=\pi(Y)$, only $(\A_5,4,2)$, $(\A_7,4,2)$ and $(\A_{13},12,2)$ are left. Note  the triple $(\A_5,4,2)$ is included in part (1) of this lemma.
Now let $(m,q)=(6,2)$. By Lemma  \ref{m-factor}, there is a factorization
of $\SL_6(2)$ with a factor having a composition factor $S$. Checking the subgroups of $\SL_6(2)$, we have
$S\cong \A_7$, $\SL_3(2)$ or $\SL_5(2)$, as in part (2) or (3) of this lemma.

{\bf Case 2}.  Let $B_0=\Sp_m(q)\wr\Sy_2$.
Note that $A_0\cap B_0\cong (\GO^-_m(q)\times \GO^+_m(q))$, see \cite[pp. 50, 3.2.4(b)]{max-factor}.
We have $A_0\cap B_0\le B_0'$ and $B_0'=(A_0\cap B_0)(Y\cap B_0')$. Let $N$ be a direct factor of
$B_0'$ such that $N$ contains the direct factor $\GO^-_m(q)$ of $A_0\cap B_0$. Let $\overline{B_0}=B_0'/N$, $\overline{Y}=(Y\cap B_0')N/N$ and $\overline{H}=(A_0\cap B_0)N/N$.
Then $\overline{B}\cong \Sp_m(q)$, $\overline{H}\cong \GO^+_m(q)$ and $\overline{B}=\overline{H}\,\overline{Y}$, and thus
$|\overline{Y}|$ is divisible by $|\overline{B}:\overline{H}|={1\over 2}q^{m\over 2}(q^{m\over 2}+1)$.

Assume that  $\overline{Y}$ is insolvable. Then $S\preccurlyeq\overline{Y}$.
Applying Lemma \ref{linear-stab-1} to the pair $(\overline{B},\overline{Y})$, by a similar argument as in Case 1, we have that either  $S\cong \PSL_2(q^{m\over 2})$ with $mf$ a power of $2$, or $(S,m,q)$ is one of $(\A_5,4,2)$,  $(\A_7,6,2)$, $(\PSL_3(2),6,2)$, $(\A_7,6,2)$ and $(\A_{13},12,2)$.
Assume that  $\overline{Y}$ is solvable. It follows from \cite[Theorem 1.1]{Li-Xia} that
$(m,q)=(4,2)$ or $(6,2)$, and $N\cong \Sy_6$ or $\Sp_6(2){.}2$, respectively. Noting that $\overline{Y}=(Y\cap B_0')N/N\cong (Y\cap B_0')/(Y\cap  N)$, we have
$S\preccurlyeq Y\cap N$. Since $\pi(Y\cap N)\subseteq\pi(S)=\pi(Y)$, checking the subgroups of $N$, we have that
$(S,m,q)$ is one of $(\A_5,4,2)$,  $(\A_7,6,2)$ and  $(\PSL_3(2),6,2)$.
Then our lemma holds in this case.

{\bf Case 3}.  Let $B_0\cong \GO^+_{2m}(2)$, $m\ge 4$ and $q=2$.
Then $A_0\cap B_0=\Sp_{2m-2}(2)\times 2$ is maximal in both $A_0$ and $B_0$, see \cite[pp. 52, 3.2.4(e)]{max-factor}.
Noting that $B_0=(A_0\cap B_0)Y_0$ is a core-free factorization, by Lemma \ref{m-factor}, choose a core-free maximal factorization $G_1=A_1B_1$ and a core-free factorization $G_1=A_1Y_1$ such that
$\soc(G_1)\cong \Omega^+_{2m}(2)$, $A_1\gtrsim \Sp_{2m-2}(2)$ and $S\preccurlyeq Y_1\le B_1$.
Then $S$ is described as in  Lemma \ref{linear-stab-6}. Recalling $|G_\a|$ is divisible by $2^{m-1}(2^m-1)$, except for $S\cong \PSL_3(2)$ and $m=6$, all
the triples in Lemma \ref{linear-stab-6}~(2) are easily excluded.
Thus $S\cong \PSL_2(2^{m\over 2})$, or $m=6$ and $S\cong \PSL_3(2)$ or $\A_7$, desired as  part (1) or (2) of this lemma.

{\bf Case 4}. Let $B_0\cong \Sp_{2a}(q^b){.}b$, where $b$ is a prime and $m=ab$. In this case, $A_0\cap B_0\cong \GO^-_{2a}(q^b){.}b$, which is maximal in $A_0$ (if $a\ne 1$) and $B_0$.
Noting that $B_0=(A_0\cap B_0)Y$, if $a\ge 3$ then part (1) of this lemma occurs by an induction on $m$.

Assume that $a=1$. Then $m$ is an odd prime and $B_0\cong \PSL_{2}(q^m){.}m$. Checking the subgroups of $\PSL_{2}(q^m)$ with order divisible by $q^m-1$, we conclude that $S\cong\PSL_{2}(q^m)$. Then $r=q^m+1$, which is not a prime, a contradiction. Thus $a=2$, ${m\over 2}$ is a prime and $B_0\cong \Sp_{4}(q^{m\over 2}){.}{m\over 2}$. By Lemmas \ref{m-factor} and \ref{sub-factor}, since $B_0=(A_0\cap B_0)Y$,
choose a core-free maximal factorization $G_1=A_1B_1$ and a core-free factorization $G_1=A_1Y_1$ such that
$\soc(G_1)\cong \Sp_{4}(q^{m\over 2})$, $A_1\gtrsim \Omega^-_{4}(q^{m\over 2})$ and $S\preccurlyeq Y_1\le B_1$.
By \cite[Theorem A]{max-factor}, we conclude that $B_1^\infty\Rad(B_1)/\Rad(B_1)\cong \PSL_2(q^m)$, $\PSL_2(q^{m\over 2})$ or $\PSL_2(q^{m\over 2})\times \PSL_2(q^{m\over 2})$. It follows that $S\cong \PSL_2(q^m)$ or $\PSL_2(q^{m\over 2})$. This implies that $r=q^m+1$ or $q^{m\over 2}+1$, and so $r$ is not a prime, a contradiction.
\qed

\vskip 20pt

\section{The classical case}\label{sect=classical}

In this section, we will work under Hypothesis \ref{hypo-graphs-classical}. Let $A$ and $B$ be maximal subgroups of $G$ such that one of them contains $T$ and the other one contains $G_\a$.
Then
$G=AB$, and both    $A$ and $B$ are insolvable. Without loss of generality, we let $A$ and $B$ be described as in the third and the fourth columns of Tables 1-4 in \cite{max-factor}, respectively. For convenience, we define $X$ and $Y$ as follows:
\[\{X,Y\}=\{T,G_\a\},\, X\le A \mbox{ and } Y\le B.\]
Recall that $S=\soc(G_\a^{\Ga(\a)})$.

\subsection{Unitary case}

Let $G=\PSU_{2m}(q)$, where $m\ge 2$  and $q=p^f$ for some  prime $p$. Then $A$ and $B$ are given
in  \cite[Tables 1 and 3]{max-factor}.

\begin{lemma}\label{U-2m-1}
 $A$ and $B$ are described as
in  \cite[Table 1]{max-factor}.
\end{lemma}
\proof
Suppose that $A$, $B$ are given as
in  \cite[Table 3]{max-factor}, and we list $A$ and $B$ as follows:
\[\tiny
\begin{array}{l|l|l|l|l}
(m,q)&A&|G:A|&B&|G:B|\\ \hline
(2,3)&\PSL_3(4)&2\cdot 3^4&\PSp_4(3)&2\cdot 3^2\cdot 7\\ \hline
(3,2)&\PSU_5(2)&2^5\cdot3\cdot7&\PSU_4(3).2&2^7\cdot 11\\
     &\PSU_5(2)&2^5\cdot3\cdot7&\M_{22}&2^8\cdot 3^4\\ \hline
(6,2)&\Suz&\Phi_{22}^*(2)\Phi_{18}^*(2)\mid |G:A|&\PSU_{11}(2)&2^{11}\cdot3\cdot5\cdot7\cdot13
\end{array}
\]
Inspect the insolvable groups contained in $A$ and $B$.
For $(m,q)=(2,3)$, combining  Lemma \ref{Order}, we conclude that
$\{T,G_\a\}=\{A,B\}$; however, neither $\PSL_3(4)$ nor $\PSp_4(3)$ has a subgroup of   prime index. Similarly, the pair $(3,2)$ is excluded for $(m,q)$.
Then only $(m,q)=(6,2)$ is left.
By \cite[Table 10.1]{transubgroup}, $\PSU_{11}(2)$ has no proper
subgroup with order divisible by $\Phi_{22}^*(2)\Phi_{18}^*(2)$.
By the Atlas \cite{Atlas}, $\Suz$ has no maximal subgroup with order divisible by $2^{11}\cdot3\cdot5\cdot7\cdot13$. It follows that $\{T,G_\a\}=\{A,B\}$. By Lemma \ref{stab-lem-insol}, none of $\PSU_{11}(2)$ and $\Suz$ is a possible vertex-stabilizer of some symmetric graph of prime degree, a contradiction. This completes the proof.
\qed

\vskip 10pt

In the following we let $A=\sN_1$. Then
$B$ and  $|G:B|$ are listed as follows:
\[\tiny
\begin{array}{l|l|l|l}
B& \sP_m& \PSp_{2m}(q){.}{(2,q-1)(m,q+1)\over (2m,q+1)}&  {}^{\bf\hat{}}\SL_m(4).2, q=2\\ \hline
|G:B|&\prod_{i=1}^m(q^{2i-1}+1)& {1\over (m,q+1)}q^{m(m-1)}\prod_{i=2}^m(q^{2i-1}+1)&
{1\over 2}q^{m^2}\prod_{i=1}^m(q^{2i-1}+1)\\
\end{array}
\]
By Lemma \ref{Order}, $|X|$ is divisible by $|G:B|$.
Since $A$ is insolvable, $(m,q)\ne (2,2)$.

\begin{lemma}\label{U-small-(m,q)}
If $(m,q)=(3,2)$ then $T\cong \PSU_5(2)$, $\SL_3(2)\unlhd G_\a^{\Ga(\a)}$ and $\O_2(G_\a)\ne 1$.
\end{lemma}
\proof
Let $(m,q)=(3,2)$.
Then $A\cong \PSU_5(2)$, $|G:A|={2^5\cdot 3\cdot 7}$, and $B$ is listed as follows:
\[\tiny
\begin{array}{l|l|l|l}
B&2^9{:}\PSL_3(4)& \PSp_{6}(2)&  \PSL_3(4).2\\ \hline
|G:B|&3^4\cdot 11& 2^6\cdot3^2\cdot11&
2^8\cdot3^4\cdot 11
\end{array}
\]
Searching the subgroups of $\PSU_5(2)$ with order divisible by $|G:B|$, we get $X=A\cong \PSU_5(2)$. Then $T=X\le A$, and so $G_\a=Y\le B$. Thus $|G_\a|$ is divisible $|G:A|$.
For  $B=\PSp_6(2)$, by calculation using GAP, we have $G_\a\cong 2^3{:}\SL_3(2)$ or $2^6{:}\PSL_3(2)$.

Assume that $B\cong \PSL_3(4).2$ or $2^9{:}\PSL_3(4)$.
Let $\overline{B}=B'/\Rad(B')$ and $\overline{Y}=
(Y\cap B')\Rad(B')/\Rad(B')$. Then $\overline{B}\cong \PSL_3(4)$, $S\preccurlyeq \overline{Y}$, and $|\overline{Y}|$ has a divisor $21$. Checking the subgroups of $\PSL_3(4)$, we get $S\cong \overline{Y}\cong \PSL_3(2)$. By Lemma \ref{stab-lem-insol},
$|G_\a|$  is a divisor of $6|\O_2(G_\a)||\GL_3(2)|$. Then $2^5\le |G_\a|_2\le 2^4|\O_2(G_\a)|$, yielding $\O_2(G_\a)\ne 1$, and the lemma follows.
\qed

\begin{lemma}\label{U-T}
 $A\ge T\cong \PSU_{2m-1}(q)$.
\end{lemma}
\proof
By Lemma \ref{U-small-(m,q)}, our lemma holds if $(m,q)=(3,2)$.
Thus assume $(m,q)\ne(3,2)$. Let $\overline{A}=A/\Rad(A)$ and
$\overline{X}=X\Rad(A)/\Rad(A)$. Then $\overline{A}\cong \PGU_{2m-1}(q)$, and $|\overline{X}|$ is divisible by $\Phi_{2(2m-1)}^*(q)$. By Lemma \ref{prime-div}, $\Phi_{2(2m-1)}^*(q)\ne 1$ as $(2(2m-1),q)\ne (6,2)$ and $2(2m-1)\ge 6$.

Assume that $m\ge 3$. Then  $\Phi_{2(2m-3)}^*(q)\ne 1$, and it is easily shown that $\overline{X}\cap \soc(\overline{A})$ has order divisible by $\Phi_{2(2m-1)}^*(q)\Phi_{2(2m-3)}^*(q)$. By \cite[Theorem 4 and Table 10.1]{transubgroup}, we conclude that $\overline{X}\cap \soc(\overline{A})\cong \PSU_{2m-1}(q)$.
This implies that $T=X\cong \PSU_{2m-1}(q)$.

Let $m=2$. Then $A=[{q+1\over (4,q+1)}]{.}\PGU_3(q)$, and $|G:A|={q^{3}(q-1)(q^2+1)}$. Suppose that $\overline{X}\cap \soc(\overline{A})\not\cong \PSU_3(q)$. Let $M$ be a maximal subgroup of $\soc(\overline{A})$ with $\overline{X}\cap \soc(\overline{A})\le M$.
Noting that $\Phi_{6}^*(q)\ge 7$, $|M|$ is divisible by $\Phi_{6}^*(q)$ and  $M$ is insolvable, by \cite[Tables 8.5 and 8.6]{Low}, $q=p$ and either $M\cong \PSL_3(2)$ with $p\ne 5$  or  $M\cong  \A_7$ with $p=5$.
In both cases,  $\Phi_{6}^*(q)=7$, and so $q\in \{3,5\}$ by \cite[Theorem 3.9]{Hering}.
Assume that $q=3$. Then $A\cong \PSU_3(3)$, $X\cong \PSL_3(2)$, and $|G:X|=2^4\cdot3^5\cdot5$. Moreover, $B\cong 3^4{:}\A_6$ or $\PSp_4(3)$. Since $G=XY$, by Lemma \ref{Order}, $|Y|$ is divisible by $2^4\cdot3^5\cdot5$. It follows that $B\cong 3^4{:}\A_6$ and $Y\cong 3^4{:}\A_5$. This forces $G_\a=Y\cong 3^4{:}\A_5$, which is impossible by Lemma \ref{stab-lem-insol}.
Assume that $q=5$. Then $A\cong 3{.}\PSU_3(5){.}3$, $X\lesssim
3{.}\A_7{.}3$, and $|G:X|$ is divisible by $2^4\cdot5^5\cdot13$.
In this case, $|Y|$ is  divisible by $2^4\cdot5^5\cdot13$.
It follows that $B\cong 5^4{:}\PSL_2(25){:}[4]$,
and then $Y\cong 5^4{:}\PSL_2(25){:}[c]$, where $c$ is a divisor of $4$. This implies that $G_\a=Y$, and so $r=25+1$, which contradicts that $r$ is a prime. Therefore,
$\overline{X}\cap \soc(\overline{A})\cong \PSU_3(q)$. Then $T=X\cong \PSU_3(q)$, and the lemma holds.
\qed

\begin{theorem}\label{U-2m}
Assume that  {\rm Hypothesis \ref{hypo-graphs-classical}} holds, and let $G=\PSU_{2m}(q)$, where $m\ge 2$, $(m,q)\ne (2,2)$  and $q=p^f$ for some  prime $p$.
Then  $T\cong \PSU_{2m-1}(q)$, $p=2$, $\O_2(G_\a)\ne 1$ and
one of the following holds.
\begin{itemize}
\item[(1)]  $(m,q)=(3,2)$ and $\SL_3(2)\unlhd G_\a^{\Ga(\a)}$.

\item[(2)]
$mf$ is a power of $2$   and $\SL_2(q^m)\unlhd G_\a^{\Ga(\a)}$.

\end{itemize}
\end{theorem}
\proof
By Lemmas \ref{U-2m-1} and \ref{U-T}, $\PSU_{2m-1}\cong T\le A=\sN_1$, and so $G_\a=Y\le B$.
Thus $S=\soc(G_\a^{\Ga(\a)})$ is described as in Lemma \ref{linear-stab-3}.
If $(m,q)=(3,2)$ then part (1)  holds by Lemma \ref{U-small-(m,q)}.
Next let $m>3$ if $q=2$. Note that $|G:T|={(2m-1,q+1)\over (2m,q+1)}q^{2m-1}(q^{2m}-1)$.

Since $|G_\a|$ is divisible by $|G:T|$, we have $|G_\a|_p\ge |G:T|_p$.
Suppose that $S$ is given in Table \ref{tab=exc-3}. Then $(S,2m,q)$ is one of $(\A_5,4,3)$, $(\A_7,4,7)$, $(\A_7,6,3)$, $(\A_7,6,5)$, $(\A_{11},10,2)$, $(\A_{13},12,2)$, $(\A_{19},18,2)$, $(\M_{11},10,2)$, $(\PSL_3(2),6,3)$,  $(\PSL_3(2),6,5)$,  $(\PSL_2(11),10,2)$ and $(\PSL_3(3),12,2)$.
For these triples, by Lemma \ref{stab-lem-insol}, calculation shows that  $|G_\a|_p<|G:T|_p$, a contradiction. Therefore,
$S\cong \PSL_2(q^m)$, $r=q^{m}+1$, $p=2$ and $mf$ is a power of $2$. Suppose that $\O_p(G_\a)=1$. Then, by Lemma \ref{stab-lem-insol}, $|G_\a|_2$ is a divisor of
$mfq^m$. Then $q^{2m-1}=|G:T|_p\le |G_\a|_2|\le mfq^m$, yielding $(m,q)=(2,2)$, a contradiction. Thus $\O_p(G_\a)\ne 1$, and then our result follows.
\qed

\subsection{Orthogonal case I}
Let  $G=\Omega_{2m+1}(q)$, where $m\ge 3$  and $q=p^f$ for some odd prime $p$.
Then $A$ and $B$ are given
in  \cite[Tables 1-3]{max-factor}. Clearly, $\Phi_{2m}^*(q)\ne 1\ne \Phi_{2m-2}^*(q)$.

\begin{lemma}\label{O-2m+1-1}
$A$ and $B$ are described  as
in  \cite[Table 1]{max-factor}.
\end{lemma}
\proof
Suppose that the lemma is false. We  deduce the contradiction in two cases.

{\bf Case 1}. Let $A$ and $B$ be as
in  \cite[Table 2]{max-factor}. We have the following three cases.

{\it 1.1}. $p=3$, $m=12$, $A\cong \F_4(q)$ and $B\cong  \Omega_{24}^-(q){.}2$. In this case, $|G:A|={q^{120}(q^4-1)(q^{10}-1)\over2}\prod_{i=7}^m(q^{2i}-1)$, and $|G:B|={1\over 2}q^{12}(q^{12}-1)$.  Clearly, $\Phi_{2m-4}^*(q)\ne 1$.

Recall that $X\le A$ and $Y\le B$. Then $|Y|$ is divisible by $\Phi_{2m}^*(q)\Phi_{2m-2}^*(q)\Phi_{2m-4}^*(q)$.    By \cite[Theorem 4 and Table 10.1]{transubgroup}, since $q$ is odd, $\POmega^-_{2m}(q)$ has no proper subgroup
with order divisible by $\Phi_{2m}^*(q)\Phi_{2m-2}^*(q)\Phi_{2m-4}^*(q)$.
It follows that $Y\succcurlyeq\POmega^-_{2m}(q)$.
Then $Y\ne G_\a$ by  Lemma \ref{stab-lem-insol}. Therefore,   $T=Y$ and $G_\a=X \le A$. Noting that  $G=BG_\a$, we have $A=(A\cap B)G_\a$. By Lemma \ref{stab-lem-insol}, $G_\a\ne A$. Then we have a factorization of $\F_4(q)$,
which is impossible by \cite[Theorem 1]{Factor-excep} as $q$ is odd here.

{\it 1.2}. $p=3$, $m=6$, $A\cong \PSp_6(q).a$ and $B\cong  \Omega_{12}^-(q){.}2$, where $a\le 2$. In this case, $|G:A|={1\over2a}q^{27}\prod_{i=4}^m(q^{2i}-1)$, and $|G:B|={1\over 2}q^{6}(q^{6}-1)$.

Similarly as in {\it 1.1}, we have $T\le B$, $G_\a\le A$, and so $|G_\a|$ is divisible $q^6(q^6-1)$.
By Lemma \ref{stab-lem-insol}, $S\not\cong  \PSp_6(q)$. It is easily shown that
$S\preccurlyeq G_\a\cap A'$, and $|G_\a\cap \A'|$ is divisible by $q^{6}\Phi_{6}^*(q)$, where  $A'$ is the derived subgroup of $A$.
Clearly, $\Phi_{6}^*(q)\ne 1$. Noting that $q$ is odd, it follows from   Lemma \ref{linear-stab-1}  that
$q\in \{3,5\}$ and $S\cong \A_7$ or $\PSL_3(2)$.  By  Lemma \ref{stab-lem-insol},
$|G_\a|$ is indivisible by $3^3$; however, here,  $|G_\a|$ is divisible by $q^6$, a contradiction.

{\it 1.3}. $m=3$, $A\cong \G_2(q)$ and $B$ is   one of $[q^5]{:}({q-1\over 2}\times \Omega_5(q)).2$, $\Omega_6^+(q){.}2$, $\Omega_6^-(q){.}2$, $({q-1\over 2}\times \Omega_5(q)).2^2$ (with $q>3$) and $({q+1\over 2}\times \Omega_5(q)).2^2$.

Suppose that $X=G_\a$. Then $G=G_\a B$, $A=G_\a(A\cap B)$, and $|G_\a|$ is divisible  $|G:B|$; in particular, $|G_\a|$ is divisible  by ${1\over 2}q^3(q^3\pm1)$ or ${q^6-1\over 2(q\pm 1)}$.
Assume that $|G_\a|$ is divisible by ${1\over 2}q^3(q^3\pm1)$.
By \cite[Theorem 1]{Factor-excep}, combining  Lemma \ref{stab-lem-insol}, we have $\SL_3(q)\le  G_\a\le \SL_3(q).2$ and $p=3$. It follows that $r={q^3-1\over q-1}$, and so $f$ is a power of $3$ by Lemma \ref{r}. In particular,  $(4,q+1)=4$. Moreover, in this case,  $|T|$ is divisible by $|G:G_\a|={1\over 2a}q^6(q^4-1)(q^3+1)$, where $a\le 2$. Then
$T\le B= \Omega_6^-(q){.}2$.
By \cite[Tables 10.3 and 10.4]{transubgroup}, $\POmega_6^-(q)$ has no proper subgroup of order divisible by $q^6(q^4-1)(q^3+1)\over q-1$. It follows that
$T= \POmega_6^-(q)$. This forces that $B$ has trivial center, and then $(4,q+1)\le 2$, a contradiction. Assume that $|G_\a|$ is divisible by ${q^6-1\over 2(q\pm 1)}$. Then $|G_\a|$ has a divisor $\Phi_6^*(q)\Phi_3^*(q)$. By \cite[Table 10.5]{transubgroup}, we have $G_\a\lesssim  \PSL_2(13)$, forcing $G_\a\cong \PSL_2(13)$ and $r=13+1$, which contradicts that $r$ is  a prime.

Now $G_\a=Y\le B$ and $T=X\le A$. Then $G=AG_\a$, yielding
$|G_\a:(A\cap G_\a)|=|G:A|={1\over 2}q^3(q^4-1)$.
Assume   that $B\succcurlyeq\PSp_4(q)$ or
$\PSL_4(q)$. Let $\overline{B}=B'\Rad(B)/\Rad(B)$. Then $\overline{B}\cong \PSp_4(q)$ or
$\PSL_4(q)$. Note that $\Phi^*_4(q)\ne 1$.
Considering the preimage of $(G_\a\cap B')\Rad(B)/\Rad(B)$ in $\GL_4(q)$, by Lemma \ref{linear-stab-1},
we conclude that $q=3$ or $7$, and $S\cong \A_5$ or $\A_7$, respectively. By  Lemma \ref{stab-lem-insol},
$|G_\a|$ is indivisible by $3^3$, which is impossible as  $|G_\a|$ has a divisor $q^3$.

The left case is that $G_\a\le B\cong \Omega_6^-(q){.}2$ and $T\le A\cong \G_2(q)$.
In particular,
$|G_\a:(A\cap G_\a)|=|G:A|={1\over 2}q^3(q^4-1)$ and $|T:(B\cap T)|=|G:B|={1\over 2}q^3(q^3-1)$. It is easily shown that $S\preccurlyeq G_\a\cap B'$. Since  $|B:B'|=2$, either $G_\a\le B'$ or
$B=G_\a B'$, yielding $|G:(G_\a\cap B')|=1$ or $2$, respectively.
It follows that $|G_\a\cap B'|$ is divisible by $q^3\Phi_4^*(q)$.
Take a maximal subgroup $M$ of $B'$ with
$G_\a\cap B'\le M$. Then $|M|$  is divisible by $q^3\Phi_4^*(q)$.
Checking the maximal subgroups of $\Omega_6^-(q)$, refer to \cite[Tables 8.33 and 8.34]{Low}, we conclude that
$M\cong \Omega_5(q){.}2\cong \PSp_4(q){.}2$.
This gives a similar contradiction as in the above paragraph.

\vskip 5pt

{\bf Case 2}. Let $A$ and $B$ be as
in  \cite[Table 3]{max-factor}. Then $G=\Omega_7(3)$ and, up to isomorphism, $A$ and $B$ are listed as follow:
\[\tiny
\begin{array}{l|l|l|l|l|l|l|l}\hline
A& \G_2(3)&\G_2(3)&\Sy_9&\Sy_9&\Sp_6(2)&\Sp_6(2)&2^6{.}\A_7\\
\hline
|G:A|&2^3\cdot3^3\cdot5&2^3\cdot3^3\cdot5&2^2\cdot3^5\cdot13&2^2\cdot3^5\cdot13&
3^5\cdot 13&3^5\cdot13 &3^7\cdot 13\\ \hline\hline
B& \Sp_6(2)&\Sy_9&\PSL_4(3).2&[3^6]{:}\PSL_3(3)&\PSL_4(3).2&[3^6]{:}\PSL_3(3)
&[3^6]{:}\PSL_3(3)\\
\hline
|G:B|&3^5\cdot13&2^2\cdot3^5\cdot13&2\cdot3^3\cdot7&2^5\cdot5\cdot7&
2\cdot3^3\cdot 7&2^5\cdot5\cdot 7 &2^5\cdot5\cdot 7\\ \hline
\end{array}
\]
Note that $|X|$ and $|Y|$ are divisible by
$|G:B|$ and $|G:A|$, respectively. Inspecting simple subgroups contained in $A$ and $B$, since neither $X$ nor $Y$ is an alternating group, we conclude that either $T=A\cong\G_2(3)$ or $\PSU_4(2)\cong T\le B\cong\Sp_6(2)$. Suppose that $T\cong \G_2(3)$. Then $G_\a\le B\cong\Sp_6(2)$ or $\Sy_9$ and, since $G=TG_\a$, we know that $|G_\a|$ is divisible by $2^3\cdot 3^3\cdot 5$.
By Lemma \ref{stab-lem-insol}, $G_\a\ne B$.
Let $M$ be a maximal subgroup of $B$ with $G_\a\le M$.
Then, by the Atlas \cite{Atlas}, $M\cong \PSU_4(2).2$, $\Sy_3\times \Sy_6$ or $(3\times \A_6){:}2$. However, such an $M$ does not contain an insolvable subgroup satisfying Lemma \ref{stab-lem-insol} and having order divisible by $2^3\cdot 3^3\cdot 5$.
Thus $\PSU_4(2)\cong T\le B\cong\Sp_6(2)$, and $G_\a\le A\cong \G_2(3)$. Then $|G:T|=2^3\cdot 3^5\cdot7\cdot 13$. Noting that
$|G_\a|$ is divisible by $|G:T|$, it follows that $G_\a=A\cong \G_2(3)$, which is impossible by Lemma \ref{stab-lem-insol}.
\qed

\begin{theorem}\label{O-2m+1}
Assume that {\rm Hypothesis \ref{hypo-graphs-classical}} holds  and $G=\Omega_{2m+1}(q)$, where $m\ge 3$  and $q=p^f$ for some odd prime $p$. Then $mf$ is a power of some odd prime $d$, $T\cong \POmega^-_{2m}(q)$,  $\soc(G_\a^{\Ga(\a)})\cong\PSL_{d}(q^{m\over d})$, $p\equiv 1\,(\mod 4)$ and $r={q^m-1\over q^{m\over d}-1}$.
\end{theorem}
\proof By Lemma \ref{O-2m+1-1}, we let $A$ and $B$ be as
in  \cite[Table 1]{max-factor}. Then
$A\cong\Omega^-_{2m}(q).2$, $B\cong [q^{m(m+1)\over 2}]{:}{1\over 2}\GL_m(q)$ and $|G:B|=\prod_{i=1}^m(q^i+1)$.
In particular, $|G:B|$ and hence $|X|$ is divisible by $\Phi_{2m}^*(q)\Phi_{2m-2}^*(q)\Phi_{2m-4}^*(q)$.

{\bf Case 1}. Assume that $\Phi_{2m-4}^*(q)=1$. Then $m=3$ and $q=p=2^s-1$ for some  prime $s\ge 2$. In this case, we have
$A\cong\Omega_6^-(q){.}2\cong 2{.}\PSU_4(q){.}2$.

Suppose that $T=X\le A$. Then $T\lesssim \Omega_6^-(q)$.
Noting that $(4,q+1)=4$, we have $T\not \cong  \PSU_4(q)$.
Recall that  $|X|$ is divisible
by $(q+1)(q^2+1)(q^3+1)$,  and hence $|T|$ is divisible by $\Phi_{2m}^*(q)\Phi_{2m-2}^*(q)$. If  $q\ge 5$ and $|T|$ is divisible by $p$ then $T\cong \PSU_4(q)$ by \cite[Table 10.3]{transubgroup},  which is impossible. If $q=3$ then $|T|$ is divisible by $35$, and so
$T\cong \PSU_4(3)$ by \cite[Table 10.4]{transubgroup},   a contradiction. Thus
we have $q\ge 5$ and $|T|$ is indivisible by $p$.
Note that $|T|$ is divisible by $\Phi_{6f}^*(p)\Phi_{4}^*(q)$.
Considering the preimage of $T$ in
 $\GL_6(q)$, by \cite[Theorem 3.1]{Bam-Pen}, we conclude that $T\le \GU_3(q)$. Then $T\le  \PSU_3(q)$, and thus
  $|T|$ is indivisible by $(q+1)(q^2+1)(q^3+1)$,
a contradiction.

Now $T\le B\cong [q^{6}]{:}{1\over 2}\GL_3(q)$ and $G_\a\le A\cong \Omega_6^-(q){.}2$. Then $G=G_\a B$, and $|G_\a:(B\cap G_\a)|=|G:B|=\prod_{i=1}^3(q^i+1)$. Thus $|G_\a|$ is divisible
$\Phi^*_6(q)$. By Lemma \ref{linear-stab-1}, we conclude that
$q=3$ and $S\cong \A_7$ or $\PSL_3(2)$.  By Lemma \ref{stab-lem-insol}, calculation shows that $|G_\a|_3\le 3^3$.
Then $3^9=|G|_3\le |T|_3|G_\a|_3\le 3^6$,
a contradiction.

{\bf Case 2}. Assume that $\Phi_{2m-4}^*(q)\ne 1$.
By \cite[Theorem 4 and Table 10.1]{transubgroup}, since $q$ is odd, $\POmega^-_{2m}(q)$ has no proper subgroup
with order divisible by $\Phi_{2m}^*(q)\Phi_{2m-2}^*(q)\Phi_{2m-4}^*(q)$. It follows that $X\succcurlyeq\POmega^-_{2m}(q)$.
Therefore, $T=X\cong\POmega^-_{2m}(q)$, and
$G_\a=Y\le B$. This implies that $(4,q^m+1)=2$ and $|G:T|=q^m(q^m-1)$. In particular,
$q^m(q^m-1)$ is a divisor of $|G_\a|$, and so $|G_\a|$ is divisible by $\Phi_{mf}^*(p)$.
Clearly, $\Phi_{mf}^*(p)\ne 1$.

Let $\overline{G_\a}=G_\a\O_p(B)/\O_p(B)$. Then $\overline{G_\a}\le {1\over 2}\GL_m(q)$, and $|\overline{G_\a}|$ is divisible by $q^m-1$. Note that $S\preccurlyeq\overline{G_\a}$.
By Lemma \ref{linear-stab-1},
 either  $S\cong\PSL_{d}(q^{m\over d})$ for some divisor $d$   of $m$, or
$(S,m,q)$ is one of $(\A_5,4,3)$, $(\A_7,6,3)$, $(\A_7,6,5)$,
$(\A_7,4,7)$, $(\A_7,3,25)$, $(\PSL_3(2),6,3)$, $(\PSL_3(2),6,5)$, $(\PSL_3(2),3,9)$ and $(\PSL_3(2),3,25)$.
Since $|G_\a|$ is divisible by $q^m$,  the latter case is easily excluded by calculating $|G_\a|_p$. Assume that the former case occurs. If $(d, q^{m\over d},r)=(2,11,11)$ then we have $|G_\a|_{11}=11$ by Lemma \ref{stab-lem-insol}, which contradicts that $|G_\a|$ is divisible by $q^m$. Thus we have $r={q^m-1\over q^{m\over d}-1}$,  and so
 $d$ is an odd prime
and $mf$ is a power of $d$, see  Lemma \ref{r}. Recalling $(4,q^m+1)=2$, we have $p\equiv 1\,(\mod 4)$.
\qed


\subsection{Orthogonal case II}
Let  $G=\POmega^-_{2m}(q)$, where $m\ge 4$  and $q=p^f$ for some  prime $p$. Then, by \cite[Theorem A]{max-factor}, $m$ is odd and
one of the following holds:
\begin{itemize}
\item[(1)] $m=5$, $q=2$, $A\cong \A_{12}$ and $B\cong 2^8{:}\Omega^-_{8}(2)$;

\item[(2)] $A\cong [q^{2m-2}]{:}[{q-1\over a}]{.}\POmega^-_{2m-2}(q){.}a$ and $B\cong {q+1\over (4,q+1)}{.}\PGU_m(q)$, where $a=2$ if $q\equiv 1\,(\mod 4)$, or $a=1$ otherwise;

\item[(3)]   $A\cong \Omega_{2m-1}(q).{(4,q-1)\over 2}$ (odd $q$) or $\Sp_{2m-2}(q)$ (even $q$),  and $B\cong {q+1\over (4,q+1)}{.}\PGU_m(q)$.
\end{itemize}

For case (1), we have $|G:B|=3^2\cdot5\cdot11$ and $|G:A|=2^{11}\cdot3\cdot17$. Check the insolvable subgroups of $A$ and $B$ with order divisible by $|G:B|$ and $|G:A|$, respectively. Recalling that $X$ is not an alternating group, it follows that $X\cong \M_{11}$ or $\M_{12}$, and $Y\gtrsim \Omega^-_{8}(2)$. Then the only possibility is that $T=Y\cong \Omega^-_{8}(2)$ and $G_\a\cong \M_{11}$. Thus $|T|_2|G_\a|_2=2^{16}<|G|_2$, a contradiction.

Suppose that (2) or (3) occurs. It is easy to check that
$|G:A|$ is divisible by $q^{m-1}{q^m+1\over (4,q^m+1)}$, and then
$|Y|$ is divisible by $q^{m-1}\Phi^*_{m(2f)}(p)$. Since $m\ge 4$ and $m$ is odd, $\Phi^*_{m(2f)}(p)\ne 1$. By \cite[Theorem 3.1]{Bam-Pen},  one of the following holds:
$Y\succcurlyeq\PSU_{m_0}(q^{m\over m_0})$ for some divisor $m_0\ge 3$ of $m$; or  $m=5$, $q^2=4$ and $Y\succcurlyeq\PSL_2(11)$; or  $m=9$, $q^2=4$ and $Y\succcurlyeq\J_3$ or $\PSL_2(19)$. For the last two cases, we have $|Y|_p<q^{m-1}$, a contradiction. Thus
we have $T=Y\cong \PSU_{m_0}(q^{m\over m_0})$, and then $G_\a=X\le A$.


\begin{theorem}\label{O-2m-}
Assume that {\rm Hypothesis \ref{hypo-graphs-classical}} holds, $G=\POmega_{2m}^-(q)$ for $m\ge 4$.
Then $m=5$, $q>2$,    $T\cong \PSU_5(q)$, $\SL_2(q^4)\unlhd G_\a^{\Ga(\a)}$,  $|\O_p(G_\a)|=q^8$ and $f$ is a power of $2$.
\end{theorem}
\proof
By the argument ahead of this lemma, $m$ is odd,
 $\PSU_{m_0}(q^{m\over m_0})\cong T\le B$, $G_\a\le A$, where  $m_0$ is a divisor of $m$.
Note that $|G:B|=q^{m(m-1)\over 2}\prod_{i=1}^{m-1}(q^i-(-1)^{i-1})$. Then $|G_\a|$ is divisible by $q^{m(m-1)\over 2}\Phi_{2(m-1)}^*(q)\Phi_{2(m-3)}^*(q)$.

 Assume that $A\cong[q^{2m-2}]{:}[{q-1\over a}]{.}\POmega^-_{2m-2}(q){.}a$.
Let $\widehat{A}$ and $\widehat{G_\a}$ be the preimages of $A'\Rad(A)/\Rad(A)$ and $(A'\cap G_\a)\Rad(A)/\Rad(A)$
in $\GO^-_{2(m-1)}(q)$, respectively. Then $S\preccurlyeq \widehat{G_\a}$, and $|\widehat{G_\a}|$ is divisible by $\Phi_{2(m-1)f}^*(p)\Phi_{2(m-3)}^*(q)$. By Lemma \ref{linear-stab-1}, either $S=\PSL_d(q^{2(m-1)\over d})$ or
$(S,2(m-1),q)$ is one of $(\A_{13},12,2)$ and $(\PSL_3(3),12,2)$.
For the latter case, by Lemma \ref{stab-lem-insol}, we conclude that $|G_\a|_2\le 2^{11}$, which contradicts that $|G_\a|$ is divisible by $q^{m(m-1)\over 2}$. Suppose that $S=\PSL_d(q^{2(m-1)\over d})$. Then $r={(q^{2(m-1)}-1)\over (q^{2(m-1)\over d}-1)}$. By Lemma \ref{r}, $d=2$, $q$ is even and $m-1$ is a power of $2$. By the choice of $\widehat{G_\a}$, we conclude that
$|G_\a|_p\le |\Rad(A)|_p|\widehat{G_\a}|_p$. Then, by Lemma \ref{linear-stab-1}, we have $|G_\a|_p\le |\Rad(A)|_p2(m-1)|S|_p=2(m-1)q^{3(m-1)}$.
Thus $q^{m(m-1)\over 2}\le |G_\a|_p\le 2(m-1)q^{3(m-1)}$.
This implies that $m\le 7$, yielding $m=5$ as $m\ge 5$ and $m-1$ is a power of $2$. Then $T\cong\PSU_{5}(q)$.
Noting that $8q^{12}\ge |G_\a|_p\ge q^{10}$,  by Lemma \ref{stab-lem-insol}, we have $|\O_p(G_\a)|=q^8$.
To finish the proof, we next show that $A\not\cong \Omega_{2m-1}(q).{(4,q-1)\over 2}$ or $\Sp_{2m-2}(q)$.

Suppose  that $A\cong \Omega_{2m-1}(q).{(4,q-1)\over 2}$ for odd $q$.
Recall that $|G_\a|$ is divisible by $q^{m(m-1)\over 2}\Phi_{2(m-1)}^*(q)\Phi_{2(m-3)}^*(q)$, and $2(m-1)\ge 8$.
By Lemma \ref{linear-stab-2}, the only possibility is that
$S\cong \PSL_d(q^{2(m-1)\over d})$, $d$ is a prime and $2(m-1)f$ is a power of $d$. This implies that $d=2$, and so $r=q^{m-1}+1$. Then $q$ is even by Lemma \ref{r}, a contradiction.

Suppose that $A\cong \Sp_{2m-2}(q)$ for even $q$. Clearly, $G_\a\ne A$.  By Lemma \ref{linear-stab-1}, we conclude that
either $S\cong \PSL_2(q^{m-1})$ or
$(S,2(m-1),q)$ is one of $(\A_{13},12,2)$ and $(\PSL_3(3),12,2)$.
For the latter case,   $|G_\a|_2\le 2^{11}$ by Lemma \ref{stab-lem-insol}, and so   $2^{21}=q^{m(m-1)\over 2}\le |G_\a|_2\le 2^{11}$, a contradiction. Then $S\cong \PSL_2(q^{m-1})$ and, by Lemma \ref{linear-stab-1},  $|G_\a|_p\le 2(m-1)|S|_p$.
Thus $q^{m(m-1)\over 2}\le 2(m-1)q^{m-1}\le q^{2(m-1)}$, yielding $m\le 4$,  contradiction.
\qed


\subsection{Orthogonal case III}
Let  $G=\POmega^+_{8}(q)$,   $q=p^f$ for some  prime $p$.
 Then $A$ and $B$ are  as in \cite[Table 4]{max-factor}.
 Recall that $\{X,Y\}=\{T,G_\a\}$ with $X\le A$ and $Y\le B$.
\begin{lemma}\label{O+8-small-q}
Either $q>3$, or one of the following holds.
\begin{itemize}
\item[(1)] $q=2$, $A\cong \Sp_6(2)\cong B$, $T\cong \Sp_6(2)$ and $G_\a\cong \Sy_5$.
\item[(2)]   $q=3$, $T=A\cong \Omega^+_8(2)$, $G_\a\cong [3^6]{:}\SL_3(3)$ and $r=13$.
\end{itemize}
\end{lemma}
\proof
Let $q=2$ or $3$. Then, up to isomorphism,  $A$ and $B$ are listed as follows:
\[\tiny
\begin{array}{l|l|l||l|l|l}
X&A& |G:A|&|G:B|&B&Y\\ \hline
&\A_9& 2^6\cdot 3\cdot 5&  2^3\cdot 3\cdot 5&\Sp_6(2)&\\

& & &3^3\cdot 5& 2^6{:}\A_8&\\

&&& 2^5\cdot5\cdot7&(3{\times}\PSU_4(2)){:}2\\ \hline
& \Sp_6(2)& 2^3\cdot 3\cdot 5& 2^6\cdot 3^3\cdot 7&(\A_5\times\A_5){:}2^4&\\
&&&  2^3\cdot 3\cdot 5&\Sp_6(2)&\\
&&&3^3\cdot 5& 2^6{:}\A_8&\\
&&&2^5\cdot5\cdot7&(3{\times}\PSU_4(2)){:}2&\\ \hline
&(3{\times}\PSU_4(2)){:}2&2^5\cdot5\cdot7&3^3\cdot 5& 2^6{:}\A_8&\\ \hline \hline
&\Omega_8^+(2)&3^7\cdot13&2^3\cdot3^3\cdot5&\Omega_7(3)& \\
&&&2^5\cdot5\cdot7&3^6{:}\PSL_4(3)&\\ \hline
&\Omega_7(3)&2^3\cdot3^3\cdot5&&\Omega_7(3)&\\
&&&2^5\cdot5\cdot7&3^6{:}\PSL_4(3)&\\
&&&2^2\cdot3^6\cdot5\cdot13& 2^.\PSU_4(3){.}2^2&\\
&&&2^3\cdot3^7\cdot5\cdot7\cdot13& (\A_4\times\PSU_4(2)){.}2&\\ \hline
&2^.\PSU_4(3){.}2^2&2^2\cdot3^6\cdot5\cdot13&2^5\cdot5\cdot7&3^6{:}\PSL_4(3)&

\end{array}
\]


{\bf Case 1}.
Let $q=2$. With the help of GAP, we search the subgroups of $A$ with order divisible by $|G:B|$.
Up to isomorphism, either $A\cong \A_9$ and $G_\a=X\cong \Sy_5,\, 2^2{\times}\A_5, \,  \A_4{\times}\A_5, \, (2^2{\times}\A_5){.}2,\,(\A_4{\times}\A_5){.}2$ or $
\Sy_7$, or  $X$ and $Y$ are listed as follows:
\[\tiny
\begin{array}{l|l|l|l}
A&  X&B&Y\\ \hline

\Sp_6(2)&
 \Sp_6(2)& (\A_5\times\A_5){:}2^4& \Sy_5, 2^2{\times}\A_5,(2^2{\times}\A_5){.}2, \A_4{\times}\A_5,(\A_4{\times}\A_5){.}\\
& \Sy_5,
 2^4{:}\A_5, 2^4{:}\Sy_5,
  \Sy_7& \Sp_6(2)&\Omega_5(3),  \Sp_6(2)\\
 &\Omega_5(3),
 \Sp_6(2)& \Sp_6(2)&\Sy_5,
 2^4{:}\A_5, 2^4{:}\Sy_5,   \Sy_7\\
 &\Omega_5(3),
\Sp_6(2) &2^6{:}\A_8&\Sy_5, 2^2{\times}\A_5,(2^2{\times}\A_5){.}2, \A_4{\times}\A_5,(\A_4{\times}\A_5){.}2\\
&&&2^4{:}\A_5,2^4{:}\Sy_5, \A_4\times 2^4{:}\A_5,(\A_4\times 2^4{:}\A_5){.}2\\
& \Sp_6(2)&(3{\times}\PSU_4(2)){:}2&\Sy_5,2^4{:}\A_5,2^4{:}\Sy_5\\ \hline
(3{\times}\PSU_4(2)){:}2&\Omega_5(3)&2^6{:}\A_8&\mbox{No}
\end{array}
\]
Then the vertex-stabilizer $G_\a$ is easily determined.
By calculation, we have a desired graph only if $A\cong \Sp_6(2)\cong B$, $T\cong \Sp_6(2)$ and $G_\a\cong \Sy_5$.

{\bf Case 2}. Let $q=3$ and $B\not\cong\Omega_7(3)$. Set
$\overline{Y}=Y\Rad(B)\Rad(B)$ and $\overline{B}=B/\Rad(B)$.

Suppose that $B\cong 3^6{:}\PSL_4(3)$. Then $A\cong \Omega_8^+(2)$,
$\Omega_7(3)$ or $2^.\PSU_4(3){.}2^2$, and
 $|\overline{Y}|$ has a divisor $5\cdot13$. By \cite[Table 10.3]{transubgroup}, we have $\overline{Y}=\overline{B}\cong\PSL_4(3)$. Thus $Y\succcurlyeq\PSL_4(3)$. Noting that $3^4-1\over 3-1$ is not a prime, we have $S\not \cong\PSL_4(3)$, and then $B>T\cong\PSL_4(3)$. Thus $|G:T|=|G:B||B:T|=2^5\cdot5\cdot7\cdot3^6$, and so $G_\a$ is a subgroup of $A$ with order divisible by $2^5\cdot3^6\cdot5\cdot7$. By \cite[Table 10.3]{transubgroup}, we conclude that $A\not\cong \Omega_8^+(2)$. Then $A\cong \Omega_7(3)$  or $2^.\PSU_4(3){.}2^2$.
If $A\cong \Omega_7(3)$  then, checking the maximal subgroup of $A$ with order divisible by $|G:T|$, we conclude that
$G_\a$ is contained in a maximal subgroup $2^.\PSU_4(3){.}2^2$ of $A$. Thus we may let $G_\a\le A_0$, where $2^.\PSU_4(3){.}2^2\le \A_0\le A$. Considering $G_\a\Rad(A_0)/\Rad(A_0)$ and $A_0/\Rad(A_0)$, by \cite[Table 10.4]{transubgroup}, we know that
$G_\a\succcurlyeq\PSL_3(4)$ or $\A_7$.
Then $S\cong \PSL_3(4)$ or $\A_7$. Noting $4^3-1\over 4-1$ is not a prime,   $S\cong \A_7$. By Lemma \ref{stab-lem-insol},  $|G_\a|_3=3^3$, a contradiction.

Suppose that $B\cong 2^.\PSU_4(3){.}2^2$. Then $A\cong \Omega_7(3)$.
It is easily shown that $\overline{Y}$ and $\overline{Y}\cap \overline{B}^\infty$ have same insolvable composition factors, and $|\overline{Y}\cap \overline{B}^\infty|$ is divisible by $3^3\cdot5$. With the help of GAP,  searching the subgroups of $\PSU_4(3)$ with order divisible by $3^3\cdot5$, we get $\overline{Y}\cap \overline{B}^\infty\cong 3^4{:}\A_5$, $\Omega_5(3)$, $3^4{:}\A_6$ or $\PSU_4(3)$. Noting that
$\overline{Y}=Y\Rad(B)\Rad(B)\cong Y/(Y\cap \Rad(B))$, it follows that $Y$ has a normal subgroup $R{.}(3^4{:}\A_5)$, $R{.}\Omega_5(3)$, $R{.}(3^4{:}\A_6)$ or $R{.}\PSU_4(3)$, where $R=Y\cap \Rad(B)$. Considering the limitations on $Y$, we conclude that $T=Y\cong \Omega_5(3)$, and $G_\a=X\le A$. Noting that $|G:T|=2^6\cdot 3^8\cdot \cdot5\cdot7\cdot13$, this leads to  a similar contradiction as in the above paragraph.

Suppose that $B\cong (\A_4\times\PSU_4(2)){.}2$. Then $A\cong \Omega_7(3)$, and $|\overline{Y}\cap \overline{B}^\infty|$ is divisible by $3^2\cdot5$.
Searching the subgroups of $\PSU_4(2)$ with order divisible by $3^2\cdot5$, we conclude that $\overline{Y}\cap \overline{B}^\infty\cong\A_6,\,\Sy_6$ or $\Omega_5(3)$. It follows that $T=Y\cong \Omega_5(3)$, and $G_\a=X\le A$, which leads to a contradiction.

{\bf Case 3}. Let  $B\cong \Omega_7(3)$ and $A\cong \Omega^+_8(2)$ or $\Omega_7(3)$.
Suppose that $G_\a\lesssim \Omega^+_8(2)$.
Then, by Lemma \ref{linear-stab-4}, $S\cong \A_5$ or $\A_7$, which leads to a contradiction by calculating $|G_\a|_3$.
Thus $T\le A\cong \Omega^+_8(2)$, and $G_\a\le B$. In particular,  $|G:A|=3^7\cdot13$.
Then $G_\a$ is a subgroup of $B$ of order divisible by $3^7\cdot13$. Checking the orders of maximal subgroups of $\Omega_7(3)$, we conclude that $\PSL_3(3)\lesssim G_\a\lesssim [3^6]{:}\PSL_3(3)\lesssim B$. Then $\soc(G_\a^{\Ga(\a)}\cong \PSL_3(3)$, $r=13$ and $|\O_3(G_\a)|=3^6$ by Lemma \ref{stab-lem-insol}. Thus $G_\a\cong [3^6]{:}\SL_3(3)$, and $|G:G_\a|=|G:B||B:G_\a|=2^8\cdot3^3\cdot5^2\cdot7$. Since
$|T|$ is divisible by $|G:G_\a|$, by \cite[Table 10.4]{transubgroup}, we have $T=A\cong \Omega^+_8(2)$, and part (2) of this lemma follows.
\qed

\vskip 10pt

For $q>3$, we list  $A$ and $B$   as follows:
\[\tiny
\begin{array}{l|l|l|l}
A&|G:A|&B&|G:B|\\ \hline
\Omega_7(q)& {1\over e}q^3(q^4-1)&\Omega_7(q)&  {1\over e}q^3(q^4-1)\\
\Omega_7(q)& {1\over e}q^3(q^4-1)& \sP_1, \sP_3\mbox{ or }\sP_4& (q^3+1)(q^2+1)(q+1)\\
\Omega_7(q)& {1\over e}q^3(q^4-1)&{}^{\bf\hat{}}({q+1\over e}{\times} \Omega^-_6(q)){.}2^e&{1\over2}q^6(q^3-1)(q^2+1)(q-1)\\
\Omega_7(q)& {1\over e}q^3(q^4-1)&{}^{\bf\hat{}}({q-1\over e}{\times} \Omega^+_6(q)){.}2^e&{1\over2}q^6(q^3+1)(q^2+1)(q+1)\\
\Omega_7(q), q\mbox{ ood}& {1\over 2}q^3(q^4-1)&(\PSp_2(q)\otimes\PSp_4(q)){.}2&{1\over2}q^7(q^6-1)(q^2+1)\\
\Omega_7(q)& {1\over e}q^3(q^4-1)&\Omega^-_8(\sqrt{q})&{1\over e}q^6(q^4-1)(q^3+1)(q+1)\\ \hline
{}^{\bf\hat{}}({q+1\over e}{\times} \Omega^-_6(q)){.}2^e&{1\over2}q^6(q^3-1)(q^2+1)(q-1)
&\sP_1, \sP_3\mbox{ or }\sP_4& (q^3+1)(q^2+1)(q+1)
\end{array}
\]
where $e=(2,q-1)$. Note that $\Omega_m(q)=\Sp_{m-1}(q)$ for odd $m$ and even $q$.

\begin{lemma}\label{O+8-q>3-1}
Assume that $q>3$. Then $A\cong \Omega_7(q)$.
\end{lemma}
\proof
 Suppose   that $A$ is of type ${}^{\bf\hat{}}({q+1\over e}{\times} \Omega^-_6(q)){.}2^e$. Let $\overline{A}=A^\infty\Rad(A)/\Rad(A)$ and $\overline{X}=(X\cap A^\infty)\Rad(A)/\Rad(A)$. Then $\overline{A}\cong \PSU_4(q)$ and $|\overline{X}|$ is divisible by $\Phi^*_6(q)\Phi^*_4(q)$.
 By \cite[Table 10.3]{transubgroup}, either $\overline{X}=\overline{A}\cong \PSU_4(q)$, or $p\ge 5$ and $|\overline{X}|$ is indivisible by $p$.
  Assume that $\overline{X}\ne\overline{A}$.
 Checking the orders of maximal subgroups of $\PSU_4(q)$ (refer to \cite[Tables 8.10 and 8.11]{Low}), we conclude that
  $\overline{X}\cong \A_7$.  Considering the order of $\A_7$, it follows from Lemma \ref{phi} that $\Phi^*_6(q)=7$  and $\Phi^*_4(q)=5$. Then $q=3$ by \cite[Theorem 3.9]{Hering}, a contradiction. Thus  $\overline{X}=\overline{A}\cong \PSU_4(q)$, and then $T=X \cong \PSU_4(q)$, yielding $(2,q-1)=1$ and $|G:T|=q^6(q^4-1)(q^3-1)$.
  Moreover, $G_\a=Y\le B$.

  Let $\overline{B}=B^\infty\Rad(B)/\Rad(B)$ and $\overline{Y}=(Y\cap B^\infty)\Rad(B)/\Rad(B)$. Then $\overline{B}\cong \PSL_4(q)$, and $|\overline{Y}|$ has a divisor $\Phi^*_3(q)\Phi^*_4(q)$.
   By \cite[Table 10.3]{transubgroup}, either $\overline{Y}=\overline{B}\cong \PSL_4(q)$, or $p\ge 5$ and $|\overline{Y}|$ is indivisible by $p$. For the former, we conclude that $G_\a=Y\succcurlyeq\PSL_4(q)$, a contradiction. Then $|\overline{Y}|$ is indivisible by $p$.
  Checking the orders of maximal subgroups of $\PSL_4(q)$ (refer to \cite[Tables 8.10 and 8.11]{Low}), we get
  $\overline{Y}\lesssim 2^.\PSL_2(7)$ or $2^.\A_7$, and $q=p$.
  Recalling $(2,q-1)=1$, we have $q=p=2$, a contradiction.
\qed

\begin{lemma}\label{O+8-q>3-2}
Assume that $q>3$. Then $T\lesssim \Omega_7(q)$.
\end{lemma}
\proof
Suppose that $T\not\lesssim \Omega_7(q)$. Then, by Lemma \ref{O+8-q>3-1}, $T\le B\not\cong \Omega_7(q)$, yielding $G_\a\le A\cong \Omega_7(q)$. Then $|G_\a|$ is divisible by $|G:B|$.
It follows that either $T\lesssim\PSU_4(q)$ or $|G_\a|$ is divisible by $\Phi^*_6(q)\Phi^*_4(q)$.
Assume that the latter case occurs. Noting $6$ is not a prime power, it follows from
 Lemmas \ref{linear-stab-1} and \ref{linear-stab-2}  that $q=5$ and $S\cong \A_7$ or $\PSL_3(2)$. Then $\Phi^*_4(q)=13\not\in \pi(S)$, which contradicts that $\pi(G_\a)=\pi(S)$, see Lemma \ref{stab-lem-insol}. Thus we have $B\cong {}^{\bf\hat{}}({q+1\over e}{\times} \Omega^-_6(q)){.}2^e$, and
 $|G:B|={1\over2}q^6(q^3-1)(q^2+1)(q-1)$, where $e=(2,q-1)$.
 Note that $|G_\a|$ is divisible by $|G:B|$.
Take a  maximal subgroup $A_0$   of $A$ with $G_\a\le A_0$.
 Then both $|A_0|$ and $|G_\a|$  are divisible by ${1\over2}q^6\Phi^*_3(q)\Phi^*_4(q)$.
 Checking the orders of maximal subgroups of $\Omega_7(q)$ and $\Sp_6(q)$, refer to
\cite[Tables 8.28, 8.29, 8.39 and 8.40]{Low},  we conclude that $A_0\cong \Omega^+_6(q){.}2\cong \PSL_4(q){.}2$.
It follows from \cite[Table 10.3]{transubgroup} that $\PSL_4(q)\preccurlyeq G_\a$, which is impossible by Lemma \ref{stab-lem-insol}.
\qed

\begin{theorem}\label{O+8}
Assume that {\rm Hypothesis \ref{hypo-graphs-classical}} holds  and  $G=\POmega^+_{8}(q)$, where  $q=p^f$ for some  prime $p$. Then one of the following holds.
\begin{itemize}
\item[(1)] $p=2$, $f$ is a power of $2$, $T\cong \Sp_6(q)$, $\SL_2(q^2)\unlhd G_\a^{\Ga(\a)}$ and $r=q^2+1$; in this case,  $G_\a\cong \Sy_5$ if $f=1$, and $\O_2(G_\a)\ne 1$ if $f>2$.

\item[(2)]   $q=3$, $T\cong \Omega^+_8(2)$, $G_\a\cong [3^6]{:}\SL_3(3)$ and $r=13$.

\end{itemize}
\end{theorem}
\proof
If $q\le 3$ then our result holds by   Lemma \ref{O+8-small-q}.
Assume that $q>3$. By Lemmas \ref{O+8-q>3-1} and \ref{O+8-q>3-2}, $T\le A\cong \Omega_7(q)$ and $G_\a\le B$.
Then $S$ is given as in Lemma \ref{linear-stab-4}.
Note that $|G_\a|$ is divisible by $|G:A|={1\over (2,q-1)}q^3(q^4-1)$.
If $(S,q)=(\A_7,7)$ then $|G_\a|_7=7<|G:A|_7$, a contradiction.
Thus $S\cong \PSL_2(q^2)$, $r=q^2+1$, $p=2$ and $f$ is a power of $2$. If $\O_p(G_\a)=1$ then, by Lemma \ref{stab-lem-insol}, $|G_\a|_2$ is a divisor of
$2fq^2$, and so $q^3 \le 2fq^2$, yielding $f\le 2$. Since $q$ is even, $A\cong \Sp_6(q)$.

By Lemma \ref{stab-lem-insol}, $|G_\a|$ is a divisor of $4f^2q^2(q^4-1)(q^2-1)|\O_2(G_\a)|$. It follows that $|G:G_\a|$ is divisible by $(q^4-1)(q^6-1)$. Since $q>3$, by \cite[Table 10.1]{transubgroup}, we have $T=A\cong \Sp_6(q)$, and the result follows.
\qed

\subsection{Orthogonal case IV}
Let  $G=\POmega^+_{2m}(q)$,   $m\ge 5$  and $q=p^f$ for some  prime $p$.

\begin{lemma}\label{O+2m-1}
 $A$ and $B$ are not described as in \cite[Table 3]{max-factor}.
\end{lemma}
\proof
Suppose that $G\cong \Omega^+_{24}(2)$, $A\cong \Co_1$ and $B\cong \Sp_{22}(2)$.
Then $|G:B|=2^{11}(2^{12}-1)$ and $|G:A|$ is divisible by $\Phi^*_{22}(2)\Phi^*_{20}(2)\Phi^*_{18}(2)=683\cdot 41\cdot 19$.
Since $|Y|$ is divisible by $|G:A|$, by \cite[Table 10.1]{transubgroup}, either $Y=B$ or $Y$ has a normal subgroup $\Omega^-_{22}(2)$, yielding $T=Y\cong \Sp_{22}(2)$ or $\Omega^-_{22}(2)$. Then $G_\a\le A$.
Checking  the insolvable  subgroups of $\Co_1$ with order divisible by $2^{11}(2^{12}-1)$, we conclude that $G_\a\succcurlyeq\G_{2}(4)$ or $\Suz$, which contradicts Lemma \ref{stab-lem-insol}.
\qed

\begin{lemma}\label{O+2m-2}
 $A$ and $B$ are not described as in \cite[Table 2]{max-factor}.
\end{lemma}
\proof

Suppose that $G=\POmega^+_{16}(q)$,
$A\cong \Omega_9(q){.}a$ and
$B=\sN_1$, where $a\le 2$.
Then $|G:B|$ and $|G:A|$ are divisible by $q^7{q^8-1\over (4,q^8-1)}$
and $q^{40}{\prod_{i\in \{8,10,12,14\}}(q^i-1)\over a(4,q^8-1)}$.

Assume that $q$ is even. Then $B\cong \Sp_{14}(q)$, and $|Y|$ is
divisible by $\Phi^*_{14}(q)\Phi^*_{12}(q)\Phi^*_{7}(q)$. By \cite[Table 10.1]{transubgroup}, we have $Y=B\cong \Sp_{14}(q)$, and so $G_\a=X\le A$.
Assume that $q$ is  odd. Then $B\cong \Omega_{15}(q){.}b$.
By Lemmas \ref{sub-factor} and \ref{m-factor}, either $Y\succcurlyeq\Omega_{15}(q)$, or
there exists a core-free maximal factorization  $B_0=HK$  such that $K$ contains a normal subgroup $N$ of $Y$ of index at most $2$.
The former cases yields that $T=Y=\Omega_{15}(q)$.
Suppose that the latter case occurs. It is easily shown that $|N|$ is divisible by $\Phi^*_{14}(q)\Phi^*_{12}(q)\Phi^*_{10}(q)$.
Then, by \cite[Theorem A]{max-factor},
we conclude that $K$ is the stabilizer of some $1$-dimensional nonsingular subspace, and so $N\le K\cong \Omega_{14}^-(q){.}2$. By \cite[Table 10.1]{transubgroup}, we have $N\succcurlyeq \POmega_{14}^-(q)$. Then $Y\succcurlyeq \POmega_{14}^-(q)$, yielding $T=Y\cong \POmega_{14}^-(q)$.

All in all,   $T\le B=\sN_1$. By Lemma \ref{linear-stab-6},
$S\cong \PSL_2(q^4)$, $p=2$ and $f$ is power of $2$.
Since $G_\a\le A\cong \Sp_8(q){.}a$, we have $S\preccurlyeq G_\a\cap A'$ and $|G_\a:(G_\a\cap A')|\le a$. By Lemma \ref{linear-stab-1}, $|G_\a\cap A'|_2\le 4|S|_2=4q^4$, and so $|G_\a|_2\le 2^3q^4$.  Recalling that $|G_\a|$ is divisible by $|G:B|$, we have $q^7\le |G_\a|_2\le 2^3q^4$, yielding $q=2$ and $|G_\a|_2=2^7$. Then  $\O_2(G_\a)=1$, otherwise, $|G_\a|_2\ge q^8$ by Lemma \ref{stab-lem-insol}, a contradiction. Thus $G_\a\lesssim\PGammaL_2(2^4)$, and then $|G_\a|_2\le 2^6$, a contradiction.
\qed

\vskip 10pt

We next deal with the case that $A$ and $B$ are described as in \cite[Table 1]{max-factor}. In particular, either $m$ is even, or $|G:B|=\prod_{i=1}^{m-1}(q^i+1)$.

\begin{lemma}\label{O+12(2)}
Assume that $A$ and $B$ are described as in \cite[Table 1]{max-factor}, and $(m,q)=(6,2)$. Then one of the following holds.
\begin{itemize}
\item[(1)] $T=\Omega^-_{10}(2)$ or $\Sp_{10}(2)$, and $\SL_5(2)\le G_\a\le [2^{20}]{:}\SL_5(2)$;
\item[(2)] $T=\Sp_{10}(2)$, $G_\a^{\Ga(\a)}\cong \SL_3(2)$ and $\O_2(G_\a)\ne 1$.
\end{itemize}
\end{lemma}
\proof
Let $(m,q)=(6,2)$. Then $A$ and $B$ are listed as follows:
\[\tiny
\begin{array}{l|l|l|l|l|l}
X^\infty& &A& |G:A|&B&|G:B|\\ \hline
\Omega^-_{10}(2)& \SO^-_{10}(2)&\Sp_{10}(2)& 2^5\cdot 3^2\cdot 7& 2^{15}{:}\GL_6(2)& 3^4\cdot5\cdot 11\cdot 17\\
\Omega^+_{10}(2)&\SO^+_{10}(2)&&&3{.}\PSU_6(2){.}[6]& 2^{14}\cdot 5\cdot 7\cdot 17\cdot 31\\
&\mbox{No}&&&\Sp_2(2)\times\Sp_6(2)& 2^{20}\cdot3^3\cdot 5\cdot 7\cdot 11\cdot 17\cdot 31\\
\Omega^-_{10}(2)&\SO^-_{10}(2)&&&\SL_6(2){.}2&2^{14}\cdot3^4\cdot5\cdot 11\cdot 17\\
&\mbox{No}&&&\Omega^+_6(4){.}2^2&  2^{16}\cdot3^4\cdot7\cdot 11\cdot 31\\ \hline
&\mbox{No}&(3\times \Omega^-_{10}(2)){.}2& 2^9\cdot3\cdot7\cdot31&\SL_6(2){.}2&2^{14}\cdot3^4\cdot5\cdot 11\cdot 17\\
&\mbox{No}&& & 2^{15}{:}\GL_6(2)& 3^4\cdot5\cdot 11\cdot 17\\  \hline
&\mbox{No}&2^{10}{:}\Omega^+_{10}(2)& 3^3\cdot7\cdot 11& 3{.}\PSU_6(2){.}[6]& 2^{14}\cdot 5\cdot 7\cdot 17\cdot 31
\end{array}
\]
If $A^\infty \le \Rad(A)X$ then we have $X\cong \PSp_{10}(2)$, $\Omega^-_{10}(2)$ or $\Omega^+_{10}(2)$.
Suppose that $A^\infty\not\le \Rad(A)X$. Then $\Rad(A)X\cap A^\infty$ is properly contained in $A^\infty$, and $X^\infty< A^\infty$.
It is easily shown that every odd prime divisor is also a divisor of $|\Rad(A)X\cap A^\infty|$. Check the maximal subgroups of  $A^\infty$ with order divisible by these primes, refer to \cite[Tables 8.64-8.69]{Low}. Combining
\cite[Table 10.4]{transubgroup}, we list $X^\infty$ and
those maximal subgroups of $A$ which possibly contain $X^\infty$ in the first and second columns of the above table, respectively. Then $X\cong \Omega^-_{10}(2)$ or $\Omega^+_{10}(2)$. Therefore, $X\cong \PSp_{10}(2)$ or $\Omega^\pm_{10}(2)$, and so $T=X<A$ and $G_\a\le B$.

Assume that $T\cong \Omega^-_{10}(2)$. Then $|G:T|=2^{10}\cdot3^2\cdot7\cdot 31$, which is a divisor of $|G_\a|$ and $|B|$. It follows that $B\cong 2^{15}{:}\GL_6(2)$ or $\SL_6(2){.}2$. Since $2^6-1$ is not a prime, by Lemma \ref{stab-lem-insol}, $\SL_6(2)\not\preccurlyeq G_\a$. Thus $G_\a\Rad(B)/\Rad(B)$ is a proper subgroup of $B/\Rad(B)$. It is easy to see that $|G_\a\Rad(B)/\Rad(B)|$ is divisible by $3^2\cdot7\cdot 31$. By \cite[Table 10.4]{transubgroup}, $G_\a\Rad(B)/\Rad(B)\lesssim \SL_5(2)$ or
$2^5{:}\SL_5(2)$.
By the Atlas \cite{Atlas}, $\SL_5(2)$ has
no proper insolvable subgroup with order divisible by $7\cdot 31$.
It follows that $\SL_5(2)\lesssim G_\a\Rad(B)/\Rad(B)$.
Thus we have $\SL_5(2)\lesssim G_\a\lesssim [2^{20}]{:}\SL_5(2)$.

 Assume next that $T\cong \Sp_{10}(2)$ or $\Omega^+_{10}(2)$. Then $|G:T|=2^{5}\cdot3^2\cdot7$ or $2^{10}\cdot3^3\cdot7\cdot 11$, respectively. Suppose that $11\in \pi(G_\a)$.
 Then $7,\,11\in \pi(G_\a)$.  Checking those $B$ with order divisible by $77$, we have $B\cong 3{.}\PSU_6(2){.}[6]$.
By \cite[Table 10.4]{transubgroup}, $G_\a\Rad(B)/\Rad(B)$ has a normal subgroup $\PSU_6(2)$ or $\M_{22}$, which is impossible by Lemma \ref{stab-lem-insol}. Thus $11\not\in \pi(G_\a)$.
Recalling that $|G_\a|$ is divisible by $|G:T|$, we have $T\cong \Sp_{10}(2)$. If $31\in \pi(G_\a)$  then $\SL_5(2)\lesssim G_\a\lesssim [2^{20}]{:}\SL_5(2)$ by a similar argument as in above paragraph. Thus we assume that neither $11$ nor $31$ is a divisor of $|G_\a|$. Moreover, for  $B\cong \Omega^+_6(4){.}2^2\cong \SL_4(4){.}2^2$,
by \cite[Table 10.3]{transubgroup},   $B$ has no proper subgroup of order divisible by $7\cdot 17$, and thus $17\not\in \pi(G_\a)$. Then the prime divisors of $|G_\a|$ are $2$, $3$, $7$ or and $5$.
By Lemma \ref{stab-lem-insol}, we conclude that $(G_\a\Rad(B)/\Rad(B))^\infty$ is a proper subgroup of $(B/\Rad(B))^\infty$. Checking the maximal subgroups of $\SL_6(2)$, $\PSU_6(2)$ and $\SL_4(4)$ (refer to \cite{Atlas} and \cite[Tables 8.3 and 8.4]{Low}), we conclude that $(G_\a\Rad(B)/\Rad(B))^\infty\succcurlyeq\SL_3(2)$ or $\A_7$. For $\A_7$, we have $G_\a\lesssim \Sy_7$ by Lemma \ref{stab-lem-insol}, and thus $|G_\a|$ is indivisible by $2^5$, a contradiction.
Thus $G_\a\succcurlyeq\SL_3(2)$.
If $\O_2(G_\a)=1$ then $|G_\a|_2\le 2^3$ by Lemma \ref{stab-lem-insol}, and so
$|G_\a|_2|T|_2\le 2^{28}<|G|_2$, a contradiction.  Thus part (2) of this lemma follows.
\qed

\begin{lemma}\label{O+2m-A-N1}
Let $(m,q)\ne(6,2)$ and $A=\sN_1$. Then $T\le A$ and $G_\a\le B$.
\end{lemma}
\proof
By the assumption,
$|G:A|={1\over (2,q-1)}q^{m-1}(q^m-1)$, and $|G:B|$ is listed as follows:
\[\tiny
\begin{array}{l|l|l}
B&|G:B|& \\ \hline
\sP_m \mbox{ or } \sP_{m-1}& \prod_{i=1}^{m-1}(q^i+1)& \\ \hline
{}^{\bf\hat{}}\GU_m(q){.}2&  {1\over 2}q^{m(m-1)\over 2}\prod_{i=1}^{m-1}(q^i-(-1)^{i-1}) & m\mbox{ even}\\ \hline

(\PSp_2(q){\times}\PSp_{m}(q)).c& {1\over c( q^2-1)}q^{(3m+2)(m-2)\over 4}\prod_{i={m\over 2}+1}^{m-1}(q^{2i}-1) &  m\mbox{ even}, c\le 2, c=2\mbox{ for odd $q$ and } 4\div m\\ \hline

{}^{\bf\hat{}}\GL_m(q){.}2& {1\over 2}q^{m(m-1)\over 2}\prod_{i=1}^{m-1}(q^i+1) &  m\mbox{ even} \\ \hline
\Omega^+_m(4).2^2 & q^{(m+2)(m-2)\over 2}\prod_{i=1}^{m\over 2}(q^{2(2i-1)}-1) & q=2,  m\mbox{ even}
\end{array}
\]
Clearly, $|G:B|$ is divisible by $\Phi^*_{2(m-1)}(q)$ or $\Phi^*_{2(m-2)}(q)$.
Moreover, $A'=\soc(A)\cong \Omega_{2m-1}$, $|A:A'|\le 2$, and $|A:A'|=2$ if and only if both   $q$ and ${m(q-1)\over 2}$ are odd.

 Suppose that $T\le B$. Then $|G:T|=|G:B||B:T|$, and thus
 $|G:T|$ is divisible by ${1\over (2,p)}q^{m(m-1)\over 2}$.
  Since    $G=TG_\a$, by Lemma \ref{Order},
  $|G_\a|_p\ge {1\over (2,p)}q^{m(m-1)\over 2}$, and $|G_\a|$ is divisible by $\Phi^*_{2(m-1)}(q)$ or $\Phi^*_{2(m-2)}(q)$.
 Noting $G=G_\a B$  and $A=G_\a(A\cap B)$ is core-free, next we deduce the contradiction  in two cases.

{\bf Case 1}. Let $q$ be odd. By Lemma \ref{m-factor}, choose a core-free maximal factorization $A_0=A_1B_1$ such that $A_0=X_1B_1$ with
 $A_0=A$ or $A'$,
 $G_\a^\infty \le X_1=A_1\cap G_\a$ and $|G_\a:X_1|$ is divisor of $|A:A_0|$.
Noting that $2m-1\ge 9$, by \cite[Theorem A]{max-factor},
$A_1$ and $B_1$ are given in \cite[Tables 1 and 2]{max-factor}.
On the other hand, it is easy to see  $|X_1|$ is divisible by $\Phi^*_{2(m-1)}(q)$ or $\Phi^*_{2(m-2)}(q)$. Checking the orders of $A_1$ and $B_1$, we conclude that $A_1\cong \Omega^-_{2(m-1)}(q){.}[2|A_0:A'|]$.
Note that $\Omega^-_{2(m-1)}(q)\not\cong S \preccurlyeq
X_1$.
By  \cite[Table 10.1]{transubgroup}, $|X_1|$ and so $|G_\a|$ is indivisible by
 $\Phi^*_{2(m-1)}(q)\Phi^*_{2(m-2)}(q)\Phi^*_{2(m-3)}(q)$.
Then $|G:B|$ is indivisible by $(q^{m-1}+1)(q^{m-2}+1)(q^{m-3}+1)$. It follows that
$B\cong {}^{\bf\hat{}}\GU_m(q){.}2$ for some even $m$, and then $|G_\a|$ divisible by  $q^{m(m-1)\over 2}\Phi^*_{2(m-2)f}(p)$.
Since $m$ is even, we have $A=A'=A_0$, and so $G_\a=X_1\le A_1 \cong \Omega^-_{2(m-1)}(q){.}2$.

It is easily shown that $S\preccurlyeq G_\a\cap A_1'$, and
$|G_\a\cap A'_1|$  is divisible by  $q^{m(m-1)\over 2}\Phi^*_{2(m-2)f}(p)$.
Combining  Lemma \ref{stab-lem-insol}, it follows from \cite[Theorem 3.1]{Bam-Pen}  that $G_\a\cap A'_1\lesssim \GO^+_2(q)\times \GO^-_{2m-4}(q)$, $[q^{2m-4}].(q-1\times \GO^-_{2m-4}(q))$,
$\GO_{2m-3}(q)$, $\GO^-_{2{m-1\over d}}(q^d){.}[d]$, $\GU_{m-1}(q){.}2$ or $\GO_{m-1}(q^2){.}2$, where $d>1$ is a divisor of $m-1$. The last three groups are excluded as their orders are indivisible by $q^{m(m-1)\over 2}$.
If $G_\a\cap A'_1\lesssim \GO_{2m-3}(q)$ then, by Lemma \ref{linear-stab-2},
we have $(S,2m-3,q)=(\A_7,7,3)$ or $(\A_7,7,5)$, which yields that $|G_\a|$ is indivisible by  $q^{m(m-1)\over 2}$, a contradiction. Therefore, $G_\a\cap A'_1\lesssim\GO^+_2(q)\times \GO^-_{2m-4}(q)$ or $[q^{2m-4}].(q-1\times \GO^-_{2m-4}(q))$.
It follows from  Lemma \ref{linear-stab-1} that $S\cong \PSL_2(q^{m-2})$, and so $r=q^{m-2}+1$. By Lemma \ref{r}, $q$ is even, a contradiction.

{\bf Case 2}. Let $q$ be even. Then $A\cong \Sp_{2m-2}(q)$.
Since $m\ge 5$, we have $\Omega^-_{2(m-1)}(q)\not\preccurlyeq G_\a$. Suppose that $(m,q)=(5,2)$. Then
$|G:A|=2^4\cdot 31$, and $T\le B\cong 2^{10}{:}\GL_5(2)$.
Recalling $G=TA$, by Lemma \ref{Order}, $|T|$ is divisible $31$, yielding $T=\GL_5(2)$. Thus $|G:T|=2^{10}\cdot3^3\cdot5\cdot17$, which is a divisor of $|G_\a|$. Checking the maximal subgroups of $\Sp_{8}(2)$, we conclude that $G_\a\lesssim \SO^-_8(2)$.  Checking the maximal subgroups of
$\SO^-_8(2)$, we get $\Omega^-_8(2)\unlhd G_\a$, which is impossible by Lemma  \ref{stab-lem-insol}. Then $m\ge7$ if $q=2$.

Assume that $|G_\a|$ is divisible by $\Phi^*_{2m-2}(q)$.
Recall that $|G_\a|_2\ge {1\over 2}q^{m(m-1)\over 2}$.
Then $|G_\a|$ is divisible by ${1\over 2}q^{m(m-1)\over 2}\Phi^*_{2m-2}(q)$. Applying Lemma \ref{linear-stab-1} to $A$ and $G_\a$,  either $S\cong \PSL_2(q^{m-1})$ and $|G_\a|_2\le (m-1)|S|_2$ with $(m-1)f$ a power of $2$, or $(S,2m-2,q)$ is one of
 $(\A_{13}, 12,2)$, $(\A_{19}, 18,2)$
 and $(\PSL_3(3), 12,2)$.
The former case implies that ${1\over 2}q^{m(m-1)\over 2}\le (m-1)q^{m-1}$, which is impossible as $m\ge 5$. For the latter case, we have $|G_\a|_2<2^{m-2}$ by Lemma \ref{stab-lem-insol}, again a contradiction.

Assume that $|G_\a|$ is indivisible by $\Phi^*_{2m-2}(q)$. Then
$B=^{\bf\hat{}}\GU_m(q){.}2$, $m$ is even and $|G_\a|$ is divisible by ${1\over 2}q^{m(m-1)\over 2}\Phi^*_{2m-4}(q)$.  Applying \cite[Theorem 3.1]{Bam-Pen} to $A$ and $G_\a$,  we have $G_\a\le A_1\cong [q^{2m-3}]{:}((q-1){\times}\Sp_{2m-4}(q))$ or $\Sp_2(q)\times \Sp_{2m-4}(q)$.
If $S\lesssim\Sp_2(q)$ then, by Lemma \ref{stab-lem-insol}, $|G_\a|_2$ is a divisor of $q^4f^2$, and so ${1\over 2}q^{m(m-1)\over 2}\le (m-1)q^{m-1}\le |G_\a|_2\le q^5$,
yielding $m\le 4$, a contradiction.
Let $N\unlhd A_1$ with $A_1/N\cong \Sp_{2m-4}(q))$.
Then $S\preccurlyeq G_\a N/N$, and $|G_\a N/N|$ is divisible by $\Phi^*_{2m-4}(q)$.
By Lemma \ref{linear-stab-1}, either $S\cong \PSL_2(q^{m-2})$ and $|G_\a N/N|_2\le (m-2)|S|_2$ with $(m-2)f$ a power of $2$, or $(S,2m-4,q)$ is one of
 $(\A_{13}, 12,2)$  and $(\PSL_3(3), 12,2)$.
 For the latter case, $m=8$ and $|G_\a|_2\ge 2^{27}$; however, by Lemma \ref{stab-lem-insol}, $|G_\a|_2\le 2^{10}$, a contradiction. The former case implies that ${1\over 2}q^{m(m-1)\over 2}\le |G_\a|_2\le |G_\a N/N|_2|N|_2\le (m-2)q^{3m-5}$, yielding $(m,q)=(6,2)$, again a contradiction.
\qed

\begin{lemma}\label{O+2m-A-N1-B}
Let $(m,q)\ne(6,2)$ and $A=\sN_1$.
Then  one of the following holds.
\begin{itemize}
\item[(1)] $T\cong \Sp_{2m-2}(q)$,  $\SL_2(q^{m\over 2})\unlhd G_\a^{\Ga(\a)}$,  $\O_2(G_\a)\ne 1$, $m\ge 8$, $p=2$ and $mf$ is a power of $2$.
\item[(2)] $T\cong \Sp_{2m-2}(q)$ or $\Omega^-_{2m-2}(q)$, $\PSL_d(q^{m\over d})\unlhd G_\a^{\Ga(\a)}$, $p=2$, $d$ is an odd prime  and $mf$ is a power of $d$.
\item[(3)] $T\cong \Omega_{2m-1}(q)$, $\PSL_d(q^{m\over d})\unlhd G_\a^{\Ga(\a)}$,  $p$ and $d$ are odd, and $mf$ is a power of $d$.
\end{itemize}
\end{lemma}
\proof
By Lemma \ref{O+2m-A-N1}, $T\le A$ and $G_\a\le B$, and then Lemma \ref{linear-stab-6} works.
Note that $|G:A|={1\over (2,q-1)}q^{m-1}(q^m-1)$. In particular,
$|G_\a|_p\ge q^{m-1}$. For those triples in Lemma \ref{linear-stab-6}~(2), by Lemma \ref{stab-lem-insol}, we have $|G_\a|_p\le q^{m-2}$, and thus they are excluded.
Thus $S\cong \PSL_d(q^{m\over d})$, $d$ is a prime and $mf$ is a power of $d$. In particular, $m$ is even if and only if $d=p=2$.

Assume that $d=2$. Then $p=2$ and $m$ is a power of $2$, and so $A\cong \Sp_{2m-2}(q)$ and $m\ge 8$.
 Noting that $|G_\a|$ is a divisor
of $|\O_p(G_\a)|\GL_{2}(q^{m\over 2})|{m^2f^2\over 4}$, it follows that $(|G_\a|,\Phi^*_{2m-2}(q)\Phi^*_{2m-4}(q)\Phi^*_{m-1}(q))=1$.
Since $G=TG_\a=TB$, by Lemma \ref{Order}, $|T|$ is divisible by $\Phi^*_{2m-2}(q)\Phi^*_{2m-4}(q)\Phi^*_{m-1}(q)$. By \cite[Table 10.1]{transubgroup}, we have $T=A\cong \Sp_{2m-2}(q)$.
Suppose that $\O_p(G_\a)=1$. Then $|G_\a|_p$ is a divisor of ${mf\over 2}q^{m\over 2}$. It follows that $q^{m-1}\le |G_\a|_p\le {mf\over 2}q^{m\over 2}<q^{{3m\over 4}+1}$, yielding $m<8$, a contradiction. Thus $\O_p(G_\a)\ne 1$. Then part (1) of this lemma follows.

Assume   $d$ is odd.
Then  $m$ is odd, and  $B=\sP_m$ or $\sP_{m-1}$. Suppose that $q$ is even. Then $A\cong \Sp_{2m-2}(q)$.
Noting that $|T|$ is divisible by $\Phi^*_{2m-2}(q)\Phi^*_{2m-4}(q)\Phi^*_{2m-6}(q)$, by \cite[Table 10.1]{transubgroup}, we have $T=\Sp_{2m-2}(q)$ or $\Omega^-_{2m-2}(q)$, and part (2) of this lemma follows.
Suppose that $q$ is odd. Then $A\cong \Omega_{2m-1}(q){.}{4\over (4,m(q-1))}$.
Noting that $|G_\a|$ is a divisor
of $|\O_p(G_\a)|\GL_{d}(q^{m\over d})|{m^2f^2\over d^2}$,
it follows that $|B:G_\a|$ is divisible by $\prod_{m>i\ge 2, (i,d)=1}\Phi^*_i(q)$. In particular, $|B:G_\a|$ is divisible by
$\Phi^*_{m-1}(q)$. Noting that $|G:G_\a|=|G:B||B:G_\a|$, it follows that $|G:G_\a|$ is divisible by $\Phi^*_{m-1}(q)\prod_{i=1}^{m-1}(q^i+1)$.
Then $|T|$ is divisible by $\Phi^*_{m-1}(q)\prod_{i=1}^{m-1}(q^i+1)$, and thus
$|T|$ has a divisor $(\Phi^*_{m-1}(q))^2\Phi^*_{2m-2}(q)\Phi^*_{2m-4}(q)$.
Since $m\ge 5$ and $q$ is odd, $\Phi^*_{m-1}(q)\ne 1$.
By \cite[Table 10.2]{transubgroup},
$T\cong \Omega_{2m-1}$, desired as in part (3) of this lemma.
\qed

\begin{lemma}\label{O+2m-A-N2}
Let $(m,q)\ne(6,2)$ and $A=\sN_2^-$. Then
 $T\cong \POmega^-_{2m-2}(q)$,  $\PSL_m(q)\unlhd G_\a^{\Ga(\a)}$,  $m$ is an odd   prime  and $f$ is a power of $m$.
\end{lemma}
\proof
In this case, $|G:A|=q^{2m-2}{(q^m-1)(q^{m-1}-1)\over 2(q+1)}$,
and $B$ is $\sP_m$, $\sP_{m-1}$ or ${}^{\bf\hat{}}\GL_m(2){.}2$ (with $q=2$ and $m$ even). Then $|G:B|$ is divisible by $\prod_{i=1}^{m-1}(q^i+1)$.
Let $\overline{A}=A/\Rad(A)$ and $\overline{X}=X\Rad(A)/\Rad(A)$. It is easily shown that
$\soc(\overline{A})\cong \POmega^-_{2m-2}(q)$ and $|\overline{X}|$ is divisible by $\Phi^*_{2m-2}(q)\Phi^*_{2m-4}(q)\Phi^*_{2m-6}(q)$.
Noting that $(m,q)\ne(6,2)$, it follows from \cite[Table 10.1]{transubgroup} that either $\POmega^-_{2m-2}(q)\lesssim \overline{X}$, or $(m-1,q)=(4,2)$ and $\Omega^-_8(2)\not\lesssim \overline{X}$. Suppose  that the latter case occurs. Checking the maximal subgroups of $\Omega^-_8(2)$, we get $\overline{X}\lesssim \PSL_2(16){.}2$. Since $\PSL_2(16){.}2$ has no proper subgroup of order divisible by $\Phi^*_{8}(2)\Phi^*_{4}(2)$, we have $\PSL_2(16)\lesssim \overline{X}$. If $T=X\cong \PSL_2(16)$ then $|T|$ is indivisible by $|G:B|$, a contradiction. Thus
$X=G_\a$, and $T\le B$, yielding $T\lesssim \SL_5(2)$.
By   Lemma \ref{stab-lem}, it is easy to show that $|G_\a|_3\le 9$, and so $|T|_3|G_\a|_3\le 3^4<3^5=|G|_3$, a contradiction. Therefore, we have $\POmega^-_{2m-2}(q)\lesssim \overline{X}$. This implies that $A\ge T\cong \POmega^-_{2m-2}(q)$. Then our result follows from Lemma \ref{linear-stab-5}.
\qed

\begin{theorem}\label{O+2m}
Assume that {\rm Hypothesis \ref{hypo-graphs-classical}} holds and $G=\POmega_{2m}^+(q)$, where $m\ge 5$ and $q=p^f$. Then one of the following holds.
\begin{itemize}
\item[(1)] $(m,q)=(6,2)$, $T=\Omega^-_{10}(2)$ or $\Sp_{10}(2)$, and $\SL_5(2)\le G_\a\le [2^{20}]{:}\SL_5(2)$;
\item[(2)] $(m,q)=(6,2)$, $T=\Sp_{10}(2)$, $G_\a^{\Ga(\a)}\cong \SL_3(2)$ and $\O_2(G_\a)\ne 1$.

\item[(3)] $T\cong \Sp_{2m-2}(q)$,  $\SL_2(q^{m\over 2})\unlhd G_\a^{\Ga(\a)}$,  $\O_2(G_\a)\ne 1$, $m\ge 8$, $p=2$ and $mf$ is a power of $2$.
\item[(4)] $T\cong \Sp_{2m-2}(q)$ or $\Omega^-_{2m-2}(q)$, $\PSL_d(q^{m\over d})\unlhd G_\a^{\Ga(\a)}$, $p=2$, $d$ is and odd prime  and $mf$ is a power of $d$.
\item[(5)] $T\cong \Omega_{2m-1}(q)$, $\PSL_d(q^{m\over d})\unlhd G_\a^{\Ga(\a)}$,  $p$ and $d$ are odd, and $mf$ is a power of $d$.

\item[(6)] $T=\POmega^-_{2m-2}(q)$,  $\PSL_m(q)\unlhd G_\a^{\Ga(\a)}$,  $m$ is an odd   prime  and $f$ is a power of $m$.
\end{itemize}
\end{theorem}
\proof
By Lemmas \ref{O+12(2)}, \ref{O+2m-A-N1-B} and \ref{O+2m-A-N2}, we assume that $A=\sP_1$ and $(m,q)\ne (6,2)$.
In this case, $m$ is even, $|G:A|={(q^m-1)(q^{m-1}+1)\over (q-1)}$, $B={}^{\bf\hat{}}\GU_m(q){.}2$ and $|G:B|={1\over 2}q^{m(m-1)\over 2}\prod_{i=1}^{m-1} (q^i+(-1)^i)$. We next deduce a contradiction.

Suppose that $X=G_\a$.
Let $\overline{A}=A/\Rad(A)$ and $\overline{X}=X\Rad(A)/\Rad(A)$. Then $|\Rad(A)|_p=q^{2m-2}$,
$\soc(\overline{A})\cong \POmega^+_{2m-2}(q)$, $|\overline{A}:\soc(\overline{A})|\le 2$ and $|\overline{X}|$ is divisible by ${1\over 2}q^{(m-4)(m-1)\over 2}\prod_{i=2}^{m-1} (q^i+(-1)^i)$.
By Lemma \ref{sub-factor}, $\overline{A}=\overline{A}\,\overline{A\cap B}$ is core-free, where $\overline{A\cap B}=(A\cap B)\Rad(A)/\Rad(A)$.
By Lemma \ref{m-factor}, we have a core-free maximal factorization $G_0=A_0B_0$ such that $\overline{A}\ge G_0=X_0B_0$, $X_0=A_0\cap \overline{X}$, $\soc(G_0)\cong \POmega^+_{2m-2}(q)$, $\overline{X}^\infty\le X_0$ and $|X_0|$ has a divisor ${1\over o}|\overline{X}|$, where $o\le 2$.
In particular, $|X_0|$ is divisible by ${1\over (4,p^2)}q^{(m-4)(m-1)\over 2}\Phi^*_{2m-4}(q)\Phi^*_{m-1}(q)\Phi^*_{m-3}(q)$.
Noting that $m-1$ is odd and $|A_0|$ is divisible by $\Phi^*_{2m-4}(q)$, by \cite[Theorem A and Proposition 2.5]{max-factor}, $A_0$ has type $\sN_1$ or $\sN^-_2$.
 However, in this case, $|A_0|$ is indivisible by $\Phi^*_{m-1}(q)$, a contradiction.

Now, we have $T\le A$ and $G_\a\le B\cong {}^{\bf\hat{}}\GU_m(q){.}2$.
Then $|G_\a|$ is divisible by $|G:A|$; in particular,  $|G_\a|$ is divisible by $\Phi^*_m(q)\Phi^*_{2m-2}(q)$. Applying Lemma \ref{linear-stab-2} to the preimages of $G_\a\cap B'$ and $B'$ in $\GL_m(q^2)$, we conclude that $(m,q^2)=(4,9),\,(4,25)$ or $(6,4)$. Then $m<5$ or $(m,q)=(6,2)$, a contradiction.
\qed

\vskip 20pt

\subsection{Symplectic case}
Let $G=\PSp_{2m}(q)$, where $m\ge2$, $(m,q)\ne (2,2)$ and $q=p^f$ for some  prime $p$. Then $A$ and $B$ are described as in \cite[Tables 1-3]{max-factor}.

\begin{lemma}\label{Sp-1}
$A$ and $B$ are not in \cite[Table 3]{max-factor}.
\end{lemma}
\proof
Suppose that $A$ and $B$ are described as in \cite[Table 3]{max-factor}. Since $A$ and $B$ are insolvable, they are listed as follows:
\[\tiny
\begin{array}{l|l|l|l|l}
G&A&|G:A|&B&|G:B|\\ \hline
\PSp_6(3)&\PSL_2(13)& 2^7\cdot3^8\cdot5& [3^5]{:}2\PSU_4(2)& 2^2\cdot7\cdot13\\  \hline
\PSp_8(2)&\Omega^-_8(2){.}2&2^3\cdot3\cdot5&\Sy_{10}& 2^8\cdot 3\cdot17\\
&\PSL_2(17)& 2^{12}\cdot3^3\cdot5^2\cdot7&\Omega^+_8(2){.}2&2^3\cdot17
\end{array}
\]

Assume that $G=\PSp_6(3)$. Then $X=A\cong \PSL_2(13)$. Let $\overline{B}=B/\Rad(B)$ and $\overline{Y}=Y\Rad(B)/\Rad(B)$.
Then $\overline{B}\cong \PSU_4(2)$ and $|\overline{Y}|$ has a divisor $2^6\cdot3^3\cdot5$. Searching the subgroups of $\PSU_4(2)$, we conclude that $\overline{Y}\succcurlyeq\PSU_4(2)$ or $\Omega_5(3)$.
By Lemma \ref{stab-lem-insol}, neither $X$ nor $Y$ plays the role of $G_\a$, a contradiction.
Thus   $G=\PSp_8(2)$.

Suppose that $A\cong \PSL_2(17)$ and $B\cong \Omega^+_8(2){.}2$.
It is easily   shown that $X=A\cong \PSL_2(17)$, and $|Y\cap B'|$ has a divisor $5^2\cdot 7$. By \cite[Table 10.4]{transubgroup}, $Y\cap B'\cong \Omega^+_8(2)$. Then $X\ne G_\a$ and $Y\ne G_\a$, a contradiction.
Suppose that $A\cong \Omega^-_8(2){.}2$.  Then
$X\cong \Omega^-_8(2)$ due to searching the subgroups of $\Omega^-_8(2)$ with order divisible by $2^7\cdot 3\cdot17$.
Thus $X=T$, and $G_\a=Y\le B$,
which contradicts Lemma \ref{linear-stab-7}.
\qed

\begin{lemma}\label{Sp-2}
If $A$ and $B$ are as  in \cite[Table 2]{max-factor}, then $p=2$, $f>2$ is a power of $2$, $G=\Sp_6(q)$, $T\cong \G_2(q)$, $\SL_2(q^2)\unlhd G_\a^{\Ga(\a)}$  and  $\O_2(G_\a)\ne 1$.
\end{lemma}
\proof
Suppose that $A$ and $B$ are described as in \cite[Table 2]{max-factor}. Then $q$ is even, and $(G,A)$ is one of $(\Sp_4(q),\Sz(q))$ and $(\Sp_6(q),\G_2(q))$.

Assume that $G=\Sp_4(q)$ and $A\cong \Sz(q)$. Then $f>1$ is odd. Noting that $A=X(A\cap B)$, by \cite{Factor-excep}, we have $A=X$. Thus $T=X\cong \Sz(q)$, and $G_\a=Y\le B\cong \GO^+_4(q)$. Since $G=TG_\a$ and $|G:T|=q^2(q^2-1)(q+1)$, by Lemma \ref{Order}, $|G_\a|$ divisible by  $q^2(q^2-1)(q+1)$. Then $|B:G_\a|$ is a divisor of $2(q-1)$. This implies that $G_\a$ has a normal subgroup $\PSL_2(q)\times\PSL_2(q)$, which is impossible by Lemma \ref{stab-lem-insol}.

Assume that $G=\Sp_6(q)$ and $A\cong \G_2(q)$. Then
$|G:A|=q^3(q^4-1)$ and $B\in \{\GO^\pm_6(q), \sP_1, \sN_2\}$. Suppose that $X\ne A$. Noting that $A=X(A\cap B)$, by \cite{Factor-excep}, $q=4$ and $X\cong\J_2$, $\SU_3(4)$ or  $\SU_3(4){.}2$.
Then $T=X\cong\J_2$ or $\PSU_3(4)$, and $G_\a=Y\le B$. Noting that $|T|$ is divisible by $|G:B|$, the factorization $G=TB$ is listed as follows:
\[\tiny
\begin{array}{l|l|l|l}
T&|G:T|&B&|G:B|\\ \hline
\J_2& 2^{11}\cdot 3\cdot5\cdot13\cdot17& \PSU_4(4){.}2& 2^5\cdot3^3\cdot7\\ \hline
\PSU_3(4)&  2^{12}\cdot 3^3\cdot5\cdot7\cdot17&\PSL_4(4){.}2&
2^5\cdot 5\cdot 13\\
\end{array}
\]
Clearly, $S\preccurlyeq Y\cap B'$, and then $\pi(Y)=\pi(S)=\pi(Y\cap B')$. Then, by \cite[Table 10.3]{transubgroup}, we get
$Y\cap B'\cong \PSU_4(4)$ or $\PSL_4(4)$. Thus $G_\a=Y\succcurlyeq\PSU_4(4)$ or $\PSL_4(4)$, a contradiction.
 Therefore, $X=A\cong \G_2(q)$, and so $T\cong \G_2(q)$ for $q>2$.

Let $\overline{B}=B^\infty\Rad(B)/\Rad(B)$ and $\overline{Y}=(Y\cap B^\infty)\Rad(B)/\Rad(B)$.
Then $|\overline{Y}|$ has a divisor $\Phi^*_4(q)$, and $S\preccurlyeq \overline{Y}$.
If $\overline{B}\cong \PSL_4(q)$, $\PSp_4(q)$ or $\PSp_2(q)\times\PSp_4(q)$ then, by Lemma \ref{linear-stab-1}, we conclude that
$S\cong \PSL_2(q^2)$.
Assume that $\overline{B}\cong \PSU_4(q)$. Then $|\overline{Y}|$ has a divisor ${1\over 2}q^3\Phi^*_4(q)$. Checking the maximal subgroups $\PSU_4(q)$, we have $\overline{Y}\lesssim q^4{:}([a]{.}\PSL_2(q^2){.}[b])$, where $ab={(q-1)(2,q^2-1)\over (4,q+1)}$.  This implies that $S\cong \PSL_2(q^2)$. All in all,
$r=q^2+1$, and so $p=2$ and $f$ is a power of $2$ by Lemma \ref{r}. If $\O_2(G_\a)=1$ then $q^3\le |G_\a|_2\le q^22f$,
yielding $q=4$. For $q=4$, confirmed by GAP, there exists no a desired graph. Then  the lemma follows.
\qed

\vskip 10pt

Assume that $A$ and $B$ are as  in \cite[Table 1]{max-factor}.
For $(m,q)=(3,2)$ or $(2,4)$, $A$ and $B$ are  listed as follows:
 \[\tiny
\begin{array}{l|l|l|l|l}
(m,q)&A&|G:A|&B&|G:B|\\ \hline
(3,2)& \Sp_2(8){.}3& 2^6\cdot3\cdot5&2^5{:}\Sy_6&3^2\cdot7\\
& & & \Sy_8& 2^2\cdot3^2\\
& & & \Omega^-_6(2){.}2& 2^2\cdot7\\
& \Omega^-_6(2){.}2& 2^2\cdot7& 2^6{:}\PSL_3(2)&3^3\cdot 5 \\
& & & \Sy_8& 2^2\cdot3^2\\ \hline
(2,4)& \Sp_2(16){.}2&2^3\cdot3 \cdot 5& 2^6{:}(3\times \Sp_2(4))& 5\cdot 17\\
&&& \Omega^+_4(4){.}2& 2^3\cdot 17\\
&&& \Omega^-_4(4){.}2& 2^3\cdot 3\cdot 5\\

&\Omega^-_4(4){.}2& 2^3\cdot 3\cdot 5&  2^6{:}\GL_2(4)& 5\cdot 17\\
&&& (\Sp_2(4)\times\Sp_2(4)){:}2& 2^3\cdot 17
\end{array}
\]
We first search the insolvable subgroups of $A$ and $B$ with orders divisible by $|G:B|$ and $|G:A|$, respectively.
It follows that $X'=A'$, and either $Y'=B'$ or
$G_\a=Y$, and then all possible pairs $(T,G_\a)$ are easily determined. Calculation shows that no desired graph arises for $(m,q)=(3,2)$ or $(4,4)$.

\begin{lemma}\label{Sp-m=small-(m,q)}
Let $A$ and $B$ be as   in \cite[Table 1]{max-factor}.
Then $(m,q)\ne (3,2),\,(2,4)$.
\end{lemma}

For convenience, we list
 $A$, $B$ and their indices in $G$  as follows:
 \[\tiny
\begin{array}{l|l|l|l|l}
A&|G:A|&B&|G:B|&\\ \hline
\PSp_{2a}(q^b){.}b&{1\over b}q^{m(m-a)}\prod_{(i,b)=1}(q^{2i}-1)& \sP_1& {q^{2m}-1\over q-1}&\\

&& \GO^\pm_{2m}(q)& q^m{q^{m}\pm1\over 2}& p=2\\
&&\sN_2&q^{2(m-1)}{q^{2m}-1\over q^2-1} & b=q=2, 2\div m\\

\hline

\GO^-_{2m}(q)& q^m{q^{m}-1\over 2}& \sP_m&\prod_{i=1}^m(q^i+1)&p=2\\
& & \Sp_m(q)\wr\Sy_2& {1\over 2}q^{m^2\over 2}\prod_{i=1}^{m\over 2}{q^{m+2i}-1\over q^{2i}-1}&2\div m, p=2\\

& & \GO^+_{2m}(2)&q^m{q^{m}+1\over 2}&q=2

\end{array}
\]
 where $b$ is a prime such that $m=ab$.

\begin{lemma}\label{Sp-m=2}
Let $A$ and $B$ be  as in \cite[Table 1]{max-factor}, and  $m=2$. Then $p=2$, $f>2$ is a power of $2$, $T\cong\PSL_2(q^2)$, and either $G_\a\cong \PSL_2(q^2)$  or $\SL_2(q)\unlhd G_\a^{\Ga(\a)}$ and $\O_2(G_\a)\ne 1$.
\end{lemma}
\proof
Note that $A'\cong \PSL_2(q^2)$ and   $q>2$.
Assume that $B\cong \GO^-_4(q)$. Then $A\cong \Sp_2(q^2){.}2$,
$B' \cong \PSL_2(q^2)$ for even  $q$, and
both $|X\cap A'|$ and $|Y\cap B'|$ are divisible by ${1\over 4}q^2(q^2-1)$. Checking the insolvable subgroups of $\PSL_2(q^2)$, it follows that either $X\cap A'\cong Y\cap B'\cong \PSL_2(q^2)$, or  one of $X\cap A'$ and $Y\cap B'$ is isomorphic to $\A_5$ if $q=4$.
By Lemma \ref{Sp-m=small-(m,q)}, $q\ne 4$.
Thus $X\cap A'\cong Y\cap B'\cong \PSL_2(q^2)$. It follows that $T\cong \PSL_2(q^2)$, $G_\a\cong \PSL_2(q^2){.}o$  and $r=q^2+1$, where $o\le 2$.
Then, by Lemma \ref{r}, $f$ is a power of $2$. In particular, $f\ge 4$.
Suppose that $G_\a\cong  \PSL_2(q^2){.}2$.
Then $|T:(T\cap G_\a)|=|G:G_\a|={1\over 2}q^2(q^2-1)$, and thus
$T_\a=T\cap G_\a\cong \D_{2r}$, a maximal subgroup of $T$.
Thus both  $T$ and $G$ are primitive  on $V$, and $T_\a$ and $G_\a$
have a common orbit on $V$ of size $r$, which is impossible by \cite[Theorem 1 and Table 7.1]{common-suborbit}.
Thus $G_\a\cong \PSL_2(q^2)$.

Assume next that $B\not\cong \GO^-_4(q)$. Then  $|X\cap A'|$ is divisible by $\Phi^*_2(q^2)$. Note that $X\cap A'$ is insoluble. Checking the subgroups of $\PSL_2(q^2)$, it follows that either $X\cap A'=A'$, or $q=p\equiv\pm 3\,(\mod 10)$ and $X\cap A'\cong\A_5$.
For the latter, we have $\Phi^*_2(q^2)=5$ or $15$, which is impossible by \cite[Theorem 3.9]{Hering}. Thus $X=A'$ or $A$.

Let $\overline{B}=B^\infty\Rad(B)/\Rad(B)$ and $\overline{Y}=(Y\cap B^\infty)\Rad(B)/\Rad(B)$. Then $\overline{B}\cong \PSL_2(q)$ or $\PSL_2(q)\times\PSL_2(q)$, $\overline{Y}$ is insolvable and has order divisible by $\Phi^*_2(q)$. Clearly, if $q\ne p$ then $\Phi^*_2(q)\ne 1$.
Checking the insolvable subgroups of $\PSL_2(q)$,
we conclude that either $\PSL_2(q)\preccurlyeq \overline{Y}$, or $\A_5\preccurlyeq \overline{Y}$ and $q=p\equiv\pm 1\,(\mod 10)$. Then $Y\succcurlyeq\PSL_2(q)$ or $\A_5$.

If $T=Y$, recalling $T\not\cong \A_5$, then $q^4=|G|_p\le |T|_p|G_\a|=|T|_p|X|\le q(2,p)q^2$, which is impossible as $q>2$.
Thus we have $Y=G_\a$, $X=T$ and $S=\soc(G_\a^{\Ga(\a)})\cong \PSL_2(q)$ or $\A_5$. If $S\cong \A_5$, recalling  $q=p\equiv\pm 1\,(\mod 10)$, then $r=5$ and $|G_\a|_p=1$ by Lemma \ref{stab-lem-insol}, and so $q^4=|G|_p\le |T|_p|G_\a|_p=q^2$, a contradiction.
If $q=11$ then $r=11$ and  $|G_\a|_p=q$, so  $q^4=|G|_p\le |T|_p|G_\a|_p=q^3$,  a contradiction.
Therefore,   $r=q+1$, $p=2$ and $f$ is power of $2$.
If $\O_2(G_\a)=1$ then  $q^4=|G|_2\le |G_\a|_2|T|_2\le |\GammaL_2(q)|_2|T|_2=q^3f$, yielding $q\le f$, a contradiction.
Thus $\O_2(G_\a)\ne 1$, and  the lemma follows.
\qed

\begin{lemma}\label{Sp-A-O^-}
Let $A$ and $B$ be as   in \cite[Table 1]{max-factor}, and $m\ge 3$.  Assume that $A\cong \GO^-_{2m}(q)$ for even $q$. Then
one of the following holds.
\begin{itemize}
\item[(1)]  $T\cong \Omega^-_{2m}(q)$ and $S\cong \PSL_d(q^{m\over d})$, $d$ is an odd prime and $mf$ is a power of $d$.
\item[(2)]  $T\cong \Omega^-_{2m}(q)$ and $S\cong \PSL_2(q^{m\over 2})$,  $\O_2(G_\a)\ne1 $ and $mf$ is a power of $2$.
\item[(3)] $q=2$, $T\cong \Omega^+_{2m}(2)$, $S\cong \PSL_2(2^m)$ and $\O_2(G_\a)=1$, where $m$  is a power of $2$.
\item[(4)] $G=\Sp_{12}(2)$, $T\cong \Omega^-_{12}(2)$,  either $S\cong \PSL_3(2)$ and $\O_2(G_\a)\ne 1$ or $S\cong \PSL_5(2)$.
\end{itemize}
\end{lemma}
\proof
Recall that $\{X,Y\}=\{T,G_\a\}$ with $X\le A$ and $Y\le B$.
We discuss in three cases according the structure of $B$.

{\bf Case 1}.
 Let $B=\sP_m$. Then $|G:B|$ has a divisor $\Phi^*_{2m}(q)\Phi^*_{2m-2}(q)\Phi^*_{2m-4}(q)$. Since $|X|$ is divisible by $|G:B|$, it is easily shown that
$|X\cap A'|$ is divisible by $\Phi^*_{2m}(q)\Phi^*_{2m-2}(q)\Phi^*_{2m-4}(q)$. By \cite[Tables 10.1 and 10.3]{transubgroup}, either $q=2$ and $m\in \{4,5\}$ or
$X\cap A'=A'\cong \Omega^-_{2m}(q)$. For the former,
$|G:B|=3^3\cdot5\cdot 17$ or $3^4\cdot5\cdot 11\cdot 17$, by
the Atlas \cite{Atlas} and \cite[Table  10.4]{transubgroup}, also $X\cap A'=A'\cong \Omega^-_{2m}(q)$. Since $m\ge 3$,
it follows that $T=X\cong \Omega^-_{2m}(q)$, and $G_\a=Y\le B$.
 Then $S$ is described as in Lemma \ref{linear-stab-7}.

If $S\cong \PSL_d(q^{m\over d})$ for some prime $d$ then
part (1) or (2) of this lemma occurs.
For the triples given in Lemma \ref{linear-stab-7}~(2) and (3), by Lemma \ref{stab-lem-insol} and calculating $|G_\a|_2$, we conclude that $(S,m,q)$ is one of   $(\A_7,4,2)$, $(\PSL_3(2),6,2)$ and $(\PSL_5(2),6,2)$.
Suppose that $(S,m,q)=(\A_7,4,2)$. By Lemma \ref{stab-lem-insol}, noting $|G|_2\le |T|_2|G_\a|_2$, we have $G_\a\cong \Sy_7$.
Then $\Sy_7\cong G_\a\O_2(B)/\O_2(B)\le \overline{B}\cong \A_8$, which is impossible.
Thus $T$ and $S$ are desired as in part (4) of this lemma.


{\bf Case 2}.  Let $B\cong\Sp_m(q)\wr\Sy_2$.
Suppose that $T=Y\le B$. Then $G_\a=X\le A$, $T\lesssim \Sp_m(q)$, and $|G:T|$
is divisible by $q^{3m^2\over4}\prod_{i={m\over 2}+1}^{m}(q^{2i}-1)$. Since $|G_\a|$ is divisible by $|G:T|$,
it is easily shown that
$|G_\a\cap A'|$ is divisible ${1\over 2}q^{3m^2\over4}\prod_{i={m\over 2}+1}^{m}(q^{2i}-1)$. Checking \cite[Tables 10.1, 10.3 and 10.4]{transubgroup}, we conclude that
$G_\a\cap A'=A'\cong \Omega^-_{2m}(q)$, which is impossible by Lemma \ref{stab-lem-insol}. Then $T=X\lesssim\Omega^-_{2m}(q)$ and $G_\a=Y\le B$.
By a similar argument as in Case 1, $T$ and $S$ are described as in  part (2) or (4) of this lemma.

{\bf Case 3}.  Let $B\cong\GO^+_{2m}(2)$, $m\ge 4$ and $q=2$.
If $G_\a=Y\le B$ then Lemma \ref{linear-stab-7} works, and $T$ and $S$ are described as in  part (2) or (4) of this lemma.

Suppose next that $G_\a=X\le A$. Then $T\le B'\cong \Omega^+_{2m}(2)$, and
$|G:T|$ is divisible by $2^m(2^m+1)$.
Considering $G_\a\cap A'$, by Lemma \ref{linear-stab-1}, either $S\cong \PSL_2(2^m)$ with $m$ a power of $2$, or $(S,2m,q)$ is one of $(\A_{11},10,2)$, $(\A_{13},12,2)$, $(\A_{19},18,2)$, $(\M_{11},10,2)$, $(\PSL_2(11),10,2)$ and $(\PSL_3(3),12,2)$. The triples $(\M_{11},10,2)$ and $(\PSL_2(11),10,2)$ are excluded by calculation of $|G_\a|_2$. If $(S,2m,q)=(\PSL_3(3),12,2)$ then
$5\not\in\pi(G_\a)$  by Lemma \ref{stab-lem-insol},
which contradicts that $|G_\a|$ is divisible by $2^m+1$.
Assume that $(S,2m,q)=(\A_{11},10,2)$, $(\A_{13},12,2)$ or $(\A_{19},18,2)$.
By \cite[pp. 52, 3.2.4(e)]{max-factor}, $A\cap B\cong \Sp_{2m-2}(2)\times 2$, which  is maximal in both $A$ and $B$.
Since $A=G_\a(A\cap B)$, we have a core-free maximal factorization of $A$ with a factor $\Sp_{2m-2}(2)\times 2$ and the other factor containing $\A_{2m+1}$, which leads to a contradiction by checking the maximal factorizations of $\GO^-_{2m}(2)$.

Now $S\cong \PSL_2(2^m)$, where $m$ is a power of $2$.
By Lemma \ref{linear-stab-1}, $|G_\a\cap A'|_2$ is a divisor of $m|S|_2$.
Then $|G_\a|$ is a divisor of $2m2^m$.
If $\O_2(G_\a)\ne 1$ then $|G_\a|_2\ge 2^{2m}$, and so $2^m\le 2m$, yielding $m\le 2$, a contradiction. Thus  $\O_2(G_\a)=1$.
It is easily shown that $|G:G_\a|$ has a divisor ${2^{m(m-1)}\over 4m(2^m-1)}\prod_{i=1}^{m-1}(2^{2i}-1)$, and then $|T|$ is divisible by ${2^{m(m-1)}\over 4m(2^m-1)}\prod_{i=1}^{m-1}(2^{2i}-1)$.
If $m>4$ then, by \cite[Table 10.1]{transubgroup}, $T\cong \Sp_{2m-2}(2)$ or $\Omega^+_{2m}(2)$.
Assume that $m=4$. Then $|T|$ is divisible by $2^8\cdot3^3\cdot7$. Checking the subgroups of $\Omega^+_8(2)$, we have  $T\cong \Omega^+_8(2)$ or $\Sp_6(2)$.
Suppose that $T\cong \Sp_{2m-2}(2)$.
Then $|G:T|$ is divisible by $2^{2m-1}$, which is a divisor of $|G_\a|$. Since $\O_2(G_\a)=1$, we have $|G_\a|_2\le 2^mmf$.
It follows that $2^{2m-1}\le 2^mm$, and so $2^{m-1}\le m$, which is impossible as $m>2$. Thus   $T\cong \Omega^+_{2m}(2)$, and  part (3) of this lemma occurs.
\qed

\begin{lemma}\label{Sp-A=Sp-2a-B-ne-O-}
Let  $m\ge3$ and  $A\cong \PSp_{2a}(q^b){.}b$, where $b$ is a prime and $m=ab$.  Assume that $B\not \cong \GO^-_{2m}(q)$. Then $p=2$ and  one of the following holds.
\begin{itemize}
\item[(1)]    $T\cong \Omega^+_{2m}(q)$, $S\cong \PSL_2(q^{m})$,  $\O_2(G_\a)=1 $ and $mf$ is a power of $2$.

 \item[(2)]  $m=3$, $q>2$,   $T\cong \PSL_2(q^3)$, $S\cong\SL_2(q^2)$, $\O_2(G_\a)\ne 1$ and $f$ is a power of $2$.

\end{itemize}
\end{lemma}
\proof
Assume that $X=G_\a$. Then $T\le B$, and so $T\lesssim\PSp_{2m-2}(q)$ or $\Omega^+_{2m}(q)$. In particular, $|G:T|$ has a divisor $q^m(q^m+1)$. Then $|X\cap A'|$ is divisible by ${1\over (b,p)}q^m\Phi^*_{2m}(q)$. Since $(m,q)\ne (3,2)$, we have $\Phi^*_{2m}(q)\ne 1$.
Clearly, $S\preccurlyeq X\cap A'$.
Applying Lemma \ref{linear-stab-1} to
$X\cap A'$ and $A'$, we conclude
$S\cong \PSL_d(q^{b{2a\over d}})$, where $d$ is a prime divisor of $2a$ and $mf$ is power of $d$. By Lemma \ref{r},
we have $d=b=p=2$. Again by Lemma \ref{linear-stab-1},  $|X\cap A'|_2$ is a divisor of $q^ma$, and so $|G_\a|_2$ is a divisor of $q^mm$. Thus $q^mm\ge |G:T|_2$ as $|G_\a|$ is divisible by $|G:T|$. It follows that $T\le B\cong \GO^+_{2m}(q)$.
It is easily shown that $|G:G_\a|$ is divisible by
${1\over m}q^{m^2-m}\Phi^*_{2m-2}(q)\Phi^*_{2m-4}(q)\Phi^*_{m-1}(q)$, which is a divisor of $|T|$. By \cite[Table 10.1]{transubgroup},
we conclude that either $T=B'\cong  \Omega^+_{2m}(q)$, or $(m,q)=(4,2)$. For the latter case, noting that $|T|$ is divisible by $2^{10}\cdot 5\cdot 7$, we have $T\cong \Omega^+_{8}(2)$ by \cite[Table 10.4]{transubgroup}.
Then part (1) of this lemma occurs.

We assume next that $X=T$. Then $T\lesssim \PSp_{2a}(q^b)$, $G_\a=Y\le B$, and $|G_\a|$ is divisible by $q^{m(m-a)}\prod_{(i,b)=1}(q^{2i}-1)$.
Let $\overline{B}=B'\Rad(B)/\Rad(B)$ and $\overline{Y}=(Y\cap B')\Rad(B)/\Rad(B)$. Then  $S\preccurlyeq  \overline{Y}$ and $|\overline{Y}|$ is divisible by $\prod_{(i,b)=1,i>1}(q^{2i}-1)$. In particular, $|\overline{Y}|$ is divisible by
$\prod_{(i,b)=1,i>1}\Phi^*_{2i}(q)$.

Suppose that $B=\sP_1$ or $\sN_2$. Then $\overline{B}\cong \PSp_{2m-2}(q)$.
By  Lemma \ref{linear-stab-1}, either $S\cong \PSL_2(q^{m-1})$,
$p=2$ and $(m-1)f$ is a power of $2$, or $(S,2(m-1),q)$ is one of the triples
$(\A_5,4,3)$, $(\A_7,4,7)$, $(\A_7,6,3)$, $(\A_7,6,5)$, $(\A_{11},10,2)$, $(\A_{13},12,2)$, $(\A_{19},18,2)$, $(\M_{11},10,2)$, $(\PSL_3(2),6,3)$, $(\PSL_3(2),6,5)$, $(\PSL_2(11),10,2)$, $(\PSL_3(3),12,2)$.
Recalling that $|G_\a|$ is divisible by $q^{m(m-a)}\prod_{(i,b)=1}(q^{2i}-1)$, the latter case is excluded by Lemma \ref{stab-lem-insol} and calculating $|G_\a|_p$ or $\pi(G_\a)$.
Then $S\cong \PSL_2(q^{m-1})$. Assume $m>3$. Then  $|G:G_\a|$ has a divisor $\prod_{i={m+1\over 2}}^{m-2}\Phi^*_{2i}(q)$, and so
$|T|$ is divisible by $\prod_{i={m+1\over 2}}^{m-2}\Phi^*_{2i}(q)$. Considering $\pi(\PSp_{2a}(q^b))$, we have $b\ge {m+1\over 2}$, yielding $b=m$ and  $T\lesssim \PSp_{2}(q^m)$. Since $\Phi^*_{2m-4}(q)$ is  a divisor of $|T|$, it implies that $(m,q)=(5,2)$ and $S\cong \PSL_2(16)$. In this case, $q^{m(m-a)}\prod_{(i,b)=1}(q^{2i}-1)$ has a divisor $7$, and so $7\in \pi(G_\a)=\pi(S)=\{2,3,5,17\}$, a contradiction.
Thus $m=3$,
$T\lesssim\PSL_2(q^3)$. Then $T\cong\PSL_2(q^3)$ by checking the subgroups of $\PSL_2(q^3)$. If $\O_2(G_\a)=1$ then
$q^6 \le |G_\a|_2\le |\PGammaL_2(q^2)|_2$, yielding $q^4\le 2f$, a contradiction. Therefore, part (2) of this lemma occurs.

Now let $B\cong \GO^+_{2m}(q)$. Recall that $|G_\a|$ is divisible by $q^{m(m-a)}\prod_{(i,b)=1}(q^{2i}-1)$. Assume that $(m,q)=(4,2)$. Then $|G_\a|$ is divisible by $2^8\cdot 3^3\cdot 7$. Checking the maximal subgroups of $\GO^+_{8}(2)$, we have $G_\a\lesssim \Sp_6(2)\times 2$. Searching subgroups of $\Sp_6(2)$ with order divisible by $3^3\cdot 7$, we conclude that $\PSU_3(3)\preccurlyeq G_\a$, a contradiction.
Then  $(m,q)\ne (4,2)$, and  $\Phi^*_{(2m-2)f}(q)\ne 1$.
Recalling the limitations on $G_\a$, by \cite[Theorem 3.1]{Bam-Pen}, we conclude that $G_\a\lesssim  \GO^-_2(q)\times\GO^-_{2m-2}(q)$, or $\Sp_{2m-2}(q)$.
Then Lemma \ref{linear-stab-1} works for
$G_\a$ modulo a suitable solvable group.
Arguing similarly as above, we have  part (2) of this lemma.
\qed

\begin{theorem}\label{Sp-2m}
Assume that {\rm Hypothesis \ref{hypo-graphs-classical}} holds and $G=\PSp_{2m}(q)$, where $m\ge 2$ and $q=p^f$. Then $q$ is even and one of the following holds.
\begin{itemize}

\item[(1)] $m=2$, $f>2$ is a power of $2$, $T\cong\PSL_2(q^2)$, and either $G_\a\cong \PSL_2(q^2)$  or $\SL_2(q)\unlhd G_\a^{\Ga(\a)}$ and $\O_2(G_\a)\ne 1$.

\item[(2)]
$m=3$, $q>2$, $f$ is a power of $2$,   $T\cong \G_2(q)$ {\rm (}with $f>2${\rm )} or $\PSL_2(q^3)$, $\SL_2(q^2)\unlhd G_\a^{\Ga(\a)}$  and  $\O_2(G_\a)\ne 1$.

\item[(3)] $(m,q)=(6,2)$,  $T\cong \Omega^-_{12}(2)$,  either $S\cong \PSL_3(2)$ and $\O_2(G_\a)\ne 1$ or $S\cong \PSL_5(2)$.

\item[(4)] $m\ge 5$, $T\cong \Omega^-_{2m}(q)$, $S\cong \PSL_d(q^{m\over d})$, $d$ is an odd prime and $mf$ is a power of $d$.
\item[(5)] $m\ge 4$, $T\cong \Omega^-_{2m}(q)$ and $S\cong \PSL_2(q^{m\over 2})$,  $\O_2(G_\a)\ne1 $ and $mf$ is a power of $2$.
\item[(6)]
 $T\cong \Omega^+_{2m}(q)$, $S\cong \PSL_2(q^{m})$,  $\O_2(G_\a)=1$, $m\ge 4$ and $mf$ is a power of $2$.

\end{itemize}
\end{theorem}
\proof
Based on the previous discussion, we assume that
$A\cong \PSp_{2a}(q){.}b$ and $B\cong \GO^-_{2m}(q)$, where $q$ is even, $m\ge 3$, $(m,q)\ne (3,2)$, and $b$ is a prime divisor of $m$ with $m=ab$.
If $G_\a\le B$ then   part (2) of this lemma occurs by a similar argument as  in the last paragraph of the proof of Lemma \ref{Sp-A=Sp-2a-B-ne-O-}.

Now Let $G_\a\le A$ and $T \le B$. By Lemma \ref{linear-stab-7} and Case 4 in its proof, we conclude that $a\ge 3$,
$S\cong \PSL_d(q^{m\over d})$, $d$ is a prime divisor of $a$ and $mf$ is a power of $d$. Assume that $m$ is even. Then $d=b=2$, $a\ge 4$ and $S\cong \PSL_2(q^{m\over 2})$. It is easily shown that $|G:G_\a|$ is divisible by $\prod_{i={m+2\over 2}}^{m}\Phi^*_{2i}(q)$. By \cite[Table 10.1]{transubgroup}, $T\cong \Omega^-_{2m}(q)$, and then $T$ and $S$ are desired as in part (5) of this lemma. Now let $m$ be odd; in particular, $a,\,b$ and $d$ are odd.
Then $A\cap B\cong \GO_{2a}^-(q^b){.}b$, which is maximal in $B$. Suppose that $T\not\cong \Omega_{2m}^-(q)$.
Noting $B=X(A\cap B)$, by Lemma \ref{sub-factor}, we have a core-free maximal factorization $G_0=A_0B_0$ with $\soc(G_0)\cong \Omega_{2m}^-(q)$ and $\Omega^-_{2a}(q^b)\lesssim B_0$.
It follows from \cite[Theorem A]{max-factor} that $q\in \{2,4\}$ and $b=2$, a contradiction. Thus $T \cong \Omega_{2m}^-(q)$, and $T$ and $S$ are desired as in part (4) of this lemma.
\qed

\subsection{Linear case} Let $G=\PSL_{n}(q)$, where $n\ge2$ and $q=p^f$ for some  prime $p$. Since $A$ and $B$ are insolvable, by \cite[Theorem A]{max-factor} and Corollary \ref{prime-dim}, $n$ is not a prime, $A\cong {}^{\bf\hat{}}\GL_a(q^b){.}b$ or $\PSp_n(q){.}[c]$, and $B\cong [q^{n-1}]{:}({q-1\over (n,q-1)}{.}\PGL_{n-1}(q))$, where $a\ge 2$, $b$ is a prime with $n=ab$, and $c={(2,q-1)({n\over 2},q-1)\over (n,q-1)}$.
Recalling $G\not\cong \A_8$, we have $(m,q)\ne (4,2)$.

\begin{lemma}\label{L-n-X=Ga}
If $G_\a\le A$ then $T\cong \PSL_{n-1}(q)$, $S\cong \PSL_d(q^{n\over d})$, where $d$ is an odd prime and $nf$ is a power of $d$.
\end{lemma}
\proof
Let $G_\a\le A$. Then $T\lesssim\PSL_{n-1}(q)$, and $|G:T|={(n-1,q-1)\over (n,q-1)}q^{n-1}(q^n-1)$.
Let $X_0=G_\a\cap A'$. Then $S\preccurlyeq X_0$ and $|X_0|$ is divisible by ${1\over (e,p)}q^{n-1}(q^n-1)$, where $e=b$ or $c$.
Suppose that $(n,q)=(6,2)$. Then $A'\cong \GL_3(4)$ or $\Sp_6(2)$, and $|G_\a|$ is divisible by $2^4\cdot3^2\cdot7$.
 Checking the subgroups of $A'$, we conclude that
  $X_0$ has a composition factor $\PSL_3(4)$, $\PSU_3(3)$, $\A_8$ or $\Sp_6(2)$, which is impossible by Lemma \ref{stab-lem-insol}. Thus $(n,q)\ne (6,2)$, and $\Phi^*_{nf}(p)\ne 1$. Then Lemma \ref{linear-stab-1} works here.
By calculation of $\pi(S)$ or $|G_\a|_p$, we conclude that $A\cong {}^{\bf\hat{}}\GL_a(q^b){.}b$,  $S\cong \PSL_d(q^{n\over d})$, $d$ is an odd prime divisor of $a$ and $nf$ is a power of $d$.
In particular, $a\ge 3$ and $n$ is odd.
Then $|G:A|$ has a divisor $(q^{n-1}-1)(q^{n-2}-1)$. By \cite[Table 10.1]{transubgroup}, we have $T\cong \PSL_{n-1}(q)$, and the lemma holds.
\qed

\begin{lemma}\label{L-n-Y=Ga}
If  $G_\a\le B$ then one of the following holds.
\begin{itemize}
\item[(1)] $(n,q)=(4,4)$, $T\cong \Sp_4(4)$, $S\cong \PSL_3(2)$ and $\O_2(G_\a)\ne 1$.
\item[(2)] $S\cong \PSL_{n-1}(q)$, $n-1$ is an odd prime and $f$ is a power of $n-1$.
\end{itemize}
\end{lemma}
\proof
Let $G_\a\le B$, set
$\overline{B}=B/\Rad(B)$ and $\overline{Y}=Y\Rad(B)/\Rad(B)$. Then $S\preccurlyeq \overline{Y}$. Note that $|G:A|$ has a divisor $q^{n-1}-1\over q-1$, which is also a divisor of $|G_\a|$. Then $|\overline{Y}|$ is divisible by $q^{n-1}-1\over q-1$, and hence divisible by $\Phi^*_{(n-1)f}(p)$.

Assume that $(n,q)=(4,4)$. Then $\overline{B}\cong \PGL_3(4)$.
Checking the subgroups of $\PGL_3(4)$, we conclude that $S\cong \overline{Y}\cong \PSL_3(2)$. Note that $T\lesssim\PSL_2(16)$ or $\Sp_4(4)$. Checking the subgroups of $\PSL_2(16)$ and $\Sp_4(4)$, we conclude that $T\cong \PSL_2(16)$ or $\Sp_4(4)$.
By Lemma \ref{stab-lem-insol}, $|G_\a|_5=1$. If $T\cong \PSL_2(16)$ then $5^2=|G|_5\le |T|_5|G_\a|_5=5$, a contradiction.
Then $T\cong \Sp_4(4)$, and $|G:T|=2^4\cdot 3^2\cdot 7$. By
Lemma \ref{stab-lem-insol}, $\O_2(G_\a)\ne 1$, and part (1) of this lemma follows.

Now $(n,q)\ne (4,4)$, and $\Phi^*_{(n-1)f}(p)\ne 1$. By Lemma \ref{linear-stab-1}, one of the following holds:
\begin{itemize}
\item[(i)]
$S\cong \PSL_d(q^{n-1\over d})$, $d$ is a prime and $(n-1)f$ is a power of $d$;
\item[(ii)]  $(S,n-1,q)$ is one of $(\A_7, 3,25)$, $(\PSL_3(2), 3,25)$, $(\PSL_3(2), 3,9)$ and $(\PSL_2(11), 5,4)$.
\end{itemize}
The case (ii) is easily excluded by calculating $|G_\a|_p$ and $|G:A|_p$. Thus (i) occurs.

Suppose that $A\cong \PSp_n(q){.}[c]$. Then $T\lesssim\PSp_n(q)$, $n$ is even
 and $|G:T|={(2,q-1)\over (n,q-1)}q^{n(n-2)\over 4}\prod_{i=2}^{n\over 2}(q^{2i-1}-1)$. In particular, $|G_\a|$ is divisible by $\Phi^*_{n-3}(q)$. Clearly, $\Phi^*_{n-3}(q)\ne 1$. Since $\pi(G_\a)=\pi(S)$, we have
 $(\Phi^*_{n-3}(q),|S|)\ne 1$. By Lemma \ref{phi},
  $n-3$ is a divisor
 of $(n-1)i\over d$ for some $2\le i\le d$. It follows that $d=n-1$, and so $S\cong \PSL_{n-1}(q)$.

Suppose that $A\cong {}^{\bf\hat{}}\GL_a(q^b){.}b$.
If $b=2$ or $b\ge 5$ then  $|G_\a|$ is divisible by $\Phi^*_{n-3}(q)$, and so $d=n-1$ by a similar argument as above. Now let $b=3$. Then $|G_\a|$ is divisible by $\Phi^*_{n-2}(q)$.
Clearly, $n\ne 8$ and $n-2>2$. Then $\Phi^*_{n-2}(q)\ne 1$, this also yields $d=n-1$. Thus the lemma follows.
\qed

\begin{theorem}\label{L-n}
Assume that {\rm Hypothesis \ref{hypo-graphs-classical}} holds and $G=\PSL_{n}(q)$, where $n\ge 4$ and $q=p^f$. Then   one of the following holds.
\begin{itemize}
\item[(1)]
$T\cong \PSL_{n-1}(q)$, $S\cong \PSL_d(q^{n\over d})$, $d$ is an odd prime and $nf$ is a power of $d$.
\item[(2)] $(n,q)=(4,4)$, $T\cong \Sp_4(4)$, $S\cong \PSL_3(2)$ and $\O_2(G_\a)\ne 1$.
\item[(3)] $S\cong \PSL_{n-1}(q)$, $n-1$ is an odd prime, $f$ is a power of $n-1$, and
\begin{itemize}
\item[] $T\cong \PSL_{m}(q^{n\over m})$, $m$ is a proper divisor of $n$; or
\item[] $T\cong \PSp_{m}(q^{n\over m})$, $m\ge 4$ is  even, and $q$ is even if $m=n$; or
\item[] $T\cong \POmega^-_{m}(q^{n\over m})$, $4\div m$, $m\ge 8$, and $q$ is even if $m=n$; or
\item[]  $T\cong  \G_2(q^{n\over 6})$, $6\div n$, $q$ is even; or
 \item[] $T\cong \PSU_3(3)$ and $(n,q)=(6,2)$.
\end{itemize}
\end{itemize}
\end{theorem}
\proof
By Lemmas \ref{L-n-X=Ga} and \ref{L-n-Y=Ga}, we assume that
$T\le A$, $S\cong \PSL_{n-1}(q)$ and $f$ is a power of $n-1$.
We shall show that $T$ is described as in part (3) of lemma.

Since $G=TB$, by Lemma \ref{Order}, $|T|$ is divisible by $q^n-1\over q-1$.
Assume that $(n,q)=(6,2)$. Then $T\lesssim \PSL_2(8)$, $\PSL_3(4)$ or $\Sp_6(2)$. Checking the subgroups of $\PSL_2(8)$, $\PSL_3(4)$ and $\Sp_6(2)$, it follows that
$T\cong \PSL_2(8),\, \PSU_3(3)$ or $\Sp_6(2)$, desired as in part (3). Thus we assume next that $(n,q)\ne (6,2)$. Then $\Phi^*_{nf}(p)\ne 1$.

{\bf Case 1}. Assume that $A\cong {}^{\bf\hat{}}\GL_a(q^b){.}b$. Then $T\lesssim\PSL_a(q^b)$.
Choose a minimal divisor $m$ of $n$ such that $T\lesssim\PSL_m(q_0)$, where $q_0=q^{n\over m}$. Then $m\le a$,   $|T|$ is divisible by $q_0^m-1\over q-1$ and $\Phi^*_{m{nf\over m}}(p)$.
If $m=2$ then $T\lesssim\PSL_2(q^{n\over 2})$ and $q_0$ is not a prime, and so $T\cong\PSL_2(q^{n\over 2})$ by checking the subgroups of $\PSL_2(q^{n\over 2})$. Thus we let $m\ge 3$, and so $a\ge 3$.
Noting that  $T$ is not an alternating group, by \cite[Theorem 3.1]{Bam-Pen}, one of the following holds:
\begin{itemize}
\item[(i)] $T\cong \POmega^-_m(q_0)$, $\PSp_m(q_0)$ or $\PSL_m(q_0)$;
\item[(ii)] $T\cong \PSU_{m}(\sqrt{q_0})$, $q_0$ is a square and  $m$ is odd;
\item[(iii)] $(T,m,q_0)$ is one of   the following triples:
 $(\J_3,9,4)$,
  $(\PSL_3(2),3,9)$, $(\PSL_3(2),3,25)$,
  $(\PSL_2(11),5,4)$,
 $(\PSL_2(13),6,4)$,  $(\PSL_2(19),9,4)$,
    $(\Suz(q_0),4,2^{nf\over m})$  and $(\G_2(q_0),6,2^{nf\over m})$.
\end{itemize}
We only need deal with cases (ii) and (iii).

Suppose that (ii) occurs. Since $m\ge 3$ is odd, by Lemma \ref{r},
$\Phi^*_{m}(\sqrt{q_0})\ne1$. Since $m$ is odd, it is easy to see that  $|\PSU_{m}(\sqrt{q_0})|$ is indivisible by $\Phi^*_{m}(\sqrt{q_0})$; however, $|T|$ is divisible by  $q_0^m-1\over q-1$, a contradiction.
Suppose that (iii) occurs. For the first six triples, we have $b={n\over m}=2$ and $q=p$, and then they are excluded by calculating $|T|$ and ${q^n-1\over q-1}$. Similarly, the triple
$(\Suz(q_0),4,2^{nf\over m})$ is excluded as $|\Suz(q_0)|$ is indivisible by $q_0+1$. Thus $T\cong \G_2(q_0)$, desired as in part (3).

{\bf Case 2}. Assume that $A\cong \PSp_n(q){.}[c]$. Then $T\lesssim\PSp_n(q)$ for even $n$.
Noting that  $T$ is not an alternating group, $n-1$ is an odd prime and $f$ is a power of $n-1$, by \cite[Theorem 3.1]{Bam-Pen}, one of the following holds:
\begin{itemize}
\item[(iv)] $T\cong \Omega^-_n(q)$ (with $q$ even) or $\PSp_n(q)$;
\item[(v)] $T\lesssim \PSp_m(q^{n\over m})$, where $m$ is a proper divisor of $n$;
\item[(vi)] $T\lesssim \PSU_{n\over 2}(q)$, where ${n\over 2}q$ is odd;
\item[(vii)] $(T,n,q)$ is one of the following triples: \\ $(\J_2,6,5)$, $(\PSL_3(2),6,3)$,  $(\PSL_3(2),6,5)$,  $(\PSL_2(13),6,3)$,  $(\PSL_2(13),6,4)$,  $(\PSL_2(17),8,2)$, \\ $(\PSL_2(25),12,2)$, $(\PSL_2(37),18,2)$, $(\PSL_2(41),20,2)$, $(\PSL_3(3),12,2)$, $(\PSL_3(4),6,3)$,\\  $(\PSp_4(5),12,2)$,
    $(\PSU_3(9),6,5)$, $(\Suz(q),4,2^f)$  and $(\G_2(q),6,2^f)$.
\end{itemize}
Recall that $|T|$ is divisible by $q^n-1\over q-1$.
For $T\lesssim \PSU_{n\over 2}(q)$, since ${n\over 2}q$ is odd, $\Phi^*_{n\over 2}(q)\ne 1$; however, it is easy to see that $|\PSU_{n\over 2}(q)|$ is indivisible by $\Phi^*_{n\over 2}(q)$, a contradiction.

For case (vii),
only the triples $(\PSL_2(13),6,3)$ and  $(\G_2(q),6,2^f)$  are left. If $(n,q)=(6,3)$ then $r={3^5-1\over 3-1}$, which  is not a prime, a contradiction.
Thus $T\cong \G_2(q)$, desired as in part (3) of this lemma.
For case (iv), either $T$ is   desired  in part (3) of this lemma, or $T\cong \Omega^-_n(q)$  with   odd ${n\over 2}$; the latter does occur as $|T|$ has a divisor $q^{n\over 2}-1\over q-1$.

The left case is that $T\lesssim \PSp_m(q^{n\over m})$. We choose $m$ minimal as possible.
If $m=2$  then $T\lesssim \PSp_2(q^{n\over 2})$, yielding $T\cong \PSp_2(q^{n\over 2})$, desired as in   part (3) of this lemma. Thus we may let $m\ge 4$.
Clearly $q^{n\over m}$ is not a prime, and $m-1$ is odd.
Using \cite[Theorem 3.1]{Bam-Pen}, by a similar argument as above,  $T$ is described as in  part  (3).
\qed

\vskip 20pt

\section{The proof of Theorem \ref{mainth-2}}\label{sect=cay}
Let $\Ga=(V,E)$ be a connected $X$-symmetric Cayley
graph of a finite simple group $T$
of prime valency $r\ge 7$, where $T\le X\le \Aut(\Ga)$. Assume that $T$ is not normal in $X$, and
$X_\a$ is insolvable, where $\a\in V$.  By Theorem \ref{mainth-1}~(2), we choose a characteristic subgroup $G$ of $X$ such that
that $G^\infty=G\gneq \Rad(G)T$, $G/\Rad(G)$ is simple, $r\not\in \pi(\Rad(G))$, $\Rad(G)$ is intransitive on $V$ and $\Ga$ is $G$-symmetric.
By Lemma \ref{N-p}, $\Rad(G)$ is semiregular on $V$, and so $|\Rad(G)|$ is a proper divisor of $|V|=|T|$.

Set $\overline{G}=G/\Rad(G)$ and $\overline{T}=T\Rad(G)/\Rad(G)$, and
identify $\overline{G}$ with a subgroup of
  $\Aut(\Ga_{G/\Rad(G)})$. Then
  $\Ga_{G/\Rad(G)}$ is $\overline{G}$-symmetric, and $\overline{G}=\overline{T}\,\overline{G}_{\overline{a}}$, where $\overline{\a}=\a^{\Rad(G)}$. Thus $(\overline{G},\overline{T}, \overline{G}_{\overline{a}})$ is one of the triples satisfying Hypothesis \ref{hypo-graphs}.
Since $G\unlhd X$, we have $G_\a\unlhd X_\a$, and so $G_\a^{\Ga(\a)}\unlhd X_\a^{\Ga(\a)}$. Since $X_\a$ is insolvable,  $X_\a^{\Ga(\a)}$ is almost simple, and then  $G_\a^{\Ga(\a)}$ is almost simple. In particular, $G_\a$ is insolvable.
Note that $\overline{G}_{\overline{a}}\cong G_{\overline{a}}/\Rad(G)=\Rad(G)G_\a/\Rad(G)\cong G_\a$. Then $\overline{G}_{\overline{a}}$ is insolvable. Thus $(\overline{G},\overline{T}, \overline{G}_{\overline{a}})$ can be read out from Theorems \ref{G-or-T-alte}, \ref{G-sporadic} and
\ref{U-2m}-\ref{L-n}.

Note that $T\cong \overline{T}$ and ${|T||G_\a|\over |\Rad(G)|}=|G:\Rad(G)|=|\overline{G}|=|\overline{T}|
|\overline{G}_{\overline{a}}:(\overline{G}_{\overline{a}}\cap \overline{T})|$. Then $|\Rad(G)|=|\overline{G}_{\overline{a}}\cap \overline{T}|={|G_\a|\over |\overline{G}:\overline{T}|}$, and so $|\overline{G}:\overline{T}|_r=r$. Thus we only need to deal with those triples $(\overline{G},\overline{T}, \overline{G}_{\overline{a}})$ with $|\overline{G}:\overline{T}|_r=r$.
For further argument, we need more precise structural information about $G$. Note that $G$ is perfect, i.e., $G'=G$.
\begin{lemma}\label{Schur}
Let $L=\overline{G}$ and $R_s=\min\{m\mid L\lesssim\PGL_m(q_0), q_0 \mbox{ a power of } s\}$, where $s$ is a  prime. Then either $G$ is a perfect central extension of $L$, or  $|\Rad(G)|_s\ge q_0^{R_s}$ for some prime $s$ and its a power  $q_0$ with $L\lesssim\PGL_m(q_0)$. For the latter case, if $L \cong\A_n$ for  $n\ge 9$ then $n$ is a power of $2$ and $|\Rad(G)|_2\ge 2^{n-2}$.
\end{lemma}
\proof
 Assume that $G$ is not a perfect central extension of $L$.
 Recalling $G^\infty=G$, it is easily shown that $G/K$ is perfect for each normal subgroup $K$ of $G$ with $K\le \Rad(G)$.
Choose a  minimal $K$ such that $G/K$ is a perfect central extension of $L$. Then $K\ne 1$. Choose $J\unlhd G$ with $J\le K$ such that
$K/J$ is a minimal normal subgroup of $G/J$. Since $K$ is solvable, we set
$K/J\cong \ZZ_s^l$ for some prime $s$ and integer $l\ge 1$.

Since $G/J=(K/J){.}(G/K)$ and $G/J$ is perfect, by the choice of $K$, $G/J=(K/J){.}(G/K)$ is not a central extension; otherwise, $G/J$ should be a perfect central extension of $L$, refer to \cite[pp. 167, (33.5)]{Aschbacher}. Thus $G/J$ acts nontrivially on $K/J$ induced by the conjugation of $G$ on $K$, and the kernel is $\C_{G/J}(K/J)$. Set $C/J=\C_{G/J}(K/J)$. Then $C  \le \Rad(G)$, and
$G/C\cong (G/J)/(C/J)\lesssim\GL_{l_0}(q_0)$, where $s^l=q_0^{l_0}$. It follows that $l_0\ge R_s$, refer to \cite[pp.184, Proposition 5.3.2]{KL-book}.
Then $|\Rad(G)|_s\ge q_0^{l_0}\ge q_0^{R_s}$.

Suppose that $L\cong \A_n$ for some $n\ge 9$. Then, by \cite[Propositions 5.3.2 and 5.3.7]{KL-book}, $l\ge R_s\ge n-2$. Thus $|\A_n|_s\ge s^{n-2}$, and then
$n-2+(2,s)\le
\sum_{i=1}^\infty\lfloor{n\over s^i}\rfloor$. The only possibility is that $s=2$ and $n$ is a power of $2$. Then the lemma follows.
\qed


We next finish the proof of Theorem \ref{mainth-2} by checking one by one the possible candidates for  $(\overline{G},\overline{T}, \overline{G}_{\overline{a}})$   described in Theorems \ref{G-or-T-alte}, \ref{G-sporadic} and
\ref{U-2m}-\ref{L-n}. Recall that $|\Rad(G)|=|\overline{G}_{\overline{a}}\cap \overline{T}|$ and $|\overline{G}:\overline{T}|_r=r$.

Assume that $L\cong \A_n$. Then, by Theorem \ref{G-or-T-alte}, either  $T\cong \A_{n-1}$ as at Line 1 of Table \ref{exceptions}, or
 $\overline{G}\cong \A_{r+1}$ for $r\ge7$, $\overline{G}_{\overline{a}}\cong\A_r$, $\Ga_{G/\Rad(G)}\cong \K_{r+1}$ and $T\cong\overline{T}$ is isomorphic to a transitive subgroup of $\A_{r+1}$. Since $|\overline{G}:\overline{T}|_r=r$, we have $r\not\in \pi(T)$. Checking the subgroups of $\A_8$, we conclude that $\A_8$ has no insolvable transitive subgroup with order indivisible by $7$. Thus we have $r>7$, and $r+1>9$. By Lemma \ref{Schur} either $G$ is a perfect central extension of $\A_{r+1}$, or $|\Rad(G)|_2\ge 2^{r-1}$. Recalling $|\Rad(G)|$ is a divisor of $|\overline{G}_{\overline{a}}|$, if $|\Rad(G)|_2\ge 2^{r-1}$ then
 $|\A_r|_2\ge 2^{r-1}$, which is impossible. Thus  $G$ is a perfect central extension of $\A_{r+1}$, and so $|\Rad(G)|\le 2$, refer to \cite[pp. 170, (33.15)]{Aschbacher}. If $|\Rad(G)|=2$ then $\Ga$ is bipartite by \cite[Theorem 1.1]{Du}, which is not the case. Then $\Rad(G)=1$, and
 $T$ is desired at Line $2$ of Table \ref{exceptions}.

Assume that $L$ is a sporadic simple group. By Theorem
\ref{G-sporadic}, $\overline{G}\cong \M_n$, $\overline{G}_{\overline{a}}\cong \M_{n-1}$, $\Ga_{G/\Rad(G)}\cong\K_n$, $n=12$ or $24$ and $\overline{T}\cong \A_5$ or $\PSL_3(2)$. In this case, $|\Rad(G)|=|\overline{G}_{\overline{a}}\cap \overline{T}|=s\in \{5,7\}$.
It follows that $G\cong \ZZ_s\times \M_n$, which contradicts the choice of $G$.

Assume that $L$ is a classical simple group. Then Theorems \ref{U-2m}-\ref{L-n} work here, and either the pair $(L,T)$ is desired as in Table \ref{exceptions} or $(\overline{G}, \overline{T},S,r)$ is listed as follows:
\[\tiny
\begin{array}{l|l|l|l|l}
|\overline{G}:\overline{T}|&\overline{G}& \overline{T}&\overline{G}_{\overline{\a}}  &r\\ \hline

3^7\cdot13&\POmega^+_8(3)& \Omega^+_8(2)&\overline{G}_{\overline{a}}\cong [3^6]{:}\SL_3(3)& 13\\ \hline
2^{10}\cdot3^2\cdot7\cdot31&\Omega^+_{12}(2)& \Omega^-_{10}(2)&  \SL_5(2)\preccurlyeq\overline{G}_{\overline{\a}}& 31\\ \hline

2^4\cdot3^2\cdot7&\PSL_{4}(4)& \Sp_{4}(4)&  \SL_3(2)\preccurlyeq\overline{G}_{\overline{\a}},\,\O_p(\overline{G}_{\overline{\a}})\ne 1& 7\\ \hline
2^{10}\cdot3\cdot5\cdot7\cdot31&\PSL_{6}(2)& \PSU_{3}(3)&  \SL_5(2)\preccurlyeq\overline{G}_{\overline{\a}}& 31\\ \hline
2^{mf}(2^{mf}+1)&\Sp_{2m}(2^f), m\ge 4&\Omega^+_{2m}(2^f)&  \SL_2(2^{mf})\lesssim\overline{G}_{\overline{a}}\lesssim \GammaL_2(2^{mf})& 2^{mf}+1\\
\end{array}
\]
Recall that $|\Rad(G)|=|\overline{G}_{\overline{a}}\cap \overline{T}|={|\overline{G}_{\overline{\a}}|\over |\overline{G}:\overline{T}|}$.

Suppose that $\overline{G}\cong \Sp_{2m}(2^f)$ and $\overline{T}\cong \Omega^+_{2m}(2^f)$, where $m\ge 4$ and $mf$ is a power of $2$.
Then $|\Rad(G)|$ is a divisor of $mf(2^{mf}-1)^2$ and divisible by $2^{mf}-1$. If $G$ is  a perfect central extension of $\Sp_{2m}(2^f)$ then $|\Rad(G)|\le (2, 2^f-1)=1$ by \cite[pp. 173, Theorem 5.1.4]{KL-book}, a contradiction. Thus, by Lemma \ref{Schur},
$|\Rad(G)|_s\ge q_0^{R_s}$, where $s\in \pi(\Rad(G))$,  $R_s$ is defined as in Lemma \ref{Schur}, $q_0$ is a power of $s$ and $q_0=2^f$ if $s=2$. Note that $R_s$ is given as in \cite[Theorem 5.3.9 and Proposition 5.4.13]{KL-book}. If $s=2$ then $R_s=2m$, and so
$mf\ge |\Rad(G)|_s\ge (2^f)^{R_s}=2^{2mf}$, which is impossible.
Thus $s$ is odd, and so $R_s\ge {1\over 2}2^{(m-1)f}(2^{(m-1)f}-1)(2^f-1)$ by \cite[Theorem 5.3.9]{KL-book}.
It follows that $\log_s(2^{mf}-1)^2\ge\log_s|\Rad(G)|_s\ge R_s\ge {1\over 2}2^{(m-1)f}(2^{(m-1)f}-1)(2^f-1)$.
Then $4\log_s(2^{mf}-1)\ge 2^{(m-1)f}(2^{(m-1)f}-1)(2^f-1)$, and so
$\log_s(2^{mf}-1)>2^{(m-1)f}-1$ as $m\ge 4$.
We have $mf>\log_s(2^{mf}-1)>2^{(m-1)f}-1$, yielding $mf+1>2^{(m-1)f}>2(m-1)f$. This implies that $m\le2$, a contradiction.

Suppose that $\overline{G}\cong \Omega^+_{12}(2)$ and $\overline{T}\cong \Omega^-_{10}(2)$. By Lemma \ref{stab-lem-insol}, either $|\O_2(\overline{G}_{\overline{\a}})|=2^{30}$
or $|\O_2(\overline{G}_{\overline{\a}})|\le 2^{10}$. Noting that $2^{30}=|\Omega^+_{12}(2)|\ge |\O_2(\overline{G}_{\overline{\a}})||\SL_5(2)|_2$, we have $|\O_2(\overline{G}_{\overline{\a}})|\le 2^{20}$, and so $|\O_2(\overline{G}_{\overline{\a}})|\le 2^{10}$. Then $|\overline{G}_{\overline{\a}}|$ is a divisor of $2^{10}|\SL_5(2)|$ and divisible by $|\SL_5(2)|$.
It follows that $|\Rad(G)|= {|\overline{G}_{\overline{\a}}|\over |\overline{G}:\overline{T}|}=2^{i}\cdot 3\cdot 5$ for some $i\le 10$.
By \cite[Theorem 5.3.9 and Proposition 5.4.13]{KL-book}, we have $|\Rad(G)|_s<s^{R_s}$ for $s\in \{2,3,5\}$, and then $G$ is a perfect central extension of $\POmega^+_{12}(2)$ by Lemma \ref{Schur}.
Thus $|\Rad(G)|\le (4,2^6-1)=1$ by \cite[pp. 173, Theorem 5.1.4]{KL-book}, a contradiction.

Finally, if $(\overline{G},\overline{T})$ is one of $(\POmega^+_8(3),\POmega^+_8(2))$, $(\PSL_{4}(4),\Sp_{4}(4))$ and $(\PSL_{6}(2),\PSU_{3}(3))$ then a similar contradiction arises as above, we omit the details here. This completes the proof of Theorem \ref{mainth-2}.

\end{document}